\documentclass[11pt]{amsart}
\usepackage{amssymb}
\usepackage{graphicx}
\usepackage{pinlabel}
\newtheorem{thm} {Theorem}
\newtheorem{claim}{Claim}
\newtheorem{conj} {Conjecture}
\newtheorem{cor} {Corollary}
\newtheorem{prop} {Proposition}
\newtheorem{lem} {Lemma}
\newtheorem{defn} {Definition}

\newtheorem{ex}{Example}

\title{Shape-Wilf-Ordering on Permutations of Length 3} 

\author{Zvezdelina Stankova}
\address{Zvezdelina Stankova, Dept. of Mathematics and
Computer Science, Mills College, Oakland, CA, {\tt stankova@mills.edu}}

\begin{document}
\thispagestyle{empty}
\maketitle
\vspace*{-5mm}
\centerline{August 2006}

\begin{abstract} The research on pattern-avoidance
has yielded so far limited knowledge on Wilf-ordering of
permutations. The Stanley-Wilf limits $\lim_{n\rightarrow \infty}
\sqrt[n]{|S_n(\tau)|}$ and further works suggest asymptotic ordering
of layered versus monotone patterns. Yet, B\'{o}na has provided
essentially the only known up to now result of its type on ordering of
permutations: $|S_n(1342)|<|S_n(1234)|<|S_n(1324)|$ for $n\geq 7$.  We
give a different proof of this result by ordering $S_3$ up to the
stronger {\it shape-Wilf-order}: $|S_Y(213)|\leq |S_Y(123)|\leq
|S_Y(312)|$ for any Young diagram $Y$, derive as a consequence that
$|S_Y(k+2,k+1,k+3,\tau)|\leq |S_Y(k+1,k+2,k+3,\tau)|\leq
|S_Y(k+3,k+1,k+2,\tau)|$ for any $\tau\in S_k$, and find out when
equalities are obtained. (In particular, for specific $Y$'s we find
out that $|S_Y(123)|=|S_Y(312)|$ coincide with every other Fibonacci
term.) This strengthens and generalizes B\'{o}na's result to arbitrary
length permutations. While all length-3 permutations have been shown
in numerous ways to be Wilf-equivalent, the current paper {\it
distinguishes} between and orders these permutations by employing all
{\it Young diagrams}. This opens up the question of whether
shape-Wilf-ordering of permutations, or some generalization of it, is
not the ``true'' way of approaching pattern-avoidance ordering.
\end{abstract}


\section{Introduction}

We review first basic concepts and results that are crucial to the
present paper, and direct the reader to \cite{SS,St1,We0,We1} for
further introductory definitions and examples on pattern-avoidance.

\smallskip A permutation $\tau$ of length $k$ is written as
$(a_1,a_2,\ldots,a_k)$ where $\tau(i)=a_i,\,1\leq i \leq k$. For
$k<10$ we suppress the commas without causing confusion. As usual,
$S_n$ denotes the symmetric group on $[n]=\{1,2, ...,n\}$.

\begin{defn}
{\rm Let $\tau$ and $\pi$ be two permutations of lengths $k$ and $n$,
respectively. We say that $\pi$ is $\tau$-{\it avoiding} if there is
no subsequence $i_{\tau(1)}, i_{\tau(2)}, ..., i_{\tau(k)}$ of $[n]$
such that $\pi(i_1)<\pi(i_2)<\ldots<\pi(i_k)$.  If there is such a
subsequence, we say that it is {\it of type} $\tau$, and denote this
by $\big(\pi(i_{\tau(1)}),$ $\pi(i_{\tau(2)})$,...,
$\pi(i_{\tau(k)})\big)\approx\tau$.}
\end {defn}

The following reformulation in terms of matrices is probably more
insightful. In it, and throughout the paper, we coordinatize all
matrices from {\it the bottom left corner} in order to keep the
resemblance with the ``shape'' of permutations.

\smallskip
\begin{defn}
{\rm Let $\pi \in S_n$. The {\it permutation matrix} $M(\pi)$ is the
$n\times n$ matrix $M_n$ having a $1$ in position $(i,\pi(i))$ for
$1\leq i \leq n$. Given two permutation matrices $M$ and $N$, we say
that $M$ {\it avoids} $N$ if no submatrix of $M$ is identical to $N$.}
\end{defn}

\noindent A permutation matrix is simply an arrangement, called a {\it
transversal}, of $n$ non-attacking rooks on an $n\times n$ board. We
refer to the elements of a transversal also as ``1's'' and
``dots''. Clearly, a permutation $\pi \in S_n$ contains a subsequence
$\tau \in S_k$ if and only if $M(\pi)$ contains $M(\tau)$ as a
submatrix.

\begin{defn} {\rm Let $S_n(\tau)$ denote the set of $\tau$-avoiding 
permutations in $S_n$.  Two permutations $\tau$ and $\sigma$ are {\it
Wilf-equivalent}, denoted by $\tau\sim\sigma$, if they are equally
restrictive: $|S_n(\tau)|=|S_n(\sigma)|$ for all $n \in \mathbb N$. If
$|S_n(\tau)|\leq |S_n(\sigma)|$ for all $n\in \mathbb N$, we say that
$\tau$ is {\it more restrictive} than $\sigma$, and denote this by
$\tau \preceq \sigma$.}
\end{defn} 

The classification of permutations in $S_k$ for $k\geq 7$ up to
Wilf-equivalence was completed over the last two decades by a number
of people. We refer the reader to Simion-Schmidt \cite{SS}, Rotem
\cite{Ro}, Richards \cite{Ri}, and Knuth \cite{Kn1,Kn2} for length
$k=3$; to West \cite{We1} and Stankova \cite{St1,St2} for $k=4$; to
Babson-West \cite{BW} for $k=5$; and to Backelin-West-Xin \cite{BWX}
and Stankova-West \cite{Stankova-West} for $k=6,7$. 

\medskip
However, total Wilf-{\it ordering} does not exist for a
general $S_k$. The first counterexample occurs in $S_5$
(cf. \cite{Stankova-West}): if $\tau=(53241)$ and $\sigma=(43251)$,
then $S_7(\tau)<S_7(\sigma)$ but $S_{13}(\tau)>S_{13}(\sigma)$, and
hence $\tau$ and $\sigma$ cannot be Wilf-ordered. This phenomenon
prompts 

\begin{defn} {\rm For two permutations $\tau$ and $\sigma$, we say
that $\tau$ is {\it asymptotically more restrictive} than
$\sigma$, denoted by $\tau \preceq_a \sigma$, if $|S_n(\tau)|\leq
|S_n(\sigma)|$ for all $n\gg 1$.}
\end{defn} 

Stanley-Wilf Theorem (cf. Marcus and Tardos \cite{Marcus}, Arratia
\cite{Arratia}) gives some insight into the asymptotic ordering of
permutations. Inequalities between the Stanley-Wilf limits
$L(\tau)=\lim_{n\rightarrow \infty} \sqrt[n]{|S_n(\tau)|}$ suggest
asymptotic comparisons between the corresponding permutations. For
instance, works of B\'{o}na \cite{Bon2,Bon4} and Regev \cite{Regev} show
that $L(I_k)=(k-1)^2\leq L(\tau)$, where $I_k=(12...k)$ is the {\it
identity} pattern and $\tau$ is any {\it layered} pattern in $S_k$
(cf. Definition~\ref{decomposable}), which provides strong evidence
that the identity pattern is more restrictive than all layered
patterns in $S_k$. Yet, this result will still not imply asymptotic
ordering between the above types of patterns if it happens that
$L(I_k)=L(\tau)$ for some layered $\tau$.

In \cite{Bon1,Bon3}, B\'{o}na provides essentially the only known so
far result on Wilf-ordering:
\begin{equation}
|S_n(1342)|<|S_n(1234)|<|S_n(1324)|\,\,\text{for}\,\,n\geq 7,
\label{WOS4}
\end{equation}
along with some sporadic examples on asymptotic Wilf-ordering,
e.g. $I_k\preceq_a \tau_k$ for some $\tau_k\in S_k$.  Since $S_2$ and
$S_3$ are each a single Wilf-equivalence class (cf.~\cite{SS}), the
first possibility of nontrivial Wilf-ordering arises in $S_4$. A
representative of each of the 3 Wilf-equivalence classes in $S_4$
appears in (\ref{WOS4}) (cf. \cite{We1,St1,St2}.).

\medskip
In order to prove differently and extend result (\ref{WOS4}) to
Wilf-ordering of certain permutations of arbitrary lengths, we shall
use the concept of a stronger Wilf-equivalence relation, called {\it
shape-Wilf-equivalence}. The latter was introduced in \cite{BW}, and
further explored in consequent papers \cite{BWX,Stankova-West}.

\begin{defn}{\rm  A {\it transversal} $T$ of a Young
diagram $Y$, denoted $T\in S_Y$, is an arrangement of 1's such that
every row and every column of $Y$ has exactly one 1 in it. A subset of
1's in $T$ forms a {\it submatrix} of $Y$ if all columns and rows of
$Y$ passing through these 1's intersect {\it inside} $Y$. For a
permutation $\tau\in S_k$, $T$ {\it contains the pattern} $\tau$ (in
$Y$) if some $k$ 1's of $T$ form a submatrix of $Y$ identical to
$M(\tau)$. Denote by $S_Y(\tau)$ the set of all transversals of $Y$
which avoid $\tau$.}
\end{defn}

\begin{figure}[h]
\begin{center}
\includegraphics[width=1.6in]{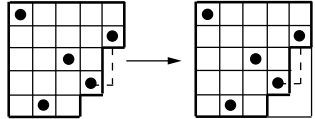}
\caption{$T\in S_Y$ versus $T^{\prime}\in S_5$}
\label{T in SY}
\end{center} 
\end{figure}

Now, suppose $T\in S_Y$ has a {\it subsequence}
$\mathcal{L}=(\alpha_1\alpha_2...\alpha_k)\approx \tau\in S_k$. From
the above definition, in order for $T\in S_Y$ to {\it contain the
pattern} $\tau$ in $Y$, it is necessary and sufficient that the column
of the rightmost element of $\mathcal{L}$ and the row of the smallest
element of $\mathcal{L}$ intersect inside $Y$. In such a case, we say
that the subsequence $\mathcal{L}$ {\it lands} inside $Y$. For
example, Figure \ref{T in SY}a shows the transversal $T\in S_Y$
representing the permutation $(51324)$. Note that $T$ contains the
patterns $(312)$ and $(321)$ because its subsequences $(513)$ and
$(532)$ land inside $Y$. However, $T$'s subsequence
$(324)\approx(213)$ does not land in $Y$, and in fact, $T$ does not
contain the {\it pattern} $(213)$; symbolically, $T\in S_Y(213)$.

When $Y$ is a square diagram of size $n$, $S_n(\tau)\equiv
S_Y(\tau)$. Let $Y(a_1,a_2,...,a_n)$ denote the Young diagram $Y$
whose $i$-th row has $a_i$ cells, for $1\leq i\leq n$.  In order for 
$Y$ to have any transversals at all, it must be {\it proper}: $Y$ must
have the same number of rows and columns and must contain the {\it
staircase} diagram $St^1=Y(n,n-1,...,2,1)$; equivalently, $Y$ must
contain its southwest-northeast $45^{\circ}$ {\it diagonal} $d(Y)$
which connects $Y$'s bottom left and top right corners. If
not specified otherwise, a Young diagram is always proper in this
paper.

Young diagrams are traditionally coordinatized from {\it the top left
corner}, meaning that their first (and largest) row and column are the
top, respectively, leftmost ones. To avoid possible confusion with the
matrix ``bottom-left-corner'' coordinatization used in this paper, one
can think of a transversal $T\in S_Y$ by first completing the (proper)
Young diagram $Y$ to a square matrix $M_n$, and then taking a
transversal $T$ of $M_n$ all of whose 1's are in the original cells of
$Y$.  Thus, whether using a matrix or a Young diagram, all
transversals resemble the ``shape'' of permutations. For instance, in
Fig.~\ref{T in SY}, the proper Young diagram $Y(5,5,4,4,3)$ is
completed to the square matrix $M_5$, and the transversal $T\in S_Y$
induces a transversal $T^{\prime}\in S_5$. As observed above, $T\in
S_Y(213)$, but $T^{\prime}\not\in S_5(213)$ because the subsequence
$(324)\approx (213)$ of $T^{\prime}$ does land in $M_5$.

\begin{defn}{\rm Two permutations $\tau$ and $\sigma$ are called 
{\it shape-Wilf-equivalent (SWE)}, denoted by 
$\tau\sim_s\sigma$, if $|S_{Y}(\tau)|=|S_{Y}(\sigma)|$ for
all Young diagrams $Y$. If $|S_{Y}(\tau)|\leq |S_{Y}(\sigma)|$ for all
such $Y$, we say that $\tau$ is {\it more shape-restrictive} than
$\sigma$, and denote this by $\tau\preceq_s\sigma$. }
\end{defn}
Clearly, $\tau\sim_s\sigma$ ($\tau\preceq_s\sigma$) imply
$\tau\sim\sigma$ ($\tau\preceq\sigma$, respectively), but the
converses are false.  Babson-West showed in \cite{BW} that SWE is
useful in establishing more Wilf-equivalences. To the best of our
knowledge, this idea of Young diagrams has not been yet been modified
or used to prove Wilf-{\it ordering}, which the present paper will
accomplish. To this end, we include below an extension of
Babson-West's proposition, replacing shape-Wilf-{\it equivalences}
``$\sim_s$'' with shape-Wilf-{\it ordering} ``$\preceq_s$''. Section
\ref{Proof of Proposition 1} presents a modification and extension of
their original proof, and introduces along the way new notation
necessary for the completion of our Wilf-ordering results.

\begin{prop} Let $A\preceq_sB$ for some permutation matrices 
$A$ and $B$. Then for any permutation matrix $C$:
\[\left( \begin{array}{c|c}
          A & 0\\ \hline
          0 & C
      \end{array} \right) \preceq_s
\left(\begin{array}{c|c}
          B & 0\\ \hline
          0 & C
      \end{array} \right)\cdot\]
\label{prop-modified-BW}
\end{prop}

\noindent If we shape-Wilf-order permutations in $S_k$ for a small
$k$, Proposition~\ref{prop-modified-BW} will enable us to
shape-Wilf-order some permutations in $S_n$ for larger $n$.  Since
$(12)\sim_s(21)$ in $S_2$, Proposition~\ref{prop-modified-BW} can
imply in this case only shape-Wilf-{\it equivalences}.

The first non-trivial shape-Wilf-{\it ordering} can occur in $S_3$,
since the latter splits into three distinct shape-Wilf-equivalence
classes: $\{(213)\sim_s (132)\}$, $\{(123)\sim_s(231)\sim_s(321)\}$,
and $\{(312)\}$. The first SWE-class was proven by Stankova-West in
\cite{Stankova-West}, and the second class was proven by
Babson-Backelin-West-Xin in \cite{BW,BWX}. The smallest Young diagram
for which all three classes differ from each other is
$Y=Y(5,5,5,5,4)$: $|S_Y(213)|=37<|S_Y(123)|=41<|S_Y(312)|=42$.
Numerical evidence suggests that such inequalities hold for {\it all}
Young diagrams $Y$, and indeed this is true:

\begin{thm}[Main Theorem] For all Young diagrams $Y$:
\label{SWOS3}
\[|S_Y(213)|\leq |S_Y(123)|\leq |S_Y(312)|.\]
\end{thm}

\begin{figure}[h]
\labellist
\small\hair 2pt
\pinlabel $\tau$ at 119 644
\pinlabel $\tau$ at 237 644
\pinlabel $\tau$ at 355 644
\pinlabel $<$ at 157 662
\pinlabel $<$ at 275 662
\endlabellist
\begin{center}
\includegraphics[width=2.5in]{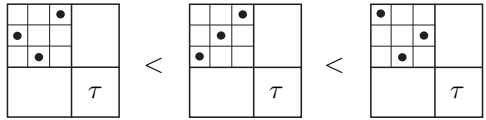}
\caption{Corollary \ref{Wilf-ordering-corollary}}
\label{Wilf-ordering corollary}
\end{center}
\end{figure}

Figure~\ref{Wilf-ordering corollary} with $\tau=\emptyset$
illustrates Theorem \ref{SWOS3}. Let $Y_n=Y(n,n,n,...,n,n-1)$ be the
Young diagram obtained by removing the right bottom cell from the
square $M_n$. Section~\ref{difference} shows 
\[|S_{Y_n}(213)|< |S_{Y_n}(123)| <
|S_{Y_n}(312)|\,\,\text{for}\,\,n\geq 5.\] These strict inequalities
preclude the possibility of the three permutations $(213)$, $(123)$,
$(312)$ to be asymptotically SWE, even though they are
Wilf-equivalent.  More precisely,
\begin{thm}
$|S_Y(213)|<|S_Y(123)|$ if and only if $Y$ contains an $i$-critical
 point with $i\geq 2$, and $|S_Y(123)|<|S_Y(312)|$ if and only if $Y$
 contains an $i$-critical point with $i\geq 3$.
\label{summary 312>321>213}
\end{thm}
The definition and a discussion of {\it critical points} can be found
in Subsection~\ref{Diagonal Properties and Critical Points}.  While
for any $\tau\in S_3$ the ``Wilf-numbers'' $|S_n(\tau)|$ equal the
Catalan numbers $c_n=\frac{1}{n+1}\binom{2n}{n}$, the
``shape-Wilf-numbers'' $|S_Y(\tau)|$ naturally vary a lot more. In
particular, for the staircases $Y=St^3_n$, $|S_Y(\tau)|$ coincide with
the odd-indexed Fibonacci terms $f_{2n-1}$, and hence involve the {\it
golden ratio} $\phi=(1+\sqrt{5})/2$ (cf. Definition~\ref{diagonals}
and Section~\ref{difference}.)

\begin{defn} {\rm We say that a permutation $\tau\in S_n$ is
{\it decomposable} into blocks $A_1$ and $A_2$ if for some $k<n$,
$\tau$ can be partitioned into two subpatterns
$A_1=(\tau_1,\tau_2,...,\tau_k)$ and
$A_2=(\tau_{k+1},\tau_{k+2},...,\tau_n)$ such that all entries of
$A_1$ are {\it bigger} than (and a priori come before) all entries of
$A_2$. We denote this by $\tau=(A_1|A_2)$.  If there is no such
decomposition into two blocks, we say that $\tau$ is {\it
indecomposable}. In particular, a {\it layered} pattern $\tau$ is a
permutation decomposable into increasing blocks. }
\label{decomposable}
\end{defn}

\noindent For example, $(4132)=(4|132)$ is decomposable, while
$(3142)$ and $(1432)$ are indecomposable; $(4123)=(4|123)$ is layered,
while $(4132)$ is not layered. Without confusion, we can also write
$(213|1)$ instead of $(3241)$.  In this notation, Proposition
\ref{prop-modified-BW} can be rewritten as $A\preceq_s
B\,\,\Rightarrow\,\,(A|C)\preceq_s(B|C)$.

\begin{cor}
For any permutation $\tau\in S_k$, $(213|\tau)\preceq_s (123|\tau)
\preceq_s (312|\tau)$. Moreover, strict asymptotic Wilf-ordering
$|S_n(213|\tau)|< |S_n(123|\tau)| <
|S_n(312|\tau)|$ occurs for $n\geq 2k+5.$
\label{Wilf-ordering-corollary}
\end{cor}

\noindent In particular, when $\tau=(1)$ Corollary
\ref{Wilf-ordering-corollary} reduces to:
\begin{eqnarray*}& &|S_n(213|1)|< |S_n(123|1)| <
|S_n(312|1)|\,\,\text{for}\,\,n\geq 7\\
&\Rightarrow& \,|S_n(3241)|\,< |S_n(2341)| \,<
\,|S_n(4231)|\,\,\,\text{for}\,\,n\geq 7.
\end{eqnarray*}

\begin{figure}[h]
\begin{center}
\includegraphics[width=4in]{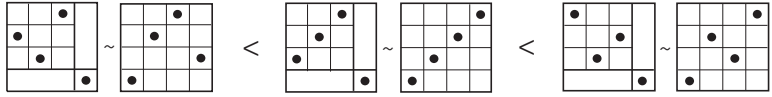}
\caption{Wilf-Ordering of $S_4$}
\label{Label SWOS4}
\end{center}
\end{figure}

Note that $(3241)\sim (1342)$ and $(4231)\sim (1324)$
(cf. Fig.~\ref{Label SWOS4}a-c) since the two permutation matrices in
each Wilf-equivalence pair can be obtained from each other by applying
symmetry operations of flipping along vertical, horizontal and/or
diagonal axes (cf. \cite{We1,St1}). Further, $(2341)\sim (1234)$ by
the SWE-relations in \cite{BWX}, or by an earlier work
\cite{St2}. Thus, choosing the second representatives of the three
Wilf-equivalence classes in $S_3$, we obtain B\'{o}na's (\ref{WOS4})
inequality as a special case of Corollary
\ref{Wilf-ordering-corollary}.

Some of the implied new shape-Wilf-orderings by
Corollary~\ref{Wilf-ordering-corollary} in $S_5$ and $S_6$ are:
\[\begin{array}{cccccccccc}
(43521)&\!\!\!\prec_s^*\!\!\!& (54321)&\!\!\!\prec_s\!\!\!& (53421)&
\quad\quad
(546231)&\!\!\!\prec_s^*\!\!\!&(654231)&\!\!\!\prec_s^*\!\!\!& (645231),\\
(546321)&\!\!\!\prec_s^*\!\!\!&(654321)&\!\!\!\prec_s\!\!\!& (645321)
& \quad\quad
(546213)&\!\!\!\prec_s^*\!\!\!& (654213)&\!\!\!\prec_s^*\!\!\!& (645213).
\end{array}\]
These inequalities imply Wilf-orderings, of which the ones
corresponding to $*$'s are new.
The paper is organized as follows.  Section~\ref{Proof of Proposition
1} presents the proof of Proposition~\ref{prop-modified-BW}, along
with a strategy for establishing strict asymptotic Wilf-orderings. In
Section~\ref{general set up}, we introduce {\it critical points},
provide the $0$- and $1$-{\it splittings} $S_Y(\sigma)\cong
S_{Y^R}(\sigma)\times S_{\phantom{}_QY}(\sigma)$ in
Proposition~\ref{0-1-splitting}, and a {\it $2$-critical splitting} in
Lemma~\ref{General 2-splitting lemma}. Subsection~\ref{general moves}
defines the {\it $\sigma\!\rightarrow\!\tau$ moves} on transversals in
$Y$, and opens up the discussion of the induced maps
$\phi:S_Y(\tau)\rightarrow S_Y(\sigma)$.  Sections~\ref{SY(321) <
SY(312)}-\ref{SY(213) < SY(123)} contain the proof of the inequalities
$|S_Y(312)|\geq |S_Y(321)|$ and $|S_Y(213)|\leq |S_Y(123)|$; a
description of the structures of $T\in S_Y(321)$ and $T\in S_Y(312)$
can be found in Subsections~\ref{structure of S_Y(321)}-\ref{structure
of S_Y(312)}. Using critical points, necessary and sufficient
conditions for strict inequalities $|S_Y(312)|>|S_Y(321)|$ and
$|S_Y(213)|<|S_Y(123)|$ are established in Sections~\ref{strict
inequalities 312>321}-\ref{strict
|SY(213)|<|SY(123)|}. Section~\ref{strict Wilf-ordering} provides the
proof of the strict Wilf-orderings $|S_n(213|\tau)|< |S_n(123|\tau)| <
|S_n(312|\tau)|$ for $n\geq 2k+5$. Finally, in
Section~\ref{difference} we calculate $|S_Y(\tau)|$ for $\tau\in S_3$
and Young diagrams $Y$ which are extreme with respect to their
critical points. The paper ends with a generalization of the
Stanley-Wilf limits and the fact that $\phi^2$ is such a limit.

\section{Proof of Proposition \ref{prop-modified-BW}}
\label{Proof of Proposition 1}

In this section we present a modified and extended version of the
original proof of Babson-West to address our new setting of shape-Wilf
ordering. Let the permutation matrices $A$, $B$ and $C$ represent
permutations $\alpha$, $\beta$ and $\gamma$, respectively. Before we
proceed with the proof, we need to introduce some definitions and
notation.

\subsection{Various subboards of $Y$}
\begin{figure}[h]
\labellist
\small\hair 2pt
\pinlabel $c$ at 134 565
\pinlabel $Y_c$ at 92 600
\pinlabel $c$ at 359 566
\pinlabel $Y_{\bar{c}}$ at 307 610
\pinlabel $Y^{\bar{c}}$ at 307 518
\pinlabel $\phantom{}^cY$ at 161 546
\pinlabel $\phantom{}^{\bar{c}}Y$ at 395 540
\pinlabel $\phantom{}_cY$ at 402 610
\endlabellist
\centering
\includegraphics[width=2.5in]{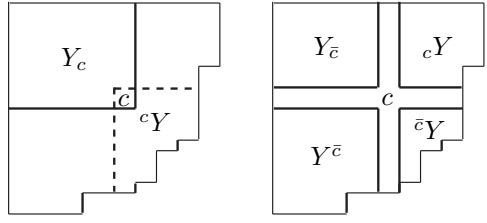}
\caption{Notation $Y_c\,\,\text{and}^{\,\,c}Y$ versus
$Y_{\bar{c}},^{\,\,\bar{c}}Y$, etc.}
\label{Notations Yc}
\end{figure}

Let $Y$ be a Young diagram, and let $c$ be a cell in $Y$. Denote by
$\phantom{}^{\bar{c}}Y$ the subboard of $Y$ to the right and below
$c$, not including $c$'s row and column; and by $Y_c$ the subboard of
$Y$ to the left and above $c$, including the corresponding cells in
$c$'s row and column. Since $Y$ is a Young diagram,
$\phantom{}^{\bar{c}}Y$ is also a Young diagram (not necessarily
proper), and $Y_c$ is a rectangle whose right bottom cell is $c$
(cf. Fig.~\ref{Notations Yc}). 
This notation is created so as to match the relative positions of $c$
and the corresponding subboard of $Y$, where exclusion of $c$'s row and
column is denoted by $\bar{c}$. In the same vein, we define $Y^c$,
$\phantom{}_{\bar{c}}Y$, etc.  We also extend the notation to (full or
partial) {\it transversals} $T$ of $Y$, to elements $\alpha\in T$, and
to grid points $P$ of $Y$; for instance,
$\phantom{}_{\bar{\alpha}}T=T|_{\phantom{}_{\bar{\alpha}}Y}$ is the
restriction of $T$ onto the subboard $\phantom{}_{\bar{\alpha}}Y$,
while $Y^P$ is the subboard $Y^c$ where $P$ is the top right corner of
cell $c$.

We use the symbols $\nearrow$ and $\searrow$ instead of the words
``increasing'' and ``decreasing''. Thus, $I_k\!\!\nearrow$, and its
transpose $J_k\!\!\searrow$.

\begin{defn}
{\rm Let $T\in S_Y$, and $\alpha,\beta\in T$. We say that $\alpha$
{\it $(21)$-dominates} $\beta$ if
$(\alpha\beta)\!\!\searrow$. Similarly, $\alpha$ {\it
$(12)$-dominates} $\beta$ if $(\beta\alpha)\!\!\nearrow$ and lands in
$Y$. We extend these definitions to any cells of and dots in $Y$.}
\end{defn}

\subsection{Coloring of $Y$ with respect to $T$ and $\gamma$} 
Fix a transversal $T\in S_Y$. With respect to the pattern $\gamma$,
$T$ induces a white/blue coloring on $Y$'s cells as follows. Color a
cell $c$ in $Y$ {\it white} if $\phantom{}^{\bar{c}}Y$ contains $C$ as
a submatrix; otherwise, color $c$ {\it blue}. Clearly, for every white
cell $w$, the rectangle $Y_w$ is also entirely white. Hence, the white
{\it subboard} $W^{\prime}$ of $Y$ is a Young subdiagram of $Y$ (not
necessarily proper), and $T$ induces a {\it partial} transversal
$T|_{W^{\prime}}$ of $W^{\prime}$.

In order for $T$ to avoid $(\alpha|\gamma)$, it is necessary and
sufficient that $T|_{W^{\prime}}$ avoids $\alpha$. However, some rows
and columns of $W^{\prime}$ cannot participate in any undesirable
$\alpha$-patterns since the 1's in them are in blue cells: recolor
these white rows and columns of $W^{\prime}$ to blue. After deletion
of the newly blue rows and columns of $W^{\prime}$, the latter is
reduced to a white proper Young {\it subdiagram} $W$ of $Y$, while
$T|_{W^{\prime}}$ is reduced to a {\it full} transversal $T|_W$ of
$W$.

\begin{defn} {\rm We say that the transversal $T$ of $Y$ {\it induces
with respect to $\gamma$} the white subdiagram $W$ of $Y$ and the
(full) transversal $T|_W$ of $W$. Let $S^W_Y(\alpha|\gamma)$ denote
the set all transversals $T\in S_Y(\alpha|\gamma)$ which induce $W$
with respect to $\gamma$. }
\end{defn}
For example, Figure~\ref{Splitting of T2}a shows a transversal $T\in
S_Y$ and the induced white subboard $W^{\prime}$ with respect to
$\gamma=(213)$: the blue subboard of $Y$ is depicted with its grid
lines, while $W^{\prime}$ is depicted without them; the dashed lines
pass through some of the blue 1's and indicate that these rows and
columns of $Y$ will be deleted from
$W^{\prime}$. Figure~\ref{Splitting of T2}c shows the final white
subdiagram $W(4,4,3,3)$ and its transversal $T|_W=(2134)$. 
Figure~\ref{Splitting of T2}a-c also illustrates that
$T=(7,6,9,2,10,1,4,5,3,8)\in S_Y$ avoids $(123|213)$ because
$T|_W=(2134)$ avoids $(123)$ on $W$, but it contains $(213|213)$
because $T|_W$ contains the pattern $(213)$ on $W$.

\smallskip
We summarize the observations in this subsection in the following
\begin{lem} Let $W$ be any Young subdiagram of $Y$. Then
\begin{enumerate}
\item
$T\in S_Y(\alpha|\gamma)\,\,\Leftrightarrow\,\,T|_W\in S_W(\alpha)$.

\smallskip
\item
${S_Y(\alpha|\gamma)=\bigsqcup_{W\subset Y}S^W_Y(\alpha|\gamma)}$.
\end{enumerate}
\end{lem}

\subsection{Splitting of transversals $T\in S_Y$ with respect to $\gamma$} 
\label{key observation}
Fix now a (white) Young subdiagram $W$ of $Y$, and let $T\in
S^W_Y(\alpha|\gamma)$. By construction of $W$, $T$ splits itself into
two disjoint subsets: the induced transversal $T|_W$ of $W$ consisting
of all ``white'' 1's, and the remainder $T_\gamma=T\backslash{T|_W}$
consisting of all ``blue'' 1's. We denote this by \[T=T|_W\oplus
T_\gamma,\,\,\text{where}\,\,T|_W\in S_W(\alpha).\]

\begin{figure}[h]
\labellist
\small\hair 2pt
\pinlabel $W^{\prime}$ at 147 633
\pinlabel $W$ at 323 633
\pinlabel $W$ at 323 525
\pinlabel $W^{\prime}$ at 456 633
\pinlabel $W^{\prime}$ at -90 633
\endlabellist
\begin{center}
\includegraphics[width=5.5in]{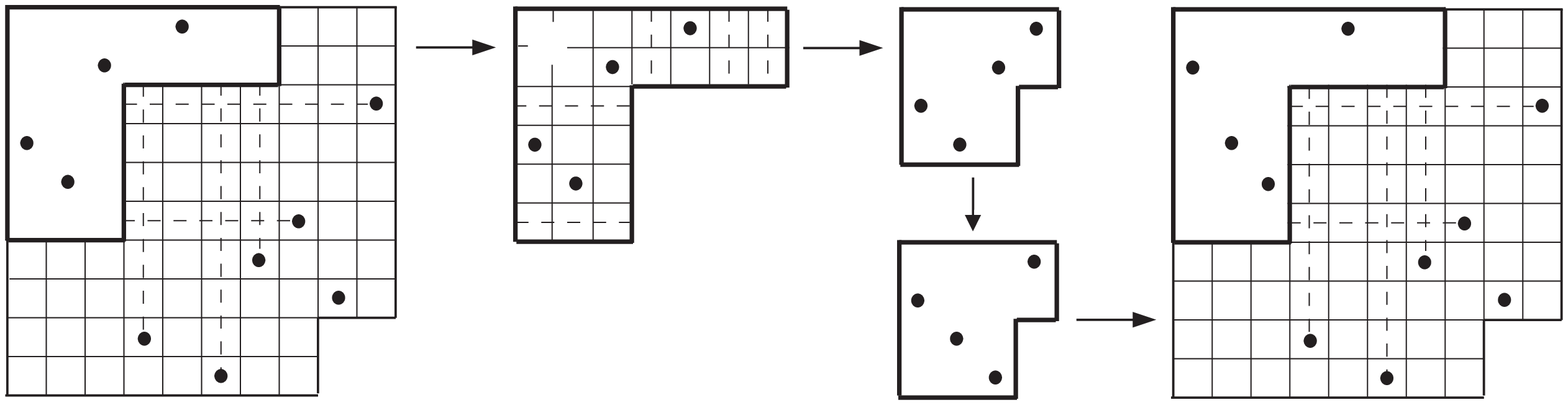}
\caption{$T=T|_W\oplus T_{(213)} \rightarrow
T^{\prime}=T^{\prime}_W\oplus T_{(213)}$ in $S_Y^W(123|213)$}
\label{Splitting of T2}
\end{center} 
\end{figure}

A {\it key observation} is that, if $T^{\prime}_W$ is another
transversals in $S_W(\alpha)$, then $T^{\prime}=T^{\prime}_W \oplus
T_\gamma\in S^W_Y(\alpha|\gamma)$. This is true because fixing
$T_\gamma$ preserves the white cells of $W$, and replacing $T|_W$ with
any other transversal of $W$ certainly does not affect the blue
colored cells in $Y\backslash W$. For example, Figure~\ref{Splitting
of T2} shows $T\in S_Y^W(123|213)$ with $W=(4,4,3,3)$, $T|_W=(2134)\in
S_W(123)$, and $T_{(213)}=(214538)\approx (214536)$. If we keep
$T_{(213)}$ and replace $T|_W$ with another $T^{\prime}_W=(3214)\in
S_W(123)$ (shown in Fig.~\ref{Splitting of T2}d), we obtain the
transversal in Fig.~\ref{Splitting of T2}e:
\[T^{\prime}=(9,7,6,2,10,1,4,5,3,8)=(9,7,6,10)\oplus
(2,1,4,5,3,8) \in S_Y^W(123|213).\] 

\noindent We conclude that all transversals $T\in
S^W_Y(\alpha|\gamma)$ whose second component is a fixed $T_\gamma$ are
obtained by adding an arbitrary transversal $T^{\prime}_W\in
S_W(\alpha)$ to $T_\gamma$:
\[T=T^{\prime}_W\oplus T_\gamma\in S^W_Y(\alpha|\gamma)\,\,\text{for any}\,\,
T^{\prime}_W\in S_W(\alpha).\]

\subsection{Description of 
the $T_\gamma$-component of $T\in S^W_Y(\alpha|\gamma)$} We can extend
the definitions of the white/blue coloring of $Y$ above to {\it
partial} transversals $T^{\prime}$ of $Y$: a {\it blue} cell $b$ in
$Y$ is such that $\phantom{}^{\bar{b}}Y$ does {\bf not} contain a
$\gamma$-subpattern of $T^{\prime}$, while a {\it white} cell $w$ in
$Y$ is such that $\phantom{}^{\bar{w}}Y$ does contain a
$\gamma$-subpattern of $T^{\prime}$.

Recall the notion of {\it reduction of $Y$ along a subset $X$} of
$Y$'s cells, introduced in \cite{Stankova-West}:
$Y\!\!\big/_{\displaystyle{\!\!X}}$ is the Young subdiagram obtained
from $Y$ by deleting all rows and columns of $Y$ which intersect $X$.
This notation should not be confused with $Y\backslash X$ - the
subboard obtained from $Y$ by removing the cells in $X$, or with
$T|_W$ - the restriction of $T$ on $W$.

\begin{figure}[h]
\labellist
\small\hair 2pt
\pinlabel $T^{\prime}$ at 127 566
\pinlabel $W$ at 46 629
\pinlabel $W$ at 289 629
\endlabellist
\begin{center}
\includegraphics[width=2.5in]{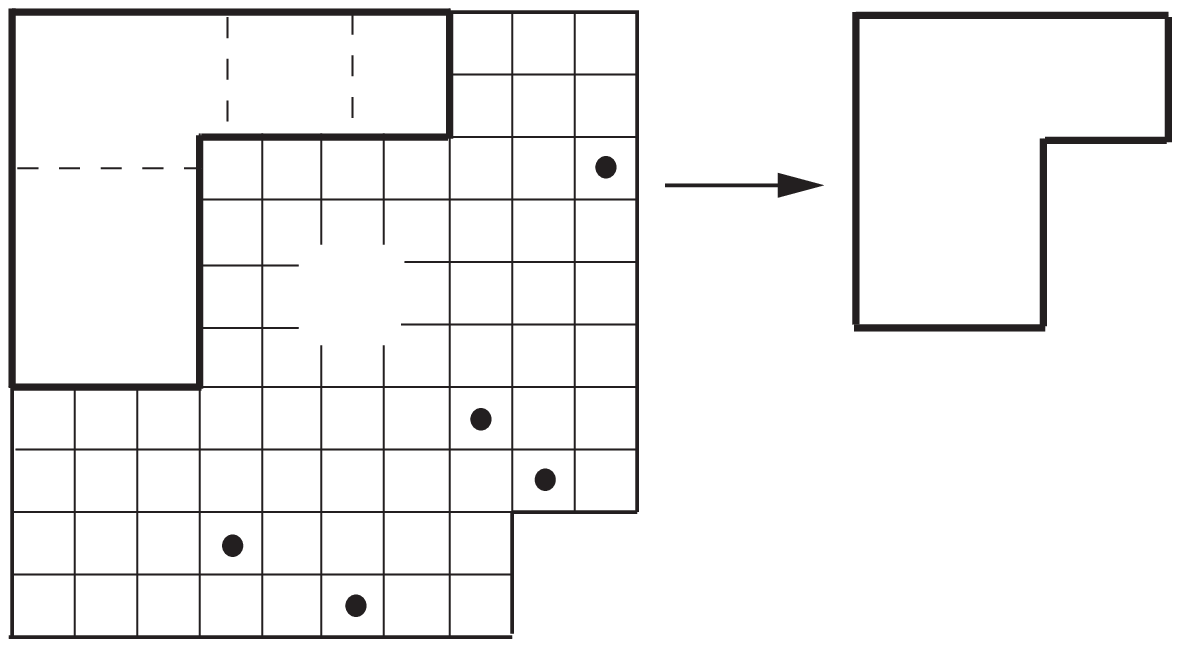}
\caption{$T^{\prime}$ saturates $W$ with respect to $(213)$}
\label{Saturation 1}
\end{center} 
\end{figure}

\begin{defn}
\label{def-saturation}
{\rm Let $W$ be a proper subdiagram of a 
Young diagram $Y$. A partial transversal $T^{\prime}$ of $Y$ {\it
saturates} $W$ {\it with respect to $\gamma$} if the induced by
$T^{\prime}$ blue/white coloring on $Y$ with respect to
$\gamma$ satisfies:

\begin{itemize}
\item[(1)] $T^{\prime}$'s elements are all placed in blue cells;
\item[(2)] Reducing $Y$ along $T^{\prime}$ and removing any leftover
blue cells results in $W$; and 
\item[(3)] $|W|+|T^{\prime}|=|Y|$, where $|U|$ is the size of a
  proper Young diagram $U$ and $|T^{\prime}|$ counts the number of
  elements in $T^{\prime}$.
\end{itemize}}
\end{defn}
\noindent Since a blue cell cannot $(21)$-dominate a white cell, no
matter which transversal of $W$ we choose to complete $T^{\prime}$ to
a (full) transversal of $Y$, the blue/white coloring of $Y$ will
remain the same (cf. Fig.~\ref{Saturation 1}.) Condition (3)
ensures that there is no entirely blue row or column without an
element of $T^{\prime}$; in fact, (3) matches the sizes of $W$ and
$T^{\prime}$ so that any transversal of $W$ will indeed complete
$T^{\prime}$ to a {\it full} transversal of $Y$.  

According to Definition~\ref{def-saturation}, for a transversal $T\in
S^W_Y(\alpha|\gamma)$ with splitting $T=T|_W\oplus T_\gamma$, the
partial transversal $T_\gamma$ of $Y$ saturates $W$ with respect to
$\gamma$.

\begin{defn}{\rm Given a subdiagram $W$ of the Young diagram $Y$, let
$\bar{S}_{Y\backslash W}(\gamma)$ denote the set of partial
transversals $T^{\prime}$ of $Y$ which saturate $W$ with respect to
$\gamma$.}
\end{defn}

\subsection{Splitting Formula for $|S_Y(\alpha|\gamma)|$}
\label{splitting formula subsection} We have seen that 
any transversal $T\in S^W_Y(\alpha|\gamma)$ splits uniquely as
$T=T|_W\oplus T_\gamma$, where $T|_W$ avoids $\alpha$ on $W$ and
$T_\gamma$ saturates $W$ in $Y$ with respect to $\gamma$. This defines
an injective map $S^W_Y(\alpha|\gamma)\hookrightarrow S_W(\alpha)
\times \bar{S}_{Y\backslash W}(\gamma)$. The key observation in
Subsection \ref{key observation} shows that this map is
surjective. Therefore,

\begin{lem}[Splitting Formula for $|S_Y(\alpha|\gamma)|$] 
For any subdiagram $W$ of the Young diagram $Y$, the 
isomorphism of sets $S^W_Y(\alpha|\gamma)\cong S_W(\alpha) \times
\bar{S}_{Y\backslash W}(\gamma)$ holds true.  Consequently,
\[|S_Y(\alpha|\gamma)|=
\sum_{W\subset Y}|S_W(\alpha)|\cdot |\bar{S}_{Y\backslash W}(\gamma)|,\]
where the sum is taken over all Young subdiagrams $W$ of $Y$.
\label{splitting lemma}
\end{lem}

Since the components $\bar{S}_{Y\backslash W}(\gamma)$ depend only on
$\gamma$ and $W$ (but not on $\alpha$), this allows for direct
comparisons between $S_Y(\alpha|\gamma)$ and $S_Y(\beta|\gamma)$. In
particular, if $\alpha\preceq_s \beta$, then $|S_W(\alpha)|\leq
|S_W(\beta)|$ for any Young diagram $W$, and the splitting formulas
for $\alpha$ and $\beta$ imply $|S_Y(\alpha|\gamma)|\leq
|S_Y(\beta|\gamma)|$.  This completes the Proof of Proposition
\ref{prop-modified-BW}. \qed

\subsection{Strategy for proving strict Wilf-ordering}
\label{strategy}
When $\alpha\preceq \beta$, the Splitting Formula can be used to prove
a strict asymptotic Wilf-ordering of the form
$|S_n(\alpha|\gamma)|\lneqq |S_n(\beta|\gamma)|$, provided that for
$n\gg 1$:

\begin{itemize}
\item[(SF1)] there is a Young diagram $W_n$ with $|S_{W_n}(\alpha)|\lneqq
  |S_{W_n}(\beta)|$; and
\item[(SF2)] there is a partial transversal $T_n$ of
$M_n$ saturating $W_n$ with respect to $\gamma$.
\end{itemize}
The existence of $W_n$ and $T_n$ ensures that $|S_{W_n}(\beta)|>0$ and
$|\bar{S}_{M_n\backslash W_n}(\gamma)|>0$, so that
\[|S_{W_n}(\alpha)|\cdot |\bar{S}_{M_n\backslash W_n}(\gamma)|\lneqq
|S_{W_n}(\beta)|\cdot |\bar{S}_{M_n\backslash W_n}(\gamma)|.\] We
shall employ this strategy in Section~\ref{strict Wilf-ordering} to
show strict asymptotic Wilf-ordering between the permutations
$(213|\tau),(123|\tau)$ and $(312|\tau)$ of Corollary
\ref{Wilf-ordering-corollary}.

\section{Critical Splittings of Young Diagrams and Transversals} 
\label{general set up}

\subsection{First and second subsequences of $T\in S_Y$.} 
Recall that $\alpha\in T$ is a {\it left-to-right maximum} of $T$ if
$\alpha$ is {\bf not} $(21)$-dominated by any other element of $T$, i.e.
$T_{\bar{\alpha}}=\emptyset$. 

\begin{defn}{\rm 
Let $T\in S_Y$. The subsequence $T^1$ of all left-to-right maxima
$\alpha_i$ of $T$ is called the {\it first subsequence} of $T$.  The
{\it second subsequence} $T^2$ of $T$ consists of all elements
$\beta_j\in T\backslash T^1$ for which $Y_{\bar{\beta_j}}$ contains
only elements of $T^1$, i.e. $\beta_j$ is $(21)$-dominated only by (a
non-empty set of) elements of $T^1$. 
}
\label{first and second subsequences}
\end{defn}

\begin{figure}[h]
\labellist
\small\hair 2pt
\pinlabel $T_{\bar{\alpha}}$ at 57 621
\pinlabel $T^1$ at 183 678
\pinlabel $T^2$ at 240 678
\pinlabel $\scriptstyle{\alpha}$ at 90 575
\pinlabel $\scriptstyle{\beta}$ at 126 555
\pinlabel $c$ at 386 567
\pinlabel $T^{\bar{c}}$ at 336 515
\pinlabel $\phantom{}_{\bar{c}}T$ at 433 613
\pinlabel $St^3_{10}$ at 610 593
\pinlabel $d_0(Y)$ at 761 652
\pinlabel $d_2(Y)$ at 761 615
\pinlabel $P$ at 656 460
\endlabellist
\begin{center}
\includegraphics[width=4.5in]{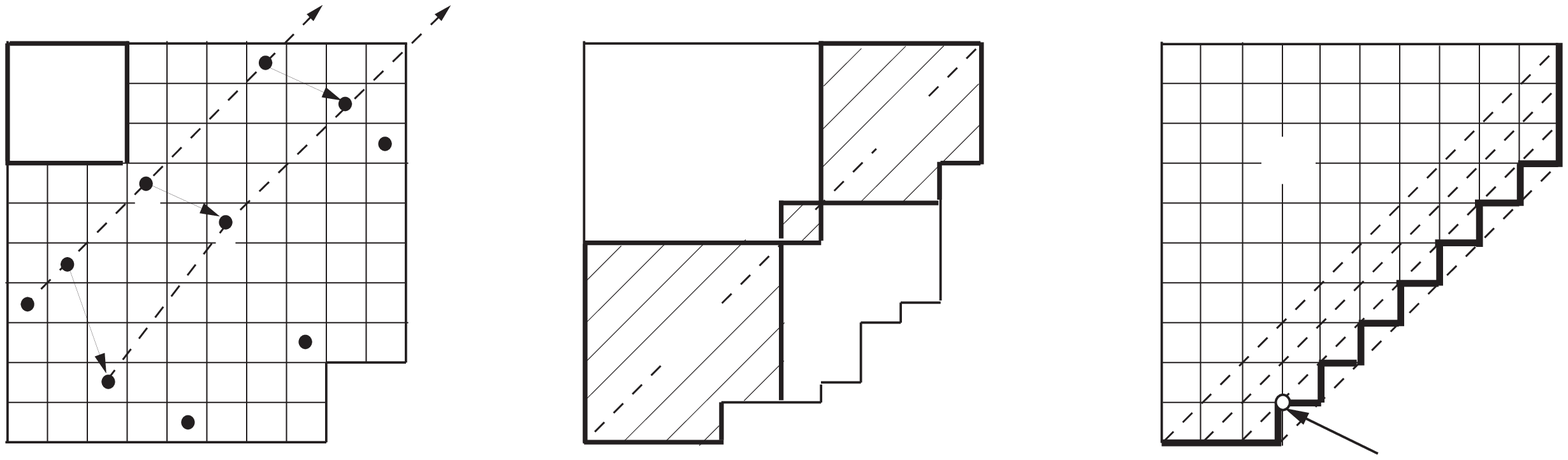}
\hspace*{-20mm}\caption{(a) $T^1$ and $T^2$ \quad\quad\quad\quad(b)
Lemma~\ref{diagonal cells} \quad\quad\quad(c) $2$-critical
$P$ in $St^3_{10}$}
\label{First and Second Subsequences1}
\end{center} 
\end{figure}
\noindent 
Observe that $T^1$ and $T^2$ are increasing subsequences of $T$.
Figure~\ref{First and Second Subsequences1}a depicts $T^1$ and $T^2$
(via dashed lines) and three instances of $\alpha_i\in T^1$
$(21)$-dominating $\beta_j\in T^2$ (via solid arrows).

\subsection{Diagonal Properties and Critical Points}
\label{Diagonal Properties and Critical Points}
We address now the relative positioning of an arbitrary transversal
within its Young diagram.

\begin{lem}
Let $T\in S_Y$ and let $c$ be a cell on the diagonal $d(Y)$. Then the
rectangle $Y_c$ contains some element of $T^1$. Consequently, all
elements of the first subsequence $T^1$ are on or above $d(Y)$.
\label{diagonal cells}
\end{lem}

\noindent{\sc Proof:} Suppose $Y_c$ contains no elements of $T$. But
there is no transversal of $Y$ to sustain such a big empty
rectangle. Indeed, since $c\in d(Y)$, $Y^{\bar{c}}$ is a proper Young
subdiagram of $Y$, say of size $k$, and there are no elements of $T$
above $Y^{\bar{c}}$. Thus, the first $k$ columns of $Y$ must have
their $1$'s within $Y^{\bar{c}}$, and $T$ induces a transversal
$T^{\bar{c}}$ of $Y^{\bar{c}}$.  Analogously, $T$ induces a transversal
$\phantom{}_{\bar{c}}T$ of $\phantom{}_{\bar{c}}Y$. Hence, $T$ must
split into $T=T^{\bar{c}}\oplus T(c) \oplus \phantom{}_{\bar{c}}T$,
where $T(c)$ is a transversal of the cell $c$ (cf. Fig.~\ref{First and
Second Subsequences1}b, where $T$ is concentrated in the 3 shaded
subboards). But cell $c$ is empty by the supposition, a
contradiction. Therefore, $Y_c$ does contain some element $\gamma\in
T$. Since either $\gamma\in T^1$ or $\gamma$ is $(21)$-dominated by
some $\alpha\in T^1$, we conclude that $Y_c$ contains an element of
$T^1$.

If some $\alpha_i\in T^1$ is below the diagonal $d(Y)$, then the
rectangle $Y_{\bar{\alpha_i}}$ contains a cell $c$ on $d(Y)$, and
$Y_c$ is empty, a contradiction with the previous
paragraph. Therefore, $T^1$'s elements are on or above $d(Y)$. \qed

\smallskip
By the {\it border} of a Young diagram $Y$ we mean the path that
starts at the bottom left corner of $Y$, follows $Y$'s outline
{\it below and to the right} of $d(Y)$, and ends at the top right
corner of $Y$.

\begin{defn}
{\rm For a Young diagram $Y$, define the {\it $i$-th diagonal
$d_i(Y)$} as follows: starting from the bottom left corner of $Y$,
move $i$ cells to the right, draw a parallel line to $d(Y)$ until it
goes through the rightmost column of $Y$; the resulting segment is
$d_i(Y)$.  For $i\geq 1$, denote by $St^i_n$ the {\it $i$-th
Staircase} Young diagram of size $n$ whose border is the stepwise path
from the bottom left corner to the top right corner of $Y$ that zigzags
between $d_{i-1}(Y)$ and $d_i(Y)$ (cf. Fig.~\ref{First and Second
Subsequences1}c for $d_i(Y)$ with $0\leq i \leq 3$, and $St^3_{10}$.)}
\label{diagonals}
\end{defn}

We distinguish between $d_0(Y)$, which is a segment going through
$Y$'s diagonal {\it grid points}, and $d(Y)$, which is the union of
all diagonal {\it cells} of $Y$. 

\begin{defn}
{\rm A grid point $P$ on $Y$'s border is called a {\it critical point
of} $Y$ if $Y$'s border goes upwards to enter $P$ and then goes to the
right to leave $P$. If in addition $P\in d_i(Y)$, then $P$ is
called an {\it i}-critical point of $Y$.}
\label{critical points}
\end{defn}

Figure~\ref{First and Second Subsequences1}c shows the {\it bottom}
$2$-critical point $P$ of $St^3_{10}$. Note that $St^n_n=M_n$ is the
only Young diagram of size $n$ with no critical points, while $St^1_n$
has the largest number of critical points. Also, for any critical
point $P$, the subboard $\phantom{}^PY$ has no cells and consists only
of the point $P$, while $Y_P$ is a rectangle.

\begin{lem} If $P$ is an $i$-critical point of $Y$ and $T\in
  S_Y$, then the rectangle $Y_P$ contains exactly $i$ elements of $T$.
\label{rectangle fill}
\end{lem}

\noindent{\sc Proof:} Let $Y$ have exactly $k$ rows above $P$.  Since
$P\in d_i(Y)$, the subboard $\phantom{}_PY$ has $k$ rows and $k-i$
columns; the latter are in fact all columns of $Y$ which are to the
right of $P$, and therefore each of these $k-i$ columns contains
exactly 1 element of $T$. Hence $k-i$ of $\phantom{}_PY$'s rows
contain an element of $T$, while $i$ rows of $\phantom{}_PY$ are empty
(cf. Fig.~\ref{General 0-1 splittings}a-b for $i=0,1$ and
Fig.~\ref{General 2-splitting}a for $i=2$.)

On the other hand, each of the top $k$ rows of $Y$ is split between
the rectangle $Y_P$ and the subboard $\phantom{}_PY$.  From the
viewpoint of $Y_P$, the above observations mean that $k-i$ rows of
$Y_P$ are empty, while exactly $i$ rows of $Y_P$ contain an element of
$T$. Thus, $|T_P|=i$. \qed

\subsection{Definition of the map $\zeta_P$}
\label{general zeta_P map}
For an $i$-critical point $P$ in $Y$, let $Q,R\in d_0(Y)$ be the
diagonal grid points of $Y$ to the left of, respectively above,
$P$. Then $\phantom{}_QY$ and $Y^R$ are proper Young subdiagrams
(cf. Fig.~\ref{General 0-1 splittings}a-b and Fig.~\ref{General
2-splitting}a.) 

Fix $T\in S_Y$. Lemma~\ref{rectangle fill} ensures that rectangle
$Y_P$ contains exactly $i$ elements of $T$, which form some
subsequence $\alpha=(\alpha_1,\alpha_2,...,\alpha_i)$. While
preserving the pattern $\alpha$, we can simultaneously pull downward
all $\alpha_i$'s until they become the top $i$ elements in a
transversal $T_1$ of $Y^R$, and we can also push all $\alpha_i$'s to
the right until they become the $i$ leftmost elements of a transversal
$T_2$ of $\phantom{}_QY$. These operations define an injective map
\[\zeta_P:S_Y\hookrightarrow S_{Y^R}\times S_{\phantom{}_QY}
\,\,\text{where}\,\,\zeta_P(T)=(T_1,T_2).\] For example,
Fig.~\ref{General 0-1 splittings}b-c show
$\zeta_P(31628547)=(3142,35214)$ with $i=1$ and $\alpha_1=6$, while
Fig.~\ref{General 2-splitting} shows
$\zeta_P(831629547)=(53142,536214)$ with $i=2$ and $\alpha_1=\alpha=8$
and $\alpha_2=\beta=6$. Since $\phantom{}^PY$ has no cells, any
subsequence of $T$ landing inside $Y$ must be contained either
entirely in the rows of $Y$ above $P$, or entirely in the columns of
$Y$ to the left of $P$. Consequently,

\begin{lem}
For any pattern $\sigma$, $T$ avoids $\sigma$ on $Y$ if and only if
the components $T_1$ and $T_2$ of $\zeta_P(T)$ avoid $\sigma$ on
${Y^R}$ and ${\phantom{}_QY}$, respectively. In particular, $\zeta_P$
respects pattern-avoidance and we can restrict
\label{zeta_P respects patterns} 
$\zeta_P:S_Y(\sigma)\hookrightarrow S_{Y^R}(\sigma)\times
S_{\phantom{}_QY}(\sigma)$.
\end{lem}

\subsection{Critical Splittings induced by $\zeta_P$}
\label{subsection on critical splittings}
\begin{prop} If $P$ is a $0$- or $1$-critical point of $Y$, then   
$S_Y(\sigma)\stackrel{\zeta_P}{\cong} S_{Y^R}(\sigma)\times
S_{\phantom{}_QY}(\sigma)$ for any $\sigma\in S_k$.  
\label{0-1-splitting}
\end{prop}

\noindent{\sc Proof:} Fix $T\in S_Y$ and let $\sigma$ be any
permutation. A $0$-critical point $P$ coincides with the points $Q$
and $R$ in the definition of $\zeta_P$, and the rectangle $Y_P$ has no
elements of $T$ by Lemma~\ref{rectangle fill} (cf. Fig.~\ref{General
0-1 splittings}a.) Thus, $\zeta_P:S_Y(\sigma)\hookrightarrow
S_{Y^R}(\sigma)\times S_{\phantom{}_QY}(\sigma)$ simply restricts
$T|_{Y^P}=T_1$ and $T|_{\phantom{}_PY}=T_2$; combined with
Lemma~\ref{zeta_P respects patterns}, this yields invertibility of
$\zeta_P$.  In this case, we say that $\zeta_P$ induces the $0$-{\it
splitting} $T=T|_{Y^P}\oplus T|_{\phantom{}_PY}$.
\begin{figure}[h]
\labellist
\small\hair 2pt
\pinlabel $\scriptstyle{\phantom{}_QY}$ at 517 603
\pinlabel $\scriptstyle{\alpha_Q}$ at 462 620
\pinlabel $\scriptstyle{Y^R}$ at 408 517
\pinlabel $\scriptstyle{\alpha_R}$ at 424 566
\pinlabel $\scriptstyle{\phantom{}_QY}$ at 270 593
\pinlabel $\zeta_P$ at 328 587
\pinlabel $\scriptstyle{\alpha}$ at 181 611
\pinlabel $\scriptstyle{R}$ at 226 585
\pinlabel $\scriptstyle{P}$ at 226 548
\pinlabel $c$ at 208 567.5
\pinlabel $\scriptstyle{Q}$ at 190 548
\pinlabel $Y_P$ at 163 633
\pinlabel $\scriptstyle{Y^R}$ at 184 526
\pinlabel $\scriptstyle{\phantom{}_QY}$ at 75 603
\pinlabel $\scriptstyle{P}$ at 10 558
\pinlabel $\scriptstyle{Y^R}$ at -34 518
\pinlabel $Y_P$ at -55 642
\endlabellist
\begin{center}
\includegraphics[width=4.5in]{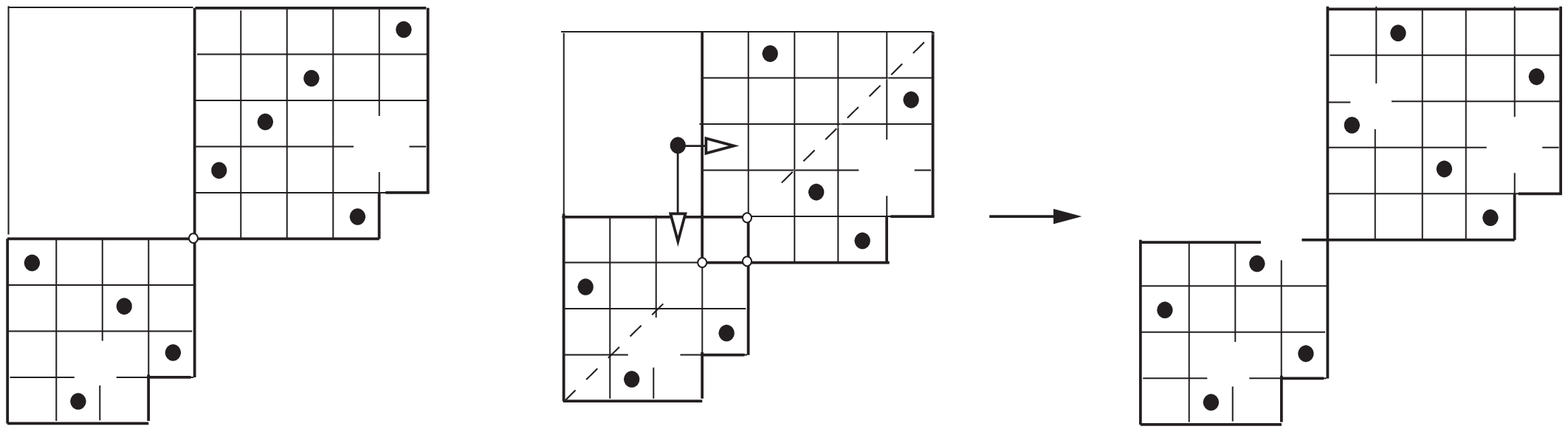}
\caption{(a) $0$-splitting \hspace*{20mm} (b)-(c) $1$-splitting \hspace*{35mm}}
\label{General 0-1 splittings}
\end{center} 
\end{figure}

Now, consider the case of a $1$-critical point $P$
(cf. Fig.~\ref{General 0-1 splittings}b-c.) Let $c$ be the cell whose
bottom right corner is $P$. Then $c$ lies on the diagonal $d(Y)$, and
$Q$ and $R$ are also respective corners of $c$. Let $c=(k,m)$ where
$k$ is $c$'s row and $m$ is $c$'s column in $Y$. By
Lemma~\ref{rectangle fill}, the rectangle $Y_P$ has exactly one
element of $T$: call it $\alpha$, and let it be in position $(i,j)$ in
$Y$. To form transversals $T_1\in S_{Y^R}$ and $T_2\in
S_{\phantom{}_QY}$, $\zeta_P$ replaces $\alpha$ by $\alpha_R$ in
position $(k,j)$ and $\alpha_Q$ in position $(i,m)$, respectively.

It is not hard to see that $\zeta_P$ is surjective.  Indeed, start
with $(T_1,T_2)\in S_{Y^R}\times S_{\phantom{}_QY}$. If $T_1$ has its
top element $\alpha_Q$ in its $j$-th column, and $T_2$ has its
leftmost element $\alpha_R$ in its $i$-th row, we can reconstruct the
unique $\alpha\in Y_P$ by replacing $(\alpha_R,\alpha_Q)$ by an
element in position $(i,j)$ and leaving the rest of $T_1$ and $T_2$
fixed. Combining this with Lemma~\ref{zeta_P respects patterns} yields
the wanted isomorphism $\zeta_P$ on $S_Y(\sigma)$. In this case, we
say that $\zeta_P$ induces the $1$-{\it splitting}
$T=T|_{Y^P}\oplus_{\scriptscriptstyle{1}} T|_{\phantom{}_PY}$. \qed

\medskip
As expected, $i$-critical points for larger $i$ complicate matters,
and in general, it is not possible to derive such nice splittings
of transversals. Below we describe the image $\zeta_P(S_Y(\sigma))$
for a $2$-critical point $P$.

\begin{defn}
{\rm Let $S_{{\scriptscriptstyle{\nearrow}}Y}(\sigma)$, respectively
$S^{\scriptscriptstyle{\nearrow}}_Y(\sigma)$, be the set of
transversals $T$ in $S_Y(\sigma)$ whose two leftmost, respectively
two top, elements form an increasing subsequence of $T$.  Define
analogously $S^{\scriptscriptstyle{\searrow}}_Y(\sigma)$ and
$S_{{\scriptscriptstyle{\searrow}}Y}(\sigma)$ with appropriate
replacement of $\nearrow$ by $\searrow$.}
\label{technical defns}
\end{defn}

We will also need the notation
$S^{\,\,\,\,\,\,\scriptscriptstyle{\searrow}}_{{\scriptscriptstyle{\nearrow}}
Y}(\sigma)=S^{\scriptscriptstyle{\searrow}}_Y(\sigma)\cap
S_{{\scriptscriptstyle{\nearrow}}Y}(\sigma)$. As with previous
notation, this one preserves the relative position of the involved
objects, in this case -- $Y$ and its two (top and/or leftmost)
subsequences of length 2. The $\nearrow$ and $\searrow$ arrows can be
arbitrarily switched to denote the corresponding other subsets of
transversals.

\begin{lem}
If $P$ is a $2$-critical point of $Y$, then for any $\sigma\in S_k$:
\begin{equation}
S_Y(\sigma) \stackrel{\zeta_P}{\cong}
S^{\scriptscriptstyle{\nearrow}}_{Y^R}(\sigma)\times
S_{{\scriptscriptstyle{\nearrow}}\phantom{}_QY}(\sigma) \sqcup
S^{\scriptscriptstyle{\searrow}}_{Y^R}(\sigma)\times
S_{{\scriptscriptstyle{\searrow}}\phantom{}_QY}(\sigma).
\label{general 2-splitting}
\end{equation}
\label{General 2-splitting lemma}
\end{lem}
\noindent{\sc Proof:} Start with $T\in S_Y$. By Lemma~\ref{rectangle fill},
we may assume that $\alpha$ and $\beta$ are the only elements of $T$
in rectangle $Y_P$, with $\alpha$ to the left of $\beta$.
\begin{figure}[h]
\labellist
\small\hair 2pt
\pinlabel $\scriptstyle{\phantom{}_QY}$ at 396 611
\pinlabel $\scriptstyle{\alpha_Q}$ at 324 665
\pinlabel $\scriptstyle{\beta}$ at 47 608
\pinlabel $\scriptstyle{\beta_R}$ at 288 557
\pinlabel $\scriptstyle{\alpha_R}$ at 234 575
\pinlabel $\scriptstyle{Y^R}$ at 272 508
\pinlabel $\scriptstyle{\phantom{}_QY}$ at 117 593
\pinlabel $\zeta_P$ at 175 589
\pinlabel $\scriptstyle{\alpha}$ at -20 648
\pinlabel $\scriptstyle{\beta_Q}$ at 345 630
\pinlabel $\scriptstyle{R}$ at 73 602
\pinlabel $\scriptstyle{P}$ at 73 548
\pinlabel $\scriptstyle{Q}$ at 20 548
\pinlabel $\scriptstyle{Y^R}$ at 28 526
\pinlabel $Y_P$ at 45 647
\endlabellist
\begin{center}
\includegraphics[width=3in]{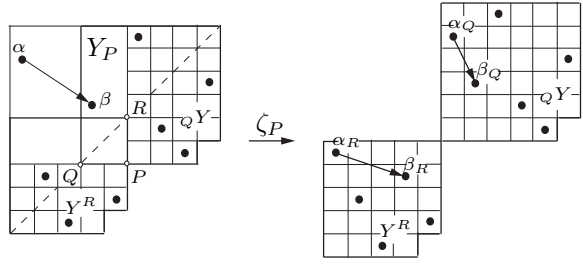}
\caption{$\zeta_P(T)=(T_1,T_2)$ on $Y^R\times_{\,Q}\!Y$
with $(\alpha\beta)\!\!\searrow$ in $Y_P$}
\label{General 2-splitting}
\end{center} 
\end{figure}
Depending on whether $(\alpha\beta)\!\nearrow$ or $\searrow$, either
component $T_1\in S_{Y^R}(\sigma)$ has its top two elements
$(\alpha_R,\beta_R)\!\!\nearrow$ {\it and} component $T_2\in
S_{\phantom{}_QY}(\sigma)$ has its two leftmost elements
$(\alpha_Q,\beta_Q)\!\nearrow$, {\it or} both of these subsequences
are decreasing. For instance, Figure~\ref{General 2-splitting} depicts
the case $(\alpha\beta)\!\searrow$.

Conversely, start with $(T_1,T_2)\in
S^{\scriptscriptstyle{\nearrow}}_{Y^R}(\sigma)\times
S_{{\scriptscriptstyle{\nearrow}}\phantom{}_QY}(\sigma)\sqcup
S^{\scriptscriptstyle{\searrow}}_{Y^R}(\sigma)\times
S_{{\scriptscriptstyle{\searrow}}\phantom{}_QY}(\sigma)$. If
$(\alpha_R,\beta_R)$ and $(\alpha_Q,\beta_Q)$ are the top two,
respectively, the leftmost two, elements of $T_1$ and $T_2$, they form
the same length-2 pattern, say, they are both decreasing. This makes
it possible to pull back $\alpha_R$ and $\alpha_Q$ to an element
$\alpha$ in rectangle $Y_P$, and pull back $\beta_R$ and $\beta_Q$ to
an element $\beta$ in rectangle $Y_P$, so that $(\alpha,\beta)$ is
also decreasing and
$\zeta_P(\alpha,\beta)=(\alpha_R,\beta_R)\times(\alpha_Q,\beta_Q)$ in
$Y^R\times_{\,Q}\!Y$. This discussion establishes the two isomorphisms
$S^{\scriptscriptstyle{\nearrow}}_Y(\sigma)\cong
S^{\scriptscriptstyle{\nearrow}}_{Y^R}(\sigma)\times
S_{{\scriptscriptstyle{\nearrow}}\phantom{}_QY}(\sigma)$ and
$S^{\scriptscriptstyle{\searrow}}_Y(\sigma)\cong
S^{\scriptscriptstyle{\searrow}}_{Y^R}(\sigma)\times
S_{{\scriptscriptstyle{\searrow}}\phantom{}_QY}(\sigma)$, and since
$S_Y(\sigma)=S^{\scriptscriptstyle{\nearrow}}_Y(\sigma)\sqcup
S^{\scriptscriptstyle{\searrow}}_Y(\sigma)$, we deduce (\ref{general
2-splitting}).  \qed

\medskip
The reader can prove a similar splitting for an
$i$-critical point $P$ with $i\geq 3$ and $\sigma\in S_k$:
\[\zeta_P:S_Y(\sigma)\cong \bigsqcup_{\tau\in S_i} 
S^{\scriptstyle{\tau}}_{Y^R}(\sigma)\times
S_{{\scriptstyle{\tau}}\phantom{}_QY}(\sigma),\] where in the
notations $S^{\scriptstyle{\tau}}_{Y^R}(\sigma)$ and
$S_{{\scriptstyle{\tau}}\phantom{}_QY}(\sigma)$ the patterns $\tau\in
S_i$ have replaced the previously used $\nearrow=(12)$ and
$\searrow=(21)$ in $S_2$. In order for this isomorphism to be useful,
one should be able to enumerate the components
$S^{\scriptstyle{\tau}}_{Y^R}(\sigma)$ and
$S_{{\scriptstyle{\tau}}\phantom{}_QY}(\sigma)$; however, for a
general pattern $\sigma$ and high critical index $i$, this question
acquires a level of difficulty at least comparable to that of
Wilf-enumeration $|S_n(\sigma)|$. Fortunately, when $i=2$ and
$\sigma=(312)$ or $(321)$, this enumeration is possible and is carried
out in Section~\ref{strict inequalities 312>321}.

\subsection{The $\sigma\rightarrow\tau$ moves}
\label{general moves}
Let $T\in S_Y$. For any two permutations $\sigma,\tau\in S_k$ we
define a {\it $\sigma\rightarrow\tau$ move} on $T$ as follows: if
$(\alpha_1\alpha_2\cdots\alpha_k)$ is a $\sigma$-subpattern of $T$ in
$Y$, we rearrange the $\alpha_i$'s within the $k\times k$ matrix they
generate so as to obtain a $\tau$-subpattern
$(\beta_1\beta_2\cdots\beta_k)$ in $Y$.  The inverse operation is
obviously a $\tau\rightarrow\sigma$ move. A sequence of
$\sigma\rightarrow\tau$ moves that starts with a transversal $T$ is
called ``a sequence of $\sigma\rightarrow\tau$ moves {\it on $T$}''.

For example, if $(\alpha\beta\gamma)$ is a $(213)$-pattern in $T$
landing in $Y$, a $(213)\rightarrow (123)$ move switches the places
of $\alpha$ and $\beta$ to obtain $(\beta\alpha\gamma)\approx (123)$
in $Y$. Throughout the paper, we will use two instances of
$\sigma\rightarrow\tau$ moves: $(213)\rightarrow (123)$ and
$(312)\rightarrow (321)$ moves, along with their inverses. In
particular, we will construct maps \[S_Y(213)\hookrightarrow
S_Y(123)\cong S_Y(321)\twoheadleftarrow S_Y(312),\] and pose
questions about the general maps $\phi:S_Y(\tau)\rightarrow
S_Y(\sigma)$ that are induced under certain circumstances by a
sequence of $\sigma\rightarrow \tau$ moves in $Y$.

\section{Proof of the Inequality $S_Y(312)\geq S_Y(321)$}
\label{SY(321) < SY(312)}

In this section we prove that $(321)\preceq_s (312)$. Since
$(321)\sim_s (123)$, this will establish the required in Theorem
\ref{SWOS3} inequalities $|S_Y(123)|\leq |S_Y(312)|$ for all Young
diagrams $Y$.  The strategy is to describe the structures of each set
$S_Y(321)$ and $S_Y(312)$, use this information to define a canonical
map $\phi:S_Y(312)\rightarrow S_Y(321)$, and finally prove that $\phi$ is
surjective.

\subsection{The structure of $T\in S_Y(321)$}
\label{structure of S_Y(321)}
$T$ is the disjoint union of its first and second subsequences: $T=T^1
\sqcup T^2$. Indeed, if there were some $\gamma\in
T\backslash\{T^1\cup T^2\}$, then $Y_{\bar{\gamma}}$ would contain
some element $\beta\in T^2$, and hence $Y_{\bar{\beta}}$ would contain
some element $\alpha\in T^1$, so that $(\alpha\beta\gamma)\approx
(321)$ in $T$ and lands in $Y$, a contradiction.

\subsection{The structure of $T\in S_Y(312)$}
\label{structure of S_Y(312)}
Compared to the previous paragraph, the structure here is considerably
more complex. We shall not need all of it in the proof of the
inequality $S_Y(312)\geq S_Y(321)$. Yet, it is enlightening as to why
the proof works and why strict inequalities $S_Y(312)>S_Y(321)$ occur
for some $Y$.  For the remainder of this subsection, we fix some
transversal $T\in S_Y(312)$.

\begin{figure}[h]
\labellist
\small\hair 2pt
\pinlabel $\scriptstyle{\alpha}$ at 61 646
\pinlabel $\scriptstyle{\beta}$ at 80 631
\pinlabel $\scriptstyle{\delta_1}$ at 102 560
\pinlabel $\scriptstyle{\delta_2}$ at 119 523
\pinlabel $\scriptstyle{G_{\beta}}$ at 136 577
\pinlabel $\scriptstyle{\beta}$ at 339 599
\pinlabel $\scriptstyle{\gamma_1}$ at 356 570
\pinlabel $\scriptstyle{\gamma_2}$ at 411 550
\pinlabel $\scriptstyle{\gamma_3}$ at 432 508
\pinlabel $\scriptstyle{\beta}$ at 479 599
\pinlabel $\scriptstyle{\gamma_1}$ at 481 534
\pinlabel $\scriptstyle{\gamma_2}$ at 552 567
\pinlabel $\scriptstyle{\gamma_3}$ at 594 518
\pinlabel $\scriptstyle{c}$ at 543 540
\endlabellist
\begin{center}
\includegraphics[width=4in]{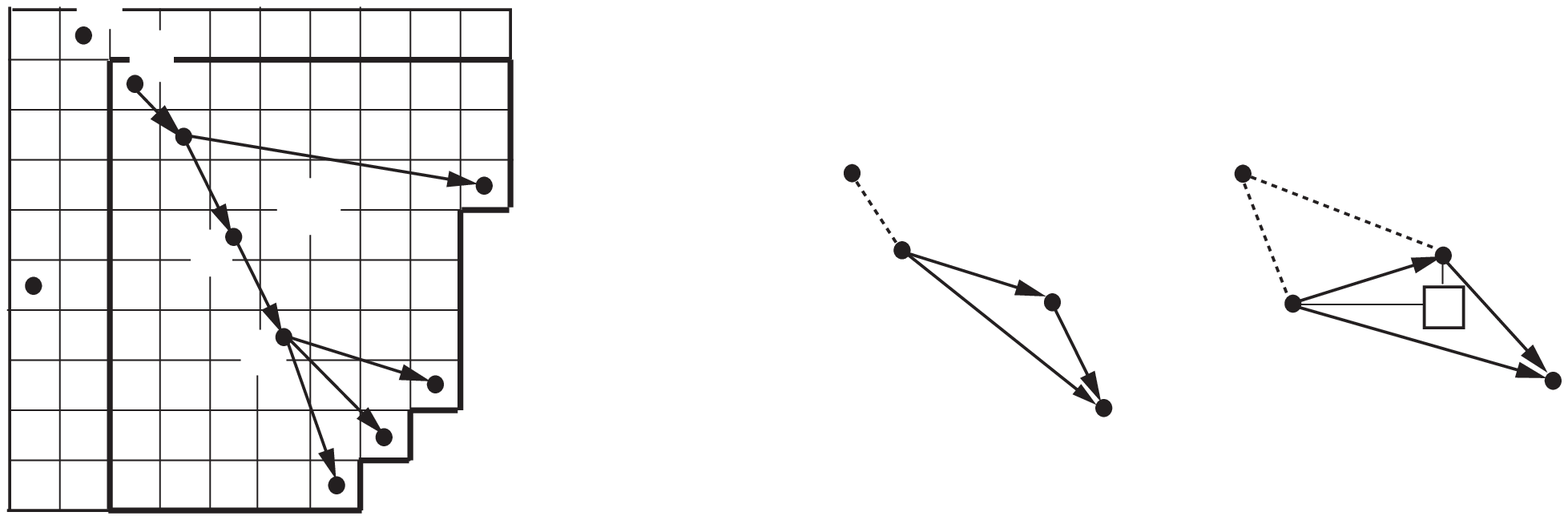}
\caption{$\!\!$(a) Graph $G_{\beta}$ for $\beta\in T^2$;
\hspace*{15mm}(b)-(c) Lemma~\ref{cycle-free}\hspace*{25mm}}
\label{Cycle-free}
\end{center} 
\end{figure}

\begin{defn} {\rm For any $\beta\in T^2$, define a {\it directed
graph} $G_{\beta}$ on the elements of $\phantom{}^{\beta}T$ as
follows: connect by a directed edge
$\stackrel{\longrightarrow}{\delta_1\delta_2}$ any two elements
$\delta_1$ and $\delta_2$ of $\phantom{}^{\beta}T$ such that
$(\delta_1\delta_2)\!\searrow$ and there is no ``intermediate''
$\delta_3\in \phantom{}^{\beta}T$ with
$(\delta_1\delta_3\delta_2)\!\searrow$ (cf. Fig.~\ref{Cycle-free}a.)}
\end{defn}
\begin{lem}{\rm For any $\beta\in T^2$, $\phantom{}^{\beta}T$ avoids
$(12)$ in $Y$. Further, $G_{\beta}$ is connected and, stripping off
the orientation of its edges, cycle-free.
\label{cycle-free}}
\end{lem}

\noindent{\sc Proof:} For the first part, by definition of $\beta\in
T^2$, there is some $\alpha\in T^1\cap T_{\bar{\beta}}$ which
$(21)$-dominates $\beta$.  To avoid the possibility of $\alpha$
playing the role of a ``3'' in a $(312)$-pattern in $T$,
$\phantom{}^{\beta}T$ must avoid $(12)$ in $Y$.

\smallskip
For the second part, $\beta$ $(21)$-dominates any $\gamma\in
\phantom{}^{\bar{\beta}}T$ so that $\gamma$ is connected to at least
one other vertex in $\phantom{}^{\beta}T\cap T_{\bar{\gamma}}$, and
eventually, there is a path starting from $\beta$ and leading to
$\gamma$. Thus, $G_{\beta}$ is connected.

Suppose that there is an (undirected) cycle $\mathcal{C}$ in
$G_{\beta}$. If we start at an arbitrary vertex $\delta\in\mathcal{C}$
and follow $\mathcal{C}$ along the orientation of its edges, we cannot
come back to $\delta$, or else we will have a decreasing sequence
$(\delta,\delta_1,\delta_2,....,\delta_k,\delta)$, which is absurd.

Therefore, going around $\mathcal{C}$ along the edge orientation leads
to a {\it smallest} vertex $\gamma_3$ in $\mathcal{C}$, at which two
edges $\stackrel{\longrightarrow}{\gamma_1\gamma_3}$ and
$\stackrel{\longrightarrow}{\gamma_2\gamma_3}$ terminate (with, say,
$\gamma_1$ before $\gamma_2$.)  If $(\gamma_1\gamma_2)\!\!\searrow$,
then $(\gamma_1\gamma_2\gamma_3)\!\!\searrow$, contradicting the
construction of $G_{\beta}$ without intermediate vertices
(cf. Fig.~\ref{Cycle-free}b.) Thus,
$(\gamma_1\gamma_2)\!\nearrow$. Since $\gamma_3\in \phantom{}^
{\bar{\gamma_1}}T \cap \phantom{}^{\bar{\gamma_2}}T$, the triangle
$\gamma_1\gamma_2\gamma_3$ contains the cell $c$ onto which
$(\gamma_1\gamma_2)$ lands as a $(12)$-pattern, and hence
$\phantom{}^{\beta}T$ also contains $c$ (cf. Fig.~\ref{Cycle-free}c.)
Yet, by the first part of this Lemma, $\phantom{}^{\beta}T$ avoids
$(12)$ in $Y$, a contradiction. Therefore, $G_{\beta}$ has no
(undirected) cycles. \qed

\smallskip 
Lemma \ref{cycle-free} allows us to think of $G_{\beta}$ as an {\it
oriented tree rooted at} $\beta$. Now consider all trees
$G_{\beta_i}$, where
$T^2=(\beta_1,\beta_2,\cdots,\beta_k)\!\!\nearrow$.  For $i<j$, if
$\gamma\in G_{\beta_i}\cap G_{\beta_j}$, then $\gamma\not =
\beta_1,\beta_2$ and $(\beta_i\gamma)\approx
(\beta_j\gamma)\!\!\searrow$. Evidently, if $m$ is between $i$ and
$j$, then $(\beta_m\gamma)\searrow$, so that $\gamma$ is also in
$G_{\beta_m}$ (cf. Fig.~\ref{Lemmas 8,9,10}a.) In other words,

\begin{lem}
Let $G=\cup_{i=1}^k G_{\beta_i}$ be the {\it union of all trees}. Then
each connected component $C_j$ of $G$ is the union of several {\it
consecutive} trees: $C_j=G_{\beta_{k_j}}\cup G_{\beta_{k_j+1}}\cup
G_{\beta_{k_j+2}}\cup \ldots \cup G_{\beta_{k_{j+1}-1}}$.
\label{connected components}
\end{lem}
\begin{figure}[h]
\labellist
\small\hair 2pt
\pinlabel $\scriptstyle{\beta_i}$ at 32 560
\pinlabel $\scriptstyle{\beta_m}$ at  64 574
\pinlabel $\scriptstyle{\beta_j}$ at  98 599
\pinlabel $\scriptstyle{\gamma}$ at 111 525
\pinlabel $\scriptstyle{\beta_i}$ at  228 585
\pinlabel $\scriptstyle{\beta_k}$ at  256 599
\pinlabel $\scriptstyle{\gamma_1}$ at 270 568
\pinlabel $\scriptstyle{\gamma_3}$ at 325 548
\pinlabel $\scriptstyle{\gamma_2}$ at 345 509
\pinlabel $\scriptstyle{\beta_k}$ at  428 585
\pinlabel $\scriptstyle{\beta_l}$ at 455 599
\pinlabel $\scriptstyle{\gamma_j}$ at 475 509
\pinlabel $\scriptstyle{\gamma_i}$ at 526 547
\endlabellist
\begin{center}
\includegraphics[width=4in]{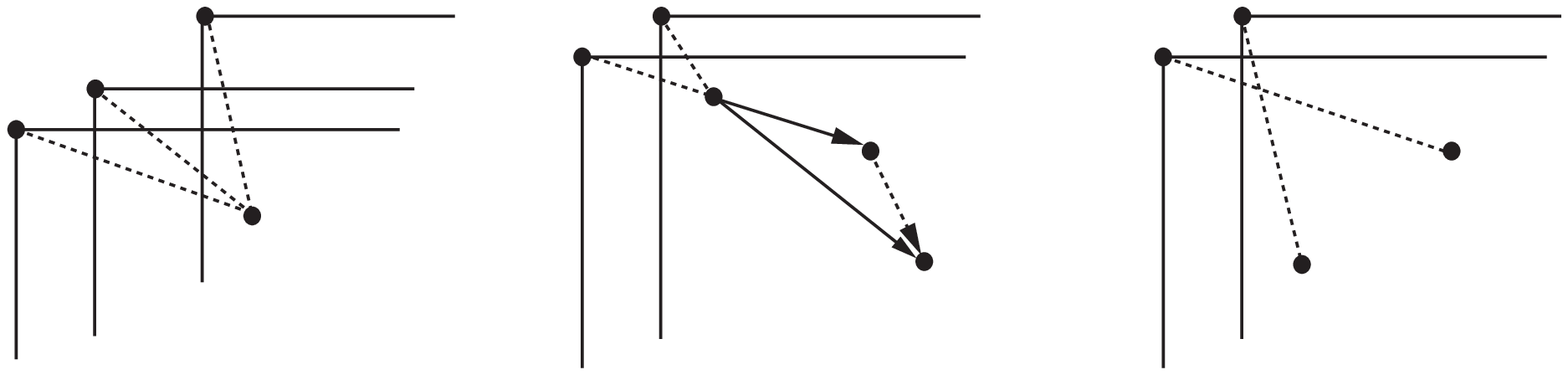}
\caption{Lemmas~\ref{connected components}, \ref{full subgraphs},
\ref{components arrangement}}
\label{Lemmas 8,9,10}
\end{center} 
\end{figure}
By construction, each edge
$\stackrel{\longrightarrow}{\gamma_1\gamma_2}$ of a connected
component $C_j$ is entirely contained in some tree
$G_{\beta_i}$. If $\gamma_1$ and $\gamma_2$ also belong to another
tree $G_{\beta_k}$, then the edge
$\stackrel{\longrightarrow}{\gamma_1\gamma_2}$ must also belong to
$G_{\beta_k}$. Indeed, if not, the (21)-pattern $(\gamma_1\gamma_2)$
requires at least one intermediate vertex $\gamma_3$ in $G_{\beta_k}$:
$(\gamma_1\gamma_3\gamma_2)\!\searrow$ (cf. Fig.~\ref{Lemmas 8,9,10}b.)
But then $\gamma_3$ is also an intermediate vertex in $G_{\beta_i}$,
hence the edge $\stackrel{\longrightarrow}{\gamma_1\gamma_2}$ does not
exist in $G_{\beta_i}$, a contradiction. We conclude that

\begin{lem} 
\label{full subgraphs}
Any tree $G_{\beta_i}$ is a {\it full} subgraph of its connected
component $C_j$.
\end{lem}

Using Lemma \ref{full subgraphs}, we can augment the proof in Lemma
\ref{cycle-free} to derive in an almost identical way that each
connected component $C_j$ has no (undirected) cycles.  Thus, we can
think of each $C_j$ as an oriented ``tree'' rooted at all of its the maximal
elements, i.e. all $\beta_i\in T^2\cap C_j$.

\begin{lem}
\label{components arrangement}
The connected components of $G$ are arranged in an increasing pattern
according to the $\beta_i$'s they contain. More precisely, choose some
$\beta_k\in C_i$ and $\beta_l\in C_j$ such that $k<l$, i.e.
$(\beta_k\beta_l)\!\nearrow$. Then $C_i$ is entirely to the left
and below $C_j$.
\end{lem}

\noindent{\sc Proof:} Consider any $\gamma_i\in C_i$ and $\gamma_j\in
C_j$. If $(\gamma_i\gamma_j)\!\searrow$ or
$(\gamma_j\gamma_i)\!\searrow$, Lemma~\ref{full subgraphs} guarantees
a path between $\gamma_i$ and $\gamma_j$, contradicting $C_i\cap
C_j=\emptyset$.  Thus, $\gamma_i$ and $\gamma_j$ form (in some order)
an increasing sequence. To complete the proof, we need to show
$(\gamma_i\gamma_j)\!\nearrow$.

To the contrary, suppose $(\gamma_j\gamma_i)\!\nearrow$. Because of
Lemma~\ref{connected components} and the arbitrary choice of
$\beta_k\in C_i$ and $\beta_l\in C_j$, we may assume that $\gamma_i\in
G_{\beta_k}\subset C_i$ and $\gamma_j\in G_{\beta_l}\subset C_j$,
i.e. $(\beta_k\gamma_i)\!\searrow$ and $(\beta_l\gamma_j)\!\searrow$
(cf. Fig.~\ref{Lemmas 8,9,10}c.) Putting together all four elements,
we arrive at the subsequence $(\beta_k\beta_l\gamma_j\gamma_i)\approx
(3412)$, which does not necessarily land in $Y$. Then
$(\beta_k\gamma_j)\!\searrow$ so that $\gamma_j\in G_{\beta_k}\subset
C_i$. Thus, $\gamma_j\in C_i\cap C_j=\emptyset$, a contradiction. If
it happens that $\gamma_i=\beta_k$, or $\gamma_j=\beta_l$, or both,
immediate contradictions in the overall arrangement arise.

We conclude that $(\gamma_i\gamma_j)\!\nearrow$, so that $C_i$ is
entirely to the left and below $C_j$. \qed

\medskip
Thus, the connected components of $G$ are arranged in a increasing
diagonal fashion, symbolically,
$G=(C_1,C_2,...,C_k)\!\nearrow$. Correspondingly, the whole transversal
$T\in S_Y(312)$ is the disjoint union of the increasing subsequence
$T^1$ and all the vertices $|C_i|$ of the $C_i$'s: 
\begin{equation}
T=T^1\sqcup |G|=T^1\sqcup_i |C_i|.
\label{Splitting of T in SY(312)}
\end{equation}
This description of a $(312)$-avoiding transversal in $Y$ is only {\it
partial} (transversals satisfying it do not necessarily avoid
$(312)$), but sufficient for our purpose to explain why $(312)$ is
easier to avoid than $(321)$ on Young diagrams $Y$ (cf. also
Section~\ref{strict inequalities 312>321}.)  In particular, the
description involves only the elements of the transversal $T$, while
it is possible to extend it to the whole Young diagram $Y$. To this
end, let $Y_j$ be the Young subdiagram of $Y$ obtained after reducing
$Y$ along all elements of $T$ not in $C_j$; one can think of $Y_j$ as
the Young subdiagram induced by the elements of $C_j$.  Since the
$C_j$'s are disjoint, the $Y_j$'s are disjoint, and we leave it to the
reader to deduce in a similar fashion as above:

\begin{lem}
The Young subdiagrams $Y_i$ are arranged in a increasing diagonal
fashion: $Y=(Y_1,Y_2,...,Y_k)\!\nearrow$.
\label{diagrams arrangement}
\end{lem}

\subsection{Definition of the map $\phi:S_Y(312)\rightarrow S_Y(321)$.}
Fix a transversal $T\in S_Y(312)$, and decompose $T=T^1\sqcup |G|$ as
in (\ref{Splitting of T in SY(312)}) (cf. Fig.~\ref{definition of map
312-321}.) Reducing $Y$ along $T^1$ leaves the pattern of $|G|$ in
a Young subdiagram $Y_0=Y\!\!\big/_{\displaystyle{\!{T^1}}}$. Since
$|G|$ represents a transversal of $Y_0$, then $Y_0$ is proper, with
diagonal $d(Y_0)$. Replacing $|G|$ by the increasing pattern
$I_s=(123...s)$ along $d(Y_0)$ produces another transversal of
$Y_0$. We reintroduce the rows and columns of the previously reduced
subsequence $T^1$ to obtain our original Young diagram $Y$ with a new
transversal $\phi(T)=T^1\sqcup I_s$. Since $\phi(T)$ is partitioned
into two increasing subsequences, $\phi(T)$ avoids $(321)$ and thus
$\phi:S_Y(312)\rightarrow S_Y(321)$ is well-defined.
\begin{figure}[h]
\labellist
\small\hair 2pt
\pinlabel $\scriptstyle{\alpha}$ at -17 647
\pinlabel $\scriptstyle{\beta}$ at 18 617
\pinlabel $\scriptstyle{T^1}$ at 16 672
\pinlabel $\scriptstyle{T^2}$ at 46 672
\pinlabel $\scriptstyle{G}$ at 34 520
\pinlabel $\scriptstyle{\beta^{\prime}}$ at 178 624
\pinlabel $\scriptstyle{|G|^1=T^2}$ at 227 642
\pinlabel $\scriptstyle{I_s\,\,\text{on}\,\,Y_0}$ at 396 642
\pinlabel $\scriptstyle{\delta^{\prime}}$ at 409 585
\pinlabel $\scriptstyle{T^1=(\phi(T))^1}$ at 628 672 
\pinlabel $\scriptstyle{(\phi(T))^2}$ at 718 672
\pinlabel $\scriptstyle{\alpha}$ at 560 648
\pinlabel $\scriptstyle{\delta}$ at 648 559
\endlabellist
\begin{center}
\includegraphics[width=5in]{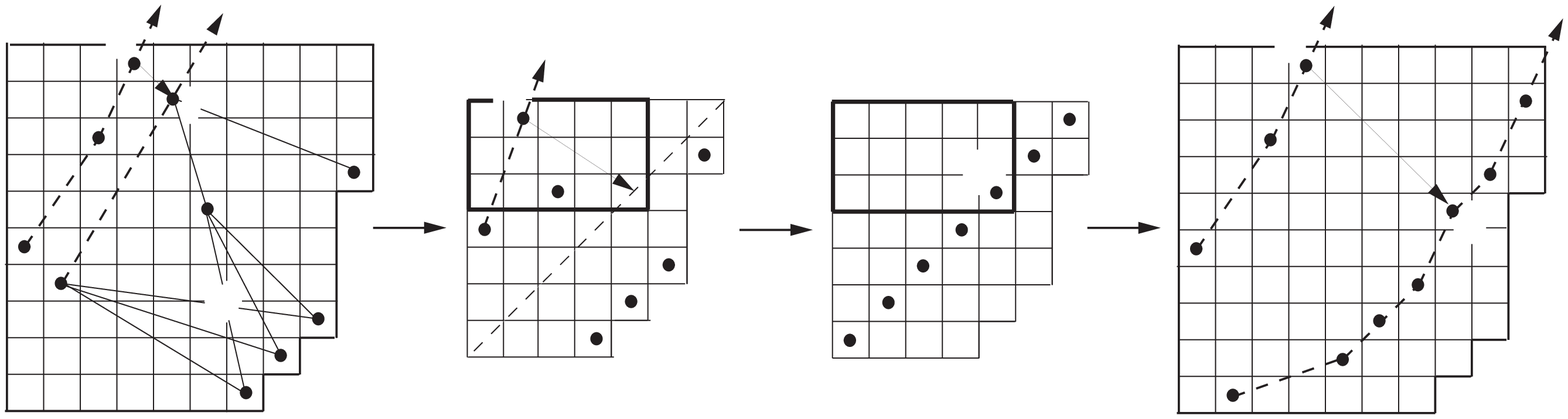}
\caption{$T\in
  S_Y(312)\,\,\rightarrow\,\,T|_{Y_0}\,\,\rightarrow\,\,I_s\,\,\rightarrow\,\,\phi(T)\in S_Y(321)$ }
\label{definition of map 312-321}
\end{center} 
\end{figure}

\subsection{Surjectivity of $\phi$.} 
To show that $\phi$ is surjective, we will first show 

\begin{lem} $\phi$ preserves $T^1$, i.e. $(\phi(T))^1=T^1$. 
\label{phi preserves T^1}
\end{lem}

\noindent{\sc Proof:} Since the elements of $T^1$ are fixed by $\phi$,
it suffices to show that any other element $\delta\in
\phi(T)\backslash T^1$, is $(21)$-dominated by some $\alpha\in T^1$,
implying $\delta\not\in (\phi(T))^1$.

Thus, start with $\delta\in \phi(T)\backslash T^1$ and pull it back to
$\delta^{\prime}\in I_s$ on $Y_0$ (cf. Fig.~\ref{definition of map
312-321}d-c.)  Consider the rectangle $(Y_0)_{\delta^{\prime}}$: since
the cell of $\delta^{\prime}$ is on the diagonal $d(Y_0)$, the proof
of Lemma \ref{diagonal cells} implies that the transversal $T|_{Y_0}$
cannot sustain such a big empty rectangle. On the other hand, in the
reduction $Y_0=Y\!\!\big/_{\displaystyle{\!{T^1}}}$, the first
sequence of the transversal $|G|$ coincides with the original second
sequence $T^2$ in $Y$: $|G|^1=T^2$ (cf. Fig.~\ref{definition of map
312-321}b.) Putting together these considerations implies the
existence of some $\beta^{\prime}\in |G|^1$ in the rectangle
$(Y_0)_{\delta^{\prime}}$.  Pulling $\beta^{\prime}$ to $\beta\in T^2$
on $Y$, we deduce that some $\alpha\in T^1$ $(21)$-dominates $\beta$
(cf. Fig.~\ref{definition of map 312-321}a.) Comparing the relative
positions of $\alpha$, $\beta$ and $\delta$ in $Y$, we conclude that
$\alpha$ $(21)$-dominates $\delta$ in $\phi(T)$.  Therefore,
$\delta\not\in (\phi(T))^1$, and as noted above, this means
$(\phi(T))^1=T^1$.  \qed

\medskip
We can also think of $\phi$ in terms of the canonical decomposition in
(\ref{Splitting of T in SY(312)}) of $T\in S_Y(312)$: replace every
connected component $C_i$ in $G$ by the increasing sequence $I_i$ in
the Young subdiagram $Y_i$. Then $I_s=I_1\sqcup I_2\sqcup \ldots
\sqcup I_k$. This works since the $C_i$'s and the $Y_i$'s are
independent of each other and arranged in an increasing sequence in
$Y$.

\begin{lem}
Given a fixed increasing sequence $L$ of dots in $Y$, there is at most
one transversal $T\in S_Y(321)$ for which $T^1=L$.
\label{L in Y}
\end{lem}

\noindent
{\sc Proof:} If $T\in S_Y(321)$ is such a transversal, then
$T=T^1\sqcup T^2$ with $T^1=L$. Reducing
$T\!\!\big/_{\displaystyle{\!{L}}}$ leaves $T^2$, which must be an
increasing sequence in and a transversal of the resulting Young
diagram $Y\!\!\big/_{\displaystyle{\!{L}}}$; yet, there is only one
such sequence in $Y\!\!\big/_{\displaystyle{\!{L}}}$, namely, its
diagonal sequence $I_s$. This uniquely defines $T^2$, and since the
rest of $T$ is the fixed $L$, it uniquely defines $T:=L\sqcup I_s$
too.  Of course, after putting back $L$ and $I_s=T^2$ to $Y$, it may
turn out that the newly added points of $T^2$ violate the definition
of $L$ by participating in $T^1$, so in this case there would be no
$T\in S_Y(321)$ with $T^1=L$.\qed

\begin{prop}
The map $\phi:S_Y(312)\rightarrow S_Y(321)$ is surjective.
\label{surjectivity1}
\end{prop}

\noindent{\sc Proof:} Let $Q\in S_Y(321)$, and decompose $Q=Q^1\sqcup
Q^2$. We will construct $T\in S_Y(312)$ such that $T^1=Q^1$. For that,
start with $Q$ and apply any sequence of
$(312)\!\rightarrow\!(321)$ moves on $Q$ until there are no more
$(312)$-patterns in $Y$.  Denote the final transversal of $Y$ by
$T$. As an example, reverse the arrow $\phi$ in Figure~\ref{Example 1}
in Section~\ref{strict inequalities 312>321}: depending on the order
of picking the $(312)$-patterns, one can get from $Q=(31524)\in
S_Y(321)$ to $T_1=(31542)$ or $T_2=(32514)$ in $S_Y(312)$.

Each move replaces a $(312)$-pattern in $Y$ with a $(321)$-pattern in
$Y$ by fixing the element playing the role of ``3'', and switching the
other two elements as in $(12)\mapsto (21)$, and thereby increasing
the number of {\it inversions} in the total transversal. Hence the
number of moves cannot exceed $\binom{n}{2}$ and the sequence of moves
eventually terminates with some $T\in S_Y(312)$.

The {\it first subsequences} of the original and of the final
permutation coincide: $T^1=Q^1$. Indeed, none of the moves
$(\alpha_1\alpha_2\alpha_3)\approx(312)\mapsto
(\alpha_1\alpha_3\alpha_2)\approx(321)$ changes the first subsequence,
because $\alpha_1$ $(21)$-dominates the other two elements, whether
before or after the move. Hence $\alpha_2$ and $\alpha_3$ are not in
and cannot land in the first subsequence via the moves, and their
switch certainly does not affect in any way the existing first
subsequence elements. We conclude that $T^1=Q^1$.

By Lemma \ref{phi preserves T^1}, $\phi$ preserves the first
subsequence, so that applying $\phi$ to $T$ yields $\phi(T)\in
S_Y(321)$ with $(\phi(T))^1=T^1=Q^1$. But by Lemma \ref{L in Y},
there is at most one transversal in $S_Y(321)$ with first subsequence
$Q^1$, namely, $Q$.  Thus, $\phi(T)=Q$ and $\phi$ is surjective. \qed

\subsection{Conclusions} 
Proposition~\ref{surjectivity1} implies that for all Young diagrams $Y$:
\[|S_Y(312)|\geq |S_Y(321)|,\] 
which is one of the two inequalities in Theorem \ref{SWOS3}. Therefore,
$(312)\succeq_s (321)$. Now Proposition~\ref{prop-modified-BW} implies
that $(312|\tau)\succeq_s (321|\tau)$ for any permutation
$\tau$; equivalently, for any Young diagram $Y$ we have
$|S_Y(312|\tau)\geq |S_Y(321|\tau)|$.  Consequently, for all $n$:
\[|S_n(312|\tau)|\geq |S_n(321|\tau)|,\]
which completes half of Corollary \ref{Wilf-ordering-corollary}. \qed

\section{Strict Inequalities $S_{Y}(312)>S_{Y}(321)$}
\label{strict inequalities 312>321}
\subsection{Examples of Strict Inequalities} 
Since $\phi:S_Y(312)\twoheadrightarrow S_Y(321)$, a strict inequality
$|S_{Y}(312)|>|S_{Y}(321)|$ occurs exactly when for some $Q\in
S_Y(321)$ the fiber $\phi^{-1}(Q)\subset S_Y(312)$ has more than 1
element. From the proof of Proposition~\ref{surjectivity1}, this
happens exactly when two distinct $T_1,T_2\in S_Y(312)$ have the same
first subsequences: $(T_1)^1=(T_2)^1$.
\begin{figure}[h]
\labellist
\small\hair 2pt
\pinlabel $\scriptstyle{T^1}$ at 148 612
\pinlabel $\scriptstyle{T^1}$ at 302 612
\pinlabel ${\text{and}}$ at 193 553
\pinlabel ${\phi}$ at 353 561
\pinlabel $\scriptstyle{T^1\!= Q^1}$ at 457 612
\pinlabel $\scriptstyle{I_3}$ at 496 610 
\endlabellist
\begin{center}
\includegraphics[width=3in]{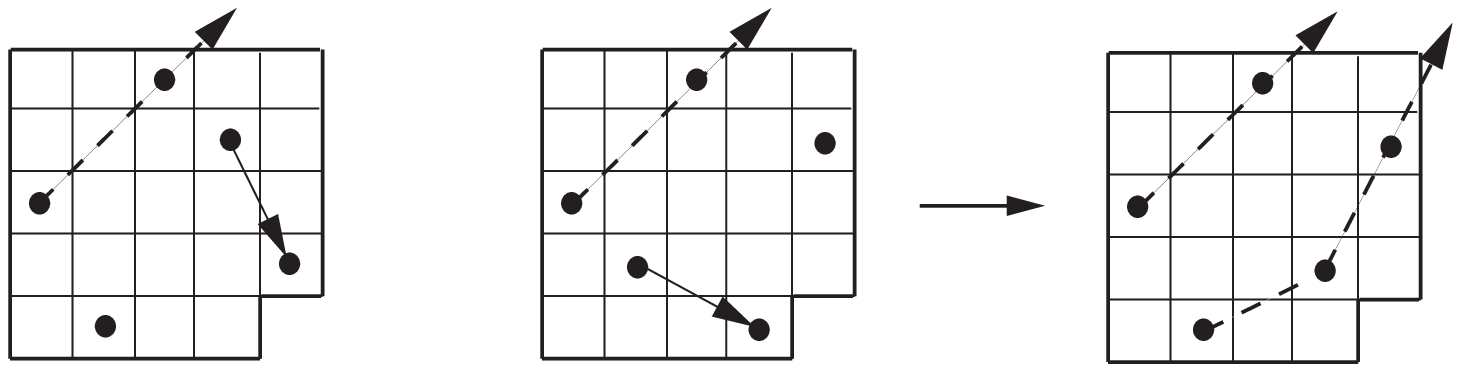}
\caption{$T_1=(31542),\,\,T_2=(32514)\in
  S_{Y_5}(312)\,\,\stackrel{\phi}{\longrightarrow}\,\, Q=(31524)\in
  S_{Y_5}(321)$}
\label{Example 1}
\end{center} 
\end{figure}

\begin{ex}{\rm We revisit
the Young diagram $Y_5=(5,5,5,5,4)$, mentioned in the Introduction. It
is the smallest Young diagram on which $(312)$ is less restrictive
than $(321)$: $|S_Y(312)|=42>41=|S_Y(321)|$. The two sets intersect in
a large subset: $|S_Y(312,321)|=21$, and
$\phi:S_Y(321)\twoheadrightarrow S_Y(312)$ acts as the identity map on
this intersection. Indeed, if $T\in S_Y(312,321)$, then $T=T^1\sqcup
T^2$, so that $I_s\equiv T^2\!\!\nearrow$ and $\phi(T)=T^1\sqcup
I_s=T$. In addition, there are 19 transversals $U\in S_Y(321)$ whose
preimages in $S_Y(312)$ consist of single elements $\phi^{-1}(U)\not
= U$.

As expected, the map $\phi$ is non-invertible only on the remaining
one transversal $Q\in S_Y(321)$, namely, $Q=(31524)$
(cf. Fig.~\ref{Example 1}, where all first subsequences are denoted by
$T^1$.) Its preimage is $\phi^{-1}(Q)=\{T_1,T_2\}$ where $T_1=(31542)$
and $T_2=(32514)$. Note that $(T_1)^1=(T_2)^1$($=\{3,5\}$), which
ensures that $\phi(T_1)=\phi(T_2)$($=Q$). Yet, the canonical
decompositions of $T_1$ and $T_2$ into connected components differ:
$T_1=T^1\sqcup \{1\} \sqcup \{4,2\}$ and $T_2=T^1\sqcup \{2,1\} \sqcup
\{4\}$, causing two preimages of $Q$.}
\label{example1}
\end{ex}
\begin{figure}[h]
\labellist
\small\hair 2pt
\pinlabel $\scriptstyle{T_1}$ at -8 522
\pinlabel $\scriptstyle{T_2}$ at 190 522
\pinlabel ${\text{and}}$ at 128 583
\pinlabel ${\phi}$ at 329 586
\pinlabel $\scriptstyle{Q}$ at 397 523
\endlabellist
\begin{center}
\includegraphics[width=4in]{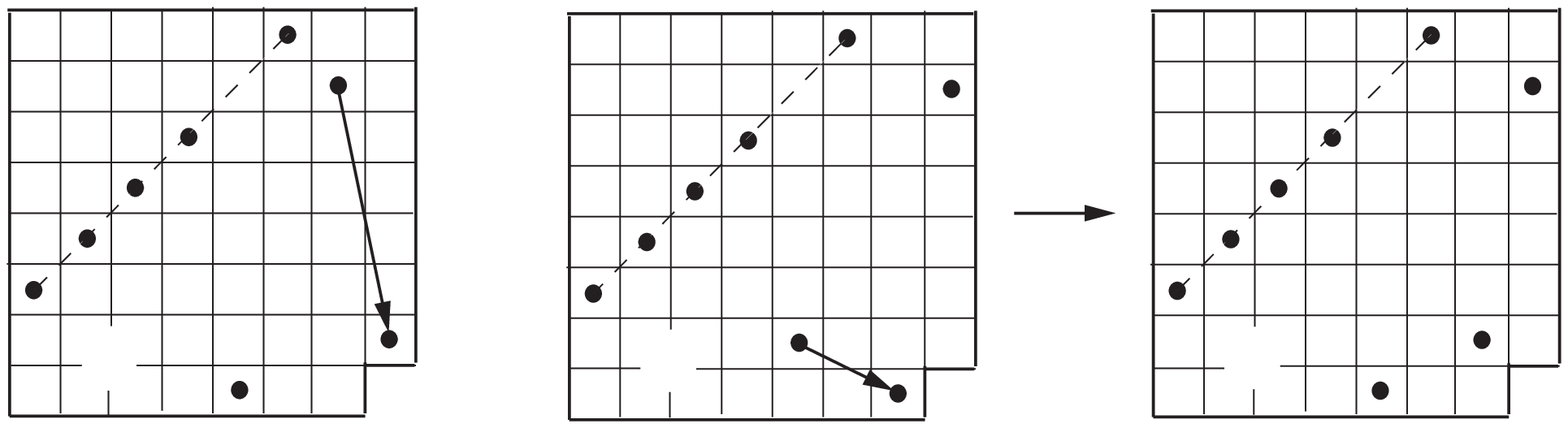}
\caption{$T_1,\,\,T_2\in
  S_{Y_n}(312)\,\,\stackrel{\phi}{\longrightarrow}\,\, Q\in
  S_{Y_n}(321)$}
\label{Example 2}
\end{center} 
\end{figure}

\begin{ex}{\rm We extend Example 1 to all $Y_n$ with $n\geq 5$. 
Let $T_1^n=(3,4,...,n-2,1,n,n-1,2)$ and
$T_2^n=(3,4,...,n-2,2,n,1,n-1)$ (cf. Fig.~\ref{Example 2}.) It is easy
to verify that $T_1^n$ and $T_2^n$ are $(312)$-avoiding on $Y_n$ with
the same first subsequence $(T_1^n)^1=(T_2^n)^1=(3,4,...,n-2,n)$, and
as such, they have the same image
$Q=\phi(T_1^n)=\phi(T_2^n)=(3,4,...,n-2,1,n,2,n-1)\in S_{Y_n}(321)$.
Hence $|S_{Y_n}(312)|>|S_{Y_n}(321)|$. Non-surprisingly, reducing
$Y_n$ along most of the first subsequence:
$Y_n\big/_{\displaystyle\{3,4,...,n-3\}}$, we recover the permutations
in $Y_5$ of Example 1.}
\label{Example Y_n}
\end{ex}

\subsection{Sufficient condition for strict inequality}
\begin{prop} If $Y$ has an $i$-critical point with $i\geq 3$, then 
$|S_Y(312)|>|S_Y(321)|$.
\label{higher critical points}
\end{prop}

\begin{figure}[h]
\labellist
\small\hair 2pt
\pinlabel $\scriptstyle{Y(P)}$ at -2 616
\pinlabel $\scriptstyle{S_1}$ at -15 530
\pinlabel $\scriptstyle{S_2}$ at 72 620
\pinlabel $\scriptstyle{T_1^6\,\,\text{or}\,\,T_2^6}$ at 33 570
\pinlabel $\scriptstyle{P}$ at 72 530
\pinlabel $\scriptstyle{Y_P}$ at 293 653
\pinlabel $\scriptstyle{\alpha=\gamma}$ at 247 638
\pinlabel $\scriptstyle{\beta=\delta}$ at 262 582
\pinlabel $\scriptstyle{Q}$ at 264 548
\pinlabel $\scriptstyle{R}$ at 316 601
\pinlabel $\scriptstyle{P}$ at 316 549
\pinlabel $\scriptstyle{\phantom{}_QY}$ at 353 583
\pinlabel $\scriptstyle{Y^R}$ at 234 515
\pinlabel $\zeta_P$ at 394 546
\pinlabel $\scriptstyle{\alpha_R}$ at 457 567
\pinlabel $\scriptstyle{\beta_R}$ at 474 546
\pinlabel $\scriptstyle{\alpha_Q}$ at 539 667
\pinlabel $\scriptstyle{\beta_Q}$ at 564 596
\pinlabel $\scriptstyle{\phantom{}_QY}$ at 605 601
\pinlabel $\scriptstyle{Y^R}$ at 449 497
\endlabellist
\begin{center}
\includegraphics[width=5in]{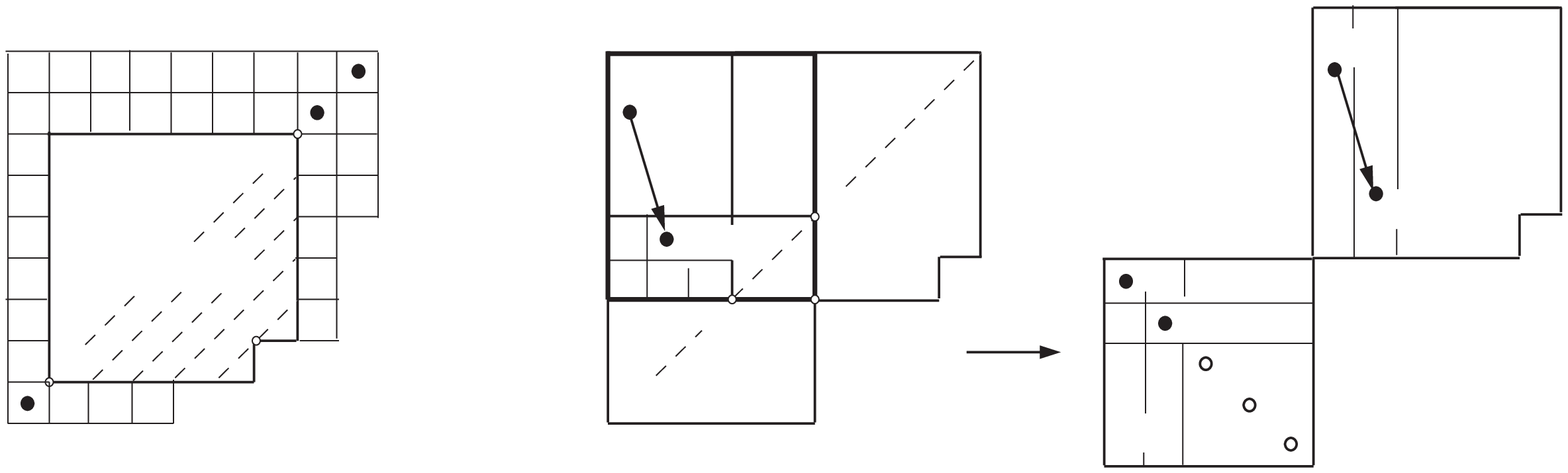}
\caption{(a) $\bar{T}_j=(1,T_j^{i+2},8,9)$ \hspace*{20mm} (b)-(c)
$(\gamma\delta)\subset Y_P$}
\label{Sufficient Condition 1}
\end{center} 
\end{figure}

\noindent{\sc Proof:} As in Example~\ref{Example Y_n}, for strict
inequality it is necessary and sufficient to exhibit two distinct
transversals $\bar{T}_1,\bar{T}_2\in S_Y(312)$ with
$(\bar{T}_1)^1=(\bar{T}_2)^1$. Let $P$ be an $i$-critical point of $Y$
with $i\geq 3$. Starting from $P$, go down (resp. right) one cell and go
left (resp. up) till hitting $d_0(Y)$: call this point $S_1$
(resp. $S_2$). With $S_1S_2$ as diagonal, we construct a subdiagram $Y(P)$
of $Y$ such that $Y(P)\cong Y_{i+2}$ and $P$ is the $i$-critical point
of $Y(P)$. For example, in Figure~\ref{Sufficient Condition 1}a the
subdiagram $Y(P)\cong Y_6$ is generated by the $4$-critical point $P$;
the dashed lines represent the diagonals $d_i(Y_6)$ for $0\leq i \leq
4$.

Now, put dots everywhere along $d(Y)$ outside of $Y(P)$. (In
Fig.~\ref{Sufficient Condition 1}a, these dots represent 1, 8 and 9.)
For $j=1,2$, insert $T_j^{i+2}$ from Example~\ref{Example Y_n} inside
$Y(P)$ in order to obtain $\bar{T}_j$ on $Y$. It is immediate that
$\bar{T}_1,\bar{T}_2\in S_Y(312)$ and they have the same first
subsequence, so that $\phi(\bar{T}_1)=\phi(\bar{T}_2)$, and hence
$|S_Y(312)|>|S_Y(321)|$. \qed

\subsection{Necessary condition for strict inequality $S_Y(312)>S_Y(321)$}
We shall prove that strict inequalities are obtained, as
Theorem~\ref{summary 312>321>213} claims, only when $Y$ has higher
critical points. To this end, we first need to establish two technical
recursive formulas for $2$-critical points when the avoided pattern is
$\sigma=(312)$ or $(321)$.

\subsubsection{Recursions for $2$-critical points} Recall the
points $R$ and $Q$ associated to $P$ in the definition of the map
$\zeta_P$. When $P$ is the {\it bottom} critical point of $Y$, $Y^R$
and $Y^Q$ are both {\it squares}, which makes the calculations below
possible (cf. Fig.~\ref{Sufficient Condition 1}b.)
Recursion~(\ref{|SY(sigma)|-recursion}) in Lemma~\ref{recursions} below
reduces calculations from the larger Young diagram $Y$ to the smaller
$\phantom{}_QY$; yet, it is not very useful on its own since it also
introduces the new sets
$S_{{\scriptscriptstyle{\searrow}}\phantom{}_QY}(\sigma)$ and
$S_{{\scriptscriptstyle{\nearrow}}\phantom{}_QY}(\sigma)$. Hence the
necessity to prove recursion~(\ref{|SY'(sigma)|-recursion}). Note the
apparent similarity between these recursive formulas for $|S_Y(\sigma)|$
and $|S_{{\scriptscriptstyle{\searrow}}Y}(\sigma)|$.

\begin{lem}
\label{recursions}
Let $Y$ be a Young diagram whose bottom critical point $P$ is
$2$-critical. If there are $k$ rows of $Y$ below $P$, for
$\sigma=(312)$ or $(321)$ we have:
\begin{eqnarray}
\quad |S_Y(\sigma)|&=&
c_{k+1}\cdot|S_{{\scriptscriptstyle{\searrow}}\phantom{}_QY}(\sigma)|+
(c_{k+2}-c_{k+1})\cdot|S_{{\scriptscriptstyle{\nearrow}}\phantom{}_QY}(\sigma)|
\label{|SY(sigma)|-recursion}\\
|S_{{\scriptscriptstyle{\searrow}}Y}(\sigma)|&=&
\quad c_k\cdot|S_{{\scriptscriptstyle{\searrow}}\phantom{}_QY}(\sigma)|+
\quad (c_{k+1}-c_{k})\cdot|S_{{\scriptscriptstyle{\nearrow}}\phantom{}_QY}(\sigma)|\label{|SY'(sigma)|-recursion}
\end{eqnarray}
\end{lem}

\noindent{\sc Proof:} The $2$-critical splitting from
Lemma~\ref{General 2-splitting lemma} implies:
\begin{equation*}
\label{summation for 2-splitting}
|S_Y(\sigma)|=|S^{\scriptscriptstyle{\nearrow}}_{Y^R}(\sigma)|\cdot
|S_{{\scriptscriptstyle{\nearrow}}\phantom{}_QY}(\sigma)|+
|S^{\scriptscriptstyle{\searrow}}_{Y^R}(\sigma)|\cdot
|S_{{\scriptscriptstyle{\searrow}}\phantom{}_QY}(\sigma)|.
\end{equation*}
Claim~\ref{avoidance in rectangles} below treats the special case of
the square $Y^R$ of size $k+2$. Substituting its results
$|S^{\scriptscriptstyle{\searrow}}_{Y^R}(\sigma)|=c_{k+1}$ and
$|S^{\scriptscriptstyle{\nearrow}}_{Y^R}(\sigma)|=c_{k+2}-c_{k+1}$, we
readily arrive at the wanted
recursion~(\ref{|SY(sigma)|-recursion}). 

\smallskip
To prove (\ref{|SY'(sigma)|-recursion}), we restrict the map $\zeta_P$ to
$S_{{\scriptscriptstyle{\searrow}}Y}(\sigma)$ in the 2-splitting
isomorphism in (\ref{general 2-splitting}):
\begin{equation}
\label{images}
\zeta_P(S_{{\scriptscriptstyle{\searrow}}Y}(\sigma))\subset
S^{\,\,\,\,\,\,\scriptscriptstyle{\nearrow}}_{{\scriptscriptstyle{\searrow}}
Y^R}(\sigma)\times
S_{{\scriptscriptstyle{\nearrow}}\phantom{}_QY}(\sigma) \sqcup
S^{\,\,\,\,\,\,\scriptscriptstyle{\searrow}}_{{\scriptscriptstyle{\searrow}}
Y^R}(\sigma)\times
S_{{\scriptscriptstyle{\searrow}}\phantom{}_QY}(\sigma).
\end{equation}
As in the definition of $\zeta_P$, we write $(\alpha\beta)$ for the
2-element subsequence of $T$ inside $Y_P$. There are three
possibilities for the initial decreasing subsequence $(\gamma\delta)$
of $T\in S_{{\scriptscriptstyle{\searrow}}Y}(\sigma)$.

\smallskip
{\it Case 1. $(\gamma\delta)\!\!\searrow$ is entirely in the rectangle
$Y_P$.} Then $(\gamma\delta)=(\alpha\beta)\!\!\searrow$ and
\[\zeta_P(T)\in 
S^{\,\,\,\,\,\,\scriptscriptstyle{\searrow}}_{{\scriptscriptstyle{\searrow}}
Y^R}(\sigma)\times
S_{{\scriptscriptstyle{\searrow}}\phantom{}_QY}(\sigma),\] with the
extra condition that $\alpha_R$ occupies cell (1,1) and $\beta_R$
occupies cell (2,2) of square $Y^R$ (cf. Fig.~\ref{Sufficient
Condition 1}b-c.)  If avoiding $\sigma=(312)$, the remainder of the
transversal in $Y^R$ is completely determined as a decreasing
subsequence (depicted in Fig.~\ref{Sufficient Condition 1}c via
``$\circ$''), while avoiding $(321)$ yields no possible completions in
$Y^R$. Thus, the images $\zeta_P(T)$ are in 1-1 correspondence with
$\{J_{k+2}\}\times
S_{{\scriptscriptstyle{\searrow}}\phantom{}_QY}(\sigma)$ if
$\sigma=(312)$, and there are 0 such if $\sigma=(321)$.

\begin{figure}[h]
\labellist
\small\hair 2pt
\pinlabel $\scriptstyle{Y_P}$ at -74 708
\pinlabel $\scriptstyle{P}$ at 10 602
\pinlabel $\scriptstyle{\phantom{}_QY}$ at 48 638
\pinlabel $\scriptstyle{\alpha}$ at -53 633
\pinlabel $\scriptstyle{\beta}$ at -9 704
\pinlabel $\scriptstyle{\gamma}$ at -83 591
\pinlabel $\scriptstyle{\delta}$ at -60 577 
\pinlabel $\scriptstyle{Q}$ at -44 602
\pinlabel $\scriptstyle{R}$ at 10 655
\pinlabel $\scriptstyle{Y^R}$ at -13 568
\pinlabel $\zeta_P$ at 85 601
\pinlabel $\scriptstyle{\alpha_R}$ at 182 593
\pinlabel $\scriptstyle{\beta_R}$ at 197 628
\pinlabel $\scriptstyle{\gamma}$ at 132.5 572
\pinlabel $\scriptstyle{\delta}$ at 154.5 561 
\pinlabel $\scriptstyle{\beta_Q}$ at 251 721
\pinlabel $\scriptstyle{\alpha_Q}$ at 236 648
\pinlabel $\scriptstyle{\phantom{}_QY}$ at 299 655
\pinlabel $\scriptstyle{Y^R}$ at 204 550
\endlabellist
\begin{center}
\includegraphics[width=3in]{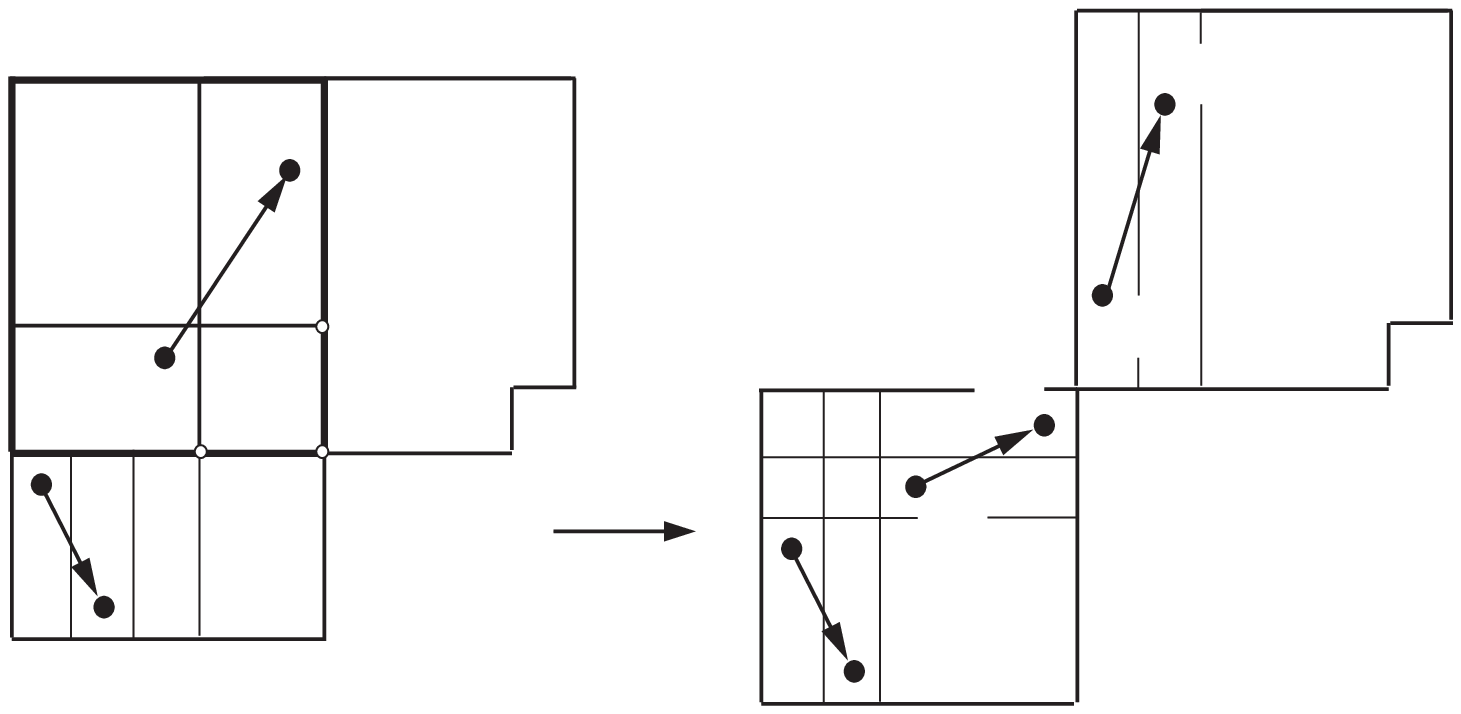}
\caption{Case 2}
\label{Case 2}
\end{center} 
\end{figure}

\smallskip
{\it Case 2. $(\gamma\delta)\!\!\searrow$ is entirely in the square
$Y^Q$} (cf. Fig.~\ref{Case 2}.) Then $(\gamma\delta)\cap
(\alpha\beta)=\emptyset$, and hence $(\alpha\beta)$ can be $\nearrow$
or $\searrow$. In either case, the four elements
$(\gamma\delta\alpha_R\beta_R)$ occupy the two leftmost columns and
two top rows of $Y^R$.  In the sub-factor
$S^{\,\,\,\,\,\,\scriptscriptstyle{\nearrow}}_{{\scriptscriptstyle{\searrow}}
Y^R}(\sigma)$ of (\ref{images}), $(\gamma\delta\alpha_R\beta_R)\approx
(2134)$, while in the sub-factor
$S^{\,\,\,\,\,\,\scriptscriptstyle{\searrow}}_{{\scriptscriptstyle{\searrow}}
Y^R}(\sigma)$, $(\gamma\delta\alpha_R\beta_R)\approx (2143)$.
Claim~\ref{more special avoidance in rectangles}a-b implies that the
number of images $\zeta_P(T)$ in these two subcases equals
respectively $(c_{k+1}-c_k-k)\cdot
|S_{{\scriptscriptstyle{\nearrow}}\phantom{}_QY}(\sigma)|$ or
$(c_k-1)\cdot
|S_{{\scriptscriptstyle{\searrow}}\phantom{}_QY}(\sigma)|$.

\begin{figure}[h]
\labellist
\small\hair 2pt
\pinlabel $\scriptstyle{Y_P}$ at 230 708
\pinlabel $\scriptstyle{P}$ at 315 602
\pinlabel $\scriptstyle{\phantom{}_QY}$ at 352 638
\pinlabel $\scriptstyle{\gamma=\alpha}$ at 233 684
\pinlabel $\scriptstyle{\beta}$ at 298 630
\pinlabel $\scriptstyle{\delta}$ at 242 573 
\pinlabel $\scriptstyle{Q}$ at 261 601
\pinlabel $\scriptstyle{R}$ at 315 655
\pinlabel $\scriptstyle{Y^R}$ at 293 569
\pinlabel $\zeta_P$ at 392 601
\pinlabel $\scriptstyle{\gamma=\alpha_R}$ at 415 619
\pinlabel $\scriptstyle{\beta_R}$ at 525 592
\pinlabel $\scriptstyle{\delta}$ at 460.5 560 
\pinlabel $\scriptstyle{\beta_Q}$ at 558 624
\pinlabel $\scriptstyle{\alpha_Q}$ at 536 692
\pinlabel $\scriptstyle{\phantom{}_QY}$ at 604 655
\pinlabel $\scriptstyle{Y^R}$ at 508 550
\endlabellist
\begin{center}
\includegraphics[width=3in]{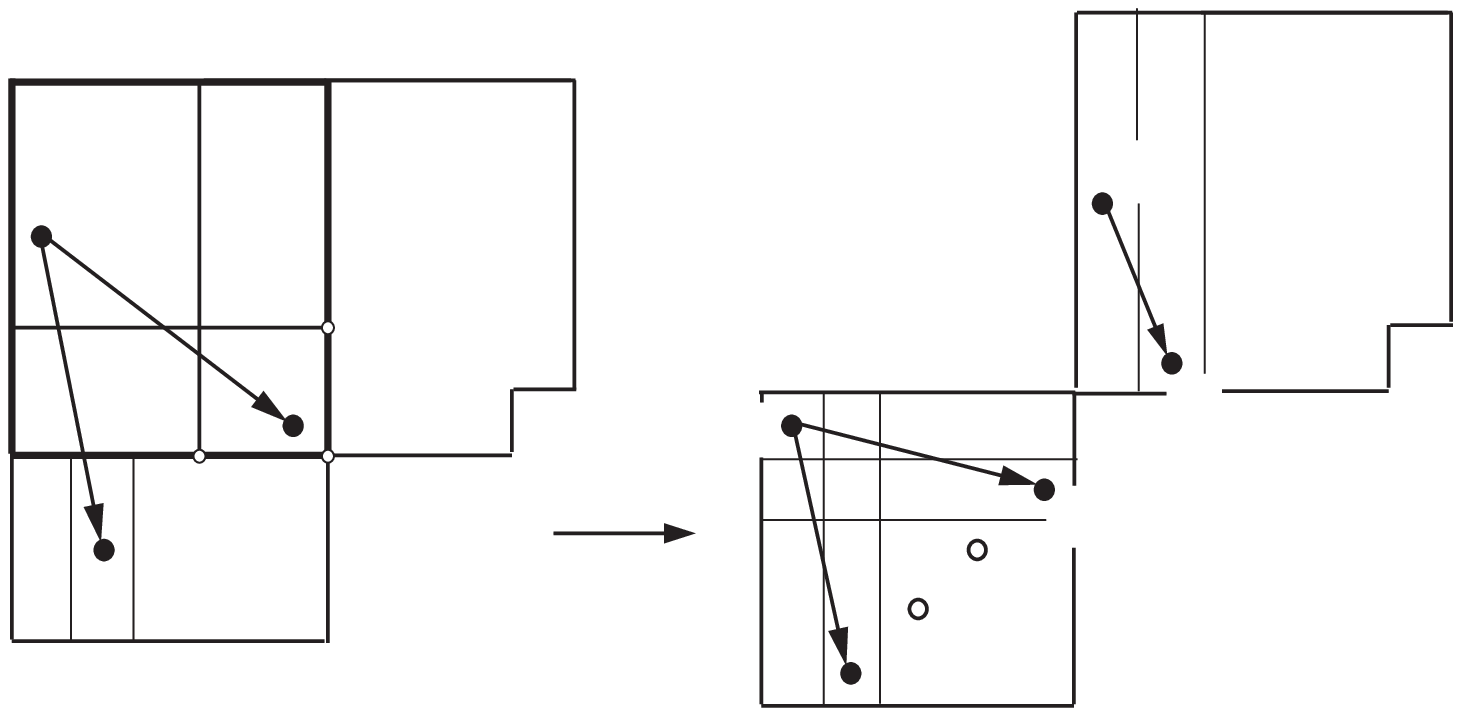}
\caption{Case 3}
\label{Case 3}
\end{center} 
\end{figure}

\smallskip
{\it Case 3. $\gamma\in Y_P$ and $\delta\in Y_Q$} (cf. Fig.~\ref{Case
3}.) Then $\gamma=\alpha$, and $(\alpha\beta)$ can be $\nearrow$ or
$\searrow$.  In the sub-factor
$S^{\,\,\,\,\,\,\scriptscriptstyle{\nearrow}}_{{\scriptscriptstyle{\searrow}}
Y^R}(\sigma)$, we have $(\alpha_R\delta\beta_R)\approx (213)$ where
$\alpha_R$ occupies cell (2,1). Claim~\ref{more special avoidance in
rectangles}c implies that the number of images in this subcase is
$k\cdot S_{{\scriptscriptstyle{\nearrow}}\phantom{}_QY}(\sigma)$. In
the sub-factor
$S^{\,\,\,\,\,\,\scriptscriptstyle{\searrow}}_{{\scriptscriptstyle{\searrow}}
Y^R}(\sigma)$, $(\alpha_R\delta\beta_R)\approx (312)$ where $\alpha_R$
occupies cell (1,1). Thus, avoiding $\sigma=(312)$ yields 0
transversals in this subcase. For $\sigma=(321)$, the position of
$\alpha_R$ allows for only one transversal on $Y^R$, namely,
$T_1=(k+2,1,2,...,k+1)$ (depicted in Fig.~\ref{Case 3}b via
``$\circ$''), and hence the images $\zeta_P(T)$ here are in 1-1
correspondence with $\{T_1\}\times
S_{{\scriptscriptstyle{\searrow}}\phantom{}_QY}(\sigma)$.

\medskip
Adding up the results in all three Cases, we obtain for each
$\sigma=(312)$ and $\sigma=(321)$:
\[|\zeta_P(S_{{\scriptscriptstyle{\searrow}}Y}(\sigma))|=
(0+1+c_k-1)\cdot 
|S_{{\scriptscriptstyle{\searrow}}\phantom{}_QY}(\sigma)|+
(c_{k+1}-c_k-k+k)\cdot
|S_{{\scriptscriptstyle{\nearrow}}\phantom{}_QY}(\sigma)|
\]
Since $\zeta_P$ is injective, we derive the wanted recursion
(\ref{|SY'(sigma)|-recursion}):
\[|S_{{\scriptscriptstyle{\searrow}}Y}(\sigma)|=|\zeta_P(S_{{\scriptscriptstyle{\searrow}}Y}(\sigma))|=c_k\cdot 
|S_{{\scriptscriptstyle{\searrow}}\phantom{}_QY}(\sigma)|+
(c_{k+1}-c_k)\cdot
|S_{{\scriptscriptstyle{\nearrow}}\phantom{}_QY}(\sigma)|.  \qed\]

\subsubsection{Calculations on the square $Y^R$} We show here all Claims
from the proof of Lemma~\ref{recursions}: they involve specific
calculations on the square $Y^R$ of size $k+2$. To simplify notation,
we shall write $\alpha$ for $\alpha_R$ and $\beta$ for
$\beta_R$. Thus, $(\alpha\beta)$ and $(\gamma\delta)$ are the
subsequences of $T$ in the top two rows, respectively leftmost two
columns, of $Y^R$.

\begin{figure}[h]
\labellist
\small\hair 2pt
\pinlabel $\scriptstyle{\alpha}$ at 431 710
\pinlabel $\scriptstyle{\beta}$ at 504 679
\pinlabel $Y^R$ at 440 642
\pinlabel $\scriptstyle{Y^R\!\!\big/_{\!\!\{\beta\}}}$ at 544 689
\pinlabel $\scriptstyle{Y^R\!\!\big/_{\!\!\{\beta\}}}$ at 680 689
\pinlabel $\scriptstyle{\alpha}$ at 756 700
\pinlabel $\scriptstyle{\beta}$ at 775 679
\pinlabel $M_{k+1}$ at 612 671
\pinlabel $Y^R$ at 730 642
\endlabellist
\begin{center}
\includegraphics[width=3in]{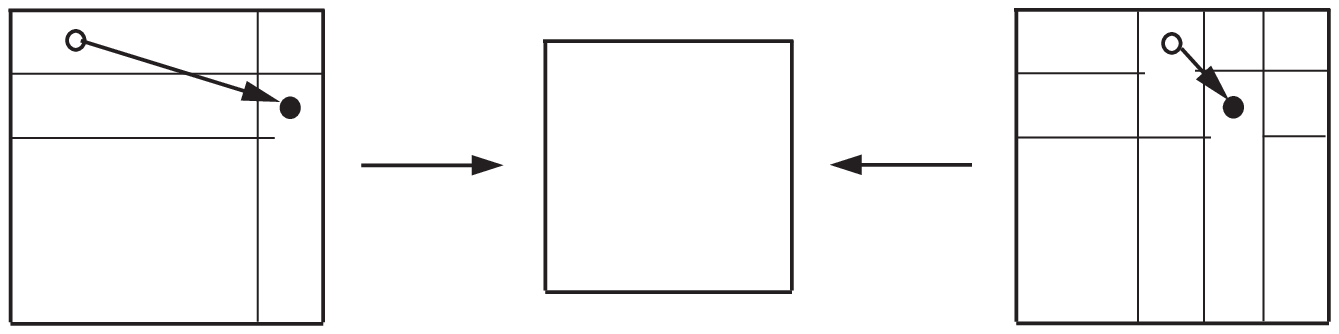}
\caption{Claim 1}
\label{Claim 1}
\end{center} 
\end{figure}

\begin{claim}
For $\sigma=(312)$ or $(321)$,
$|S^{\scriptscriptstyle{\searrow}}_{Y^R}(\sigma)|=c_{k+1}$, and hence 
$|S^{\scriptscriptstyle{\nearrow}}_{Y^R}(\sigma)|=c_{k+2}-c_{k+1}$.
\label{avoidance in rectangles}
\end{claim}
\noindent{\sc Proof:} We calculate first
$|S^{\scriptscriptstyle{\searrow}}_{Y^R}(\sigma)|$, so we assume that
$(\alpha\beta)\!\!\searrow$ in $Y^R$. If $\sigma=(321)$, to avoid the
situation of $\alpha$ and $\beta$ playing the roles of ``3'' and ``2''
in a $(321)$-pattern in $Y^R$, $\beta$ must be in the last column (and
the second row) of $Y^R$. As such, $\beta$ cannot participate in any
$(321)$-pattern on $Y^R$, so that reducing along $\beta$ we obtain a
$(321)$-avoiding transversal $T^{\prime}$ on the rectangle
$Y^R\!\!\big/_{\!\!\{\beta\}}=M_{k+1}$, without any further
restrictions (cf. Fig.~\ref{Claim 1}a-b.) The original transversal of
$Y^R$ can be reconstructed from $T^{\prime}$ by reinserting $\beta$ in
the last row and second column of $Y^R$.  We have established that
$S^{\scriptscriptstyle{\searrow}}_{Y^R}(321)\cong S_{k+1}(321)$, and
hence $|S^{\scriptscriptstyle{\searrow}}_{Y^R}(321)|=c_{k+1}$.

Similarly, if $\sigma=(312)$, $\alpha$ and $\beta$ must be in adjacent
columns in $Y^R$ in order to avoid $(312)$. But they are already in
the top two rows of $Y$, so they are situated in diagonally-adjacent
cells. Again, reducing along $\beta$ we obtain a $(312)$-avoiding
transversal $T^{\prime\prime}$ on the rectangle
$Y^R\!\!\big/_{\!\!\{\beta\}}=M_{k+1}$, without any further
restrictions (cf. Fig.~\ref{Claim 1}c-b.) The original transversal of
$Y^R$ can be reconstructed from $T^{\prime\prime}$ by reinserting
$\beta$ in the second row and the column on the right of
$\alpha$'s column in $Y^R$.  We have established that
$S^{\scriptscriptstyle{\searrow}}_{Y^R}(312)\cong S_{k+1}(312)$, and
hence $|S^{\scriptscriptstyle{\searrow}}_{Y^R}(312)|=c_{k+1}$.

To finish the argument, we note that
$S^{\scriptscriptstyle{\nearrow}}_{Y^R}(\sigma)$ is the complement of
$S^{\scriptscriptstyle{\searrow}}_{Y^R}(\sigma)$ in
${S}_{Y^R}(\sigma)$, so that for $\sigma=(321)$ of $(312)$ we have
$|S^{\scriptscriptstyle{\nearrow}}_{Y^R}(\sigma)|=
|S_{Y^R}(\sigma)|-|S^{\scriptscriptstyle{\searrow}}_{Y^R}(\sigma)|=
c_{k+2}-c_{k+1}.$
\qed

\begin{claim}
\label{more special avoidance in rectangles}
(a) There are $c_{k+1}-c_k-k$ transversals $T\in
S^{\,\,\,\,\,\,\scriptscriptstyle{\nearrow}}_{{\scriptscriptstyle{\searrow}}
Y^R}(\sigma)$ with $\{\gamma\delta\}\cap \{\alpha\beta\}=\emptyset$.

\hspace*{12.5mm} (b)  There are $c_k-1$ transversals $T\in
S^{\,\,\,\,\,\,\scriptscriptstyle{\searrow}}_{{\scriptscriptstyle{\searrow}}
Y^R}(\sigma)$ with $\{\gamma\delta\}\cap \{\alpha\beta\}=\emptyset$.

\hspace*{12.5mm} (c) There are $k$ transversals $T\in
S^{\,\,\,\,\,\,\scriptscriptstyle{\nearrow}}_{{\scriptscriptstyle{\searrow}}
Y^R}(\sigma)$ with $\alpha=\gamma$.
\end{claim}

\noindent
{\sc Proof:} All transversals in question are elements of
$S_{{\scriptscriptstyle{\searrow}} Y^R}(\sigma)$,
i.e. $(\gamma\delta)\!\!\searrow$. Similarly to the proof of
Claim~\ref{avoidance in rectangles}, either $\delta$ is in the last
row (and second column) of $Y^R$ ($\sigma=(321)$), or neighboring
$\gamma$ southeast-diagonally ($\sigma=(312)$). In either
case, we reduce $Y^R$ along $\delta$ to obtain equinumerant
subsets of $Y^R\!\!\big/_{\!\!\{\delta\}}=M_{k+1}$ with
the following restrictions (cf. Fig.~\ref{Claim 2}):
\begin{figure}[h]
\labellist
\small\hair 2pt
\pinlabel $\scriptstyle{\alpha}$ at 200 691
\pinlabel $\scriptstyle{\beta}$ at 250 708
\pinlabel $\scriptstyle{\alpha}$ at 318 710
\pinlabel $\scriptstyle{\beta}$ at 366 689
\pinlabel $M_{k+1}$ at 235 640
\pinlabel $M_{k+1}$ at 353 640
\pinlabel $M_{k+1}$ at 435 655
\pinlabel $\scriptstyle{\alpha}$ at 416 681
\pinlabel $\scriptstyle{\beta}$ at 457 708
\pinlabel $\stackrel{\scriptscriptstyle{\nearrow}}{M_{k+1}}=$ at 596 677
\pinlabel $\cong$ at 743 673
\pinlabel $={\scriptscriptstyle{\nearrow}}M_{k+1}$ at 
 903 676
\endlabellist
\hspace*{-14mm}\includegraphics[width=5in]{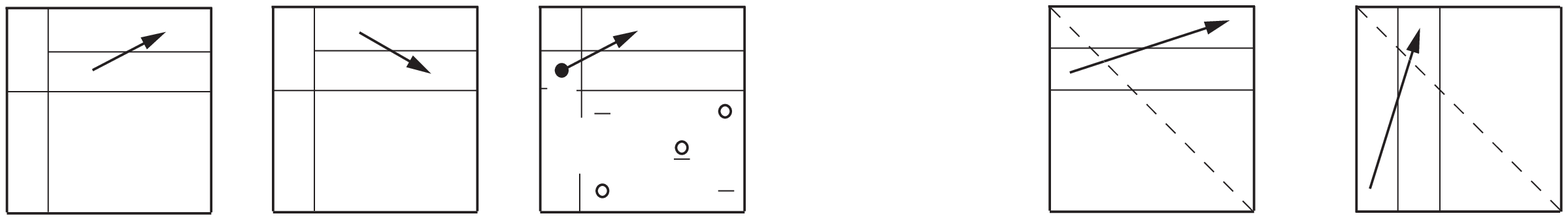}
\caption{Claim 2. (a1)-(b1)-(c1)\hspace*{27mm}
$(\stackrel{\scriptscriptstyle{\nearrow}}{M}_{k+1})^t\protect\cong
{\scriptscriptstyle{\nearrow}}M_{k+1}$\hspace*{20mm}}
\label{Claim 2}
\end{figure}

\begin{itemize}
\item[(a1)] all transversals of
$S^{\scriptscriptstyle{\nearrow}}_{M_{k+1}}(\sigma)$ for which the top
two elements $(\alpha\beta)\!\!\nearrow$ do not lie in the first
column of $M_{k+1}$ (occupied by $\gamma$);

\item[(b1)] all transversals of
$S^{\scriptscriptstyle{\searrow}}_{M_{k+1}}(\sigma)$ for which the top
two elements $(\alpha\beta)\!\!\searrow$ do not lie in the first
column of $M_{k+1}$ (occupied by $\gamma$);

\item[(c1)] all transversals of
$S^{\,\,\,\,\,\,\scriptscriptstyle{\nearrow}}_{M_{k+1}}(\sigma)$, for
which one of the top two elements $(\alpha\beta)\!\!\nearrow$ does lie
in the first column of $M_{k+1}$ ($\alpha=\gamma$ is that element.) 
\end{itemize}

Let's start with case (c1). Since $\alpha$ is in position (2,1), the
rows below $\alpha$ are filled either with an increasing (for
$\sigma=(321)$) or with a decreasing (for $\sigma=(312)$) subsequence.
In Fig.~\ref{Claim 2}c, $\circ$ and $-$ denote, respectively, these
increasing and decreasing subsequences.  At the same time, $\beta$
can be reinserted in any of the $k$ possible cells of the top row of
$Y^R$ without creating any $\sigma$-patterns. Thus, the number of
transversals in (c) is $k$. \qed

Case (a1) is the complement of (c1) inside
$S^{\,\,\,\,\,\,\scriptscriptstyle{\nearrow}}_{M_{k+1}}(\sigma)$.
Since $(321)$, $(312)$ and $Y^R$ are symmetric with respect to
transposing across the northwest/southeast diagonal, and since
$(S^{\,\,\,\,\,\,\scriptscriptstyle{\nearrow}}_{M_{k+1}}(\sigma))^t=
S_{{\scriptscriptstyle{\nearrow}}M_{k+1}}(\sigma)$, we
can use Claim~\ref{avoidance in rectangles} for $Y^R=M_{k+1}$ to calculate:
\[|S^{\,\,\,\,\,\,\scriptscriptstyle{\nearrow}}_{M_{k+1}}(\sigma)|=
|S_{{\scriptscriptstyle{\nearrow}}M_{k+1}}(\sigma)|=c_{k+1}-c_k.\]
Therefore, the number of transversals in (a) equals
$c_{k+1}-c_k-k$. \qed

Finally, case (b1) misses only one transversal of the set
$S^{\,\,\,\,\,\,\scriptscriptstyle{\searrow}}_{M_{k+1}}(\sigma)$:
namely, when $\alpha$ is in position (1,1), without any more
restrictions (cf. Fig~\ref{Claim 2}b.) In such a situation, the rest of
$M_{k+1}$ is again filled either with an increasing or with a
decreasing subsequence (respectively, for $\sigma=(321)$ and
$(312)$). Thus, case (b1) counts 1 fewer transversals than
$S^{\,\,\,\,\,\,\scriptscriptstyle{\searrow}}_{M_{k+1}}(\sigma)$.
Using again the transposing argument and Claim~\ref{avoidance in
rectangles} for $Y^R=M_{k+1}$, we conclude that the number of
transversals in (b) equals
\[|S^{\,\,\,\,\,\,\scriptscriptstyle{\searrow}}_{M_{k+1}}(\sigma)|-1=
|S_{{\scriptscriptstyle{\searrow}}M_{k+1}}(\sigma)|-1=c_k-1.\qed\]

\subsubsection{Conclusions for low-rank critical points}
\begin{lem}  
If $Y$ has only $2$-critical points, then 
$|S_{{\scriptscriptstyle{\searrow}}Y}(321)|=
|S_{{\scriptscriptstyle{\searrow}}Y}(312)|$, 
$|S_{{\scriptscriptstyle{\nearrow}}Y}(321)|=
|S_{{\scriptscriptstyle{\nearrow}}Y}(312)|$, and hence
$|S_{Y}(321)|=
|S_{Y}(312)|$.
\label{three identities}
\end{lem} 

\noindent{\sc Proof:} For the special case of a square $Y=M_{k+2}$
(which has no critical points), the first two equalities were proven
in Claim~\ref{avoidance in rectangles} for the square $Y_R$, while the
third equality is the well-known Wilf-equivalence $(312)\sim (321)$.

For the general case, we proceed by induction on the size $n$ of
$Y$. For $n\leq 3$ there are no $2$-critical points. Suppose that $Y$
is of size $n\geq 4$, not a square, and has only $2$-critical
points. Let $Y$'s bottom ($2$-)critical point be $P$. The Young
subdiagram $\phantom{}_QY$ from Lemma~\ref{recursions} is of smaller
size, and by construction, its critical points are all of $Y$'s
critical points, short of $P$.  Applying induction to $\phantom{}_QY$,
we have $|S_{{\scriptscriptstyle{\searrow}}\phantom{}_QY}(321)|=
|S_{{\scriptscriptstyle{\searrow}}\phantom{}_QY}(312)|$ and
$|S_{{\scriptscriptstyle{\nearrow}}\phantom{}_QY}(321)|=
|S_{{\scriptscriptstyle{\nearrow}}\phantom{}_QY}(312)|$.  Recursions
(\ref{|SY(sigma)|-recursion})-(\ref{|SY'(sigma)|-recursion}) then
imply $|S_Y(321)|=|S_Y(312)|$ and
$|S_{{\scriptscriptstyle{\searrow}}Y}(321)|=
|S_{{\scriptscriptstyle{\searrow}}Y}(312)|$. Since
$S_{{\scriptscriptstyle{\nearrow}}Y}(\sigma)$ is the complement of
$S_{{\scriptscriptstyle{\nearrow}}Y}(\sigma)$ in $S_Y(\sigma)$ for any
$\sigma$, it also follows that
$|S_{{\scriptscriptstyle{\nearrow}}Y}(321)|=
|S_{{\scriptscriptstyle{\nearrow}}Y}(312)|$. \qed

\begin{prop} 
$|S_Y(312)|=|S_Y(321)|$ if $Y$ has only $i$-critical points with $i\leq 2$.
\label{strict inequality 312>321}
\end{prop}

\noindent
{\sc Proof:} If $Y$ has some $0$- or $1$-critical point $P$,
Proposition~\ref{0-1-splitting} implies that there is a $0$- or
$1$-splitting for any permutation $\sigma$:
\[|S_Y(\sigma)|=|S_{U}(\sigma)|\cdot |S_{V}(\sigma)|,\]
where $U$ and $V$ are some Young subdiagrams of $Y$ of smaller sizes.
Since by construction the diagonals $d(U)$ and $d(V)$ lie on $d(Y)$,
the set of critical points of $U$ and $V$ is the same as the set of
critical points of $Y$, short of $P$. In other words, $U$ and $V$
again have only $0$-, $1$- or $2$-critical points.  Continuing the
splitting process for every $0$- or $1$-critical point of the smaller
diagrams, we arrive eventually at a splitting
\[|S_Y(\sigma)|=|S_{U_1}(\sigma)|\cdot |S_{U_2}(\sigma)|\cdots
|S_{U_k}(\sigma)|,\] where each subdiagram $U_i$ has only $2$-critical
points (or no critical points at all). Lemma~\ref{three identities}
guarantees that $|S_{U_i}(321)|=|S_{U_i}(321)|$ for all $i$, so that
the products $|S_Y(321)|$ and $|S_Y(312)|$ are also equal. \qed

\bigskip
Finally, combining the results of Propositions~\ref{higher critical
points}-\ref{strict inequality 312>321}, we derive the second 
necessary and sufficient condition in Theorem~\ref{summary
312>321>213}: $|S_Y(312)|>|S_Y(321)|$ if and only if
$Y$ contains an $i$-critical point with $i\geq 3$.\qed

\section{Proof of the Inequality $S_Y(213)\leq S_Y(123)$}
\label{SY(213) < SY(123)}

\subsection{The (213)-decomposition} In \cite{Stankova-West}, Stankova-West
show $(213)\sim_s (132)$. Their proof introduces a special
decomposition of the $(213)$-avoiding transversals on any Young
diagram. Here we modify and extend this decomposition for our
purposes, and use it later for comparing $S_Y(213)$ and $S_Y(123)$.

\begin{defn}
\label{213-decomposition}
{\rm Let $Y$ be a Young diagram and let $c$ be a cell in the bottom
row of $Y$. Start from the bottom left corner of $c$, draw a
$45^{\circ}$ ray in north-east direction until the ray intersects for
the first time the border of $Y$, and use the resulting segment as the
diagonal of a smaller subdiagram $\mathcal{A}_c$ of $Y$. Reducing $Y$
along $\mathcal{A}_c$ leaves a subdiagram
$\mathcal{B}_c=Y\!\!\big/_{\displaystyle{\!\mathcal{A}_c}}$.  Thus,
$c$ determines a pair $(\mathcal{A}_c,\mathcal{B}_c)$ of Young
subdiagrams of $Y$, called the {\it $(213)$-decomposition of $Y$
induced by $c$} and denoted by
$Y_{213}(c)=\mathcal{A}_c\otimes\mathcal{B}_c$. If a transversal $T\in
S_Y$ is concentrated in $\mathcal{A}_c$ and $\mathcal{B}_c$, we say
that $T$ respects this $(213)$-decomposition of $Y$ and we write
$T=T|_{\mathcal{A}_c}\otimes T|_{\mathcal{B}_c}$.}
\end{defn}
\begin{figure}[h]
\labellist
\small\hair 2pt
\pinlabel ${\mathcal{BL}_c}$ at 50 608
\pinlabel ${\mathcal{BR}_c}$ at 179 606
\pinlabel ${\mathcal{A}_c}$ at 102 518
\pinlabel ${\mathcal{BL}_c}$ at 363 608
\pinlabel ${\mathcal{BR}_c}$ at 406 608
\pinlabel ${\mathcal{B}_c}$ at 389 555
\pinlabel ${\mathcal{A}_c}$ at 275 518
\pinlabel ${\mathcal{B}_b}$ at 633 590
\pinlabel ${\mathcal{A}_b}$ at 500 518
\pinlabel ${\mathcal{BL}_b}$ at 614 620
\pinlabel ${\mathcal{BR}_b}$ at 648 631
\pinlabel $c$ at 81 441
\pinlabel $b$ at 63 476.5
\pinlabel $c$ at 252 441
\pinlabel $b$ at 478 476.5
\pinlabel $\text{or}$ at 235 640
\pinlabel $\alpha$ at 89 548
\pinlabel $\gamma$ at 171 573
\pinlabel $\text{or}$ at 438 555
\endlabellist
\begin{center}
\includegraphics[width=5in]{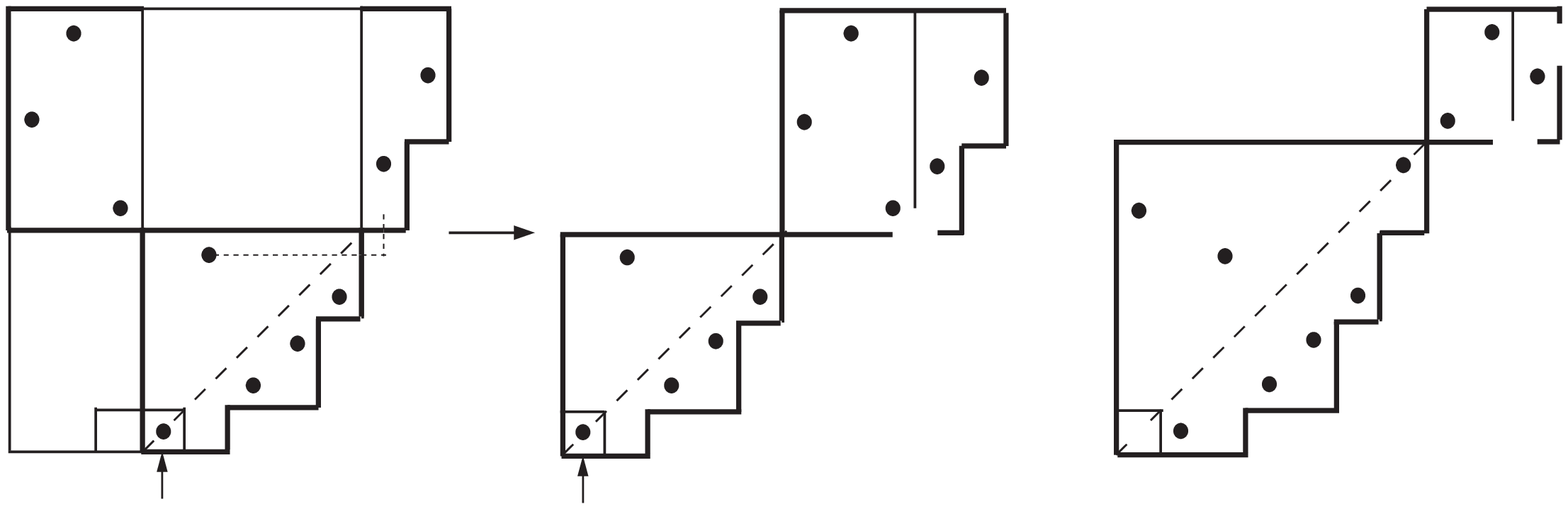}
\caption{Minimal and non-minimal $(213)$-decompositions}
\label{213-Decomposition}
\end{center} 
\end{figure}

Because of the $45^{\circ}$ angle of the diagonal $d(\mathcal{A}_c)$,
the Young subdiagram $\mathcal{A}_c$ is proper, and hence the
reduction along it, $\mathcal{B}_c$, is also proper. The smallest
$\mathcal{A}_c$ can be is the cell $c$: this happens when $c$ is the
rightmost cell in the bottom row of $Y$. The decomposition
$Y_{213}(c)$ is trivial exactly when $\mathcal{B}_c=\emptyset$: this
happens when $c$ is the bottom left corner of $Y$ and $Y$ has no
$0$-critical points. In such a case, $Y_{213}(c)=\mathcal{A}_c$.
While $\mathcal{A}_c$'s rows and columns are {\bf not} interspersed
with ``outside'' rows or columns from $Y\backslash \mathcal{A}_c$, in
general, $\mathcal{B}_c$ splits into two parts: $\mathcal{B}_c=
\mathcal{BR}_c+\mathcal{BL}_c$ where $\mathcal{BR}_c$ is to the left
and above $\mathcal{A}_c$ and $\mathcal{BL}_c$ is to the right and
above $\mathcal{A}_c$ (cf. Fig.~\ref{213-Decomposition}.)

\smallskip The name ``$(213)$-decomposition'' comes from the fact that
all $(213)$-avoiding transversals $T$ respect at least one
$(213)$-decomposition of $Y$.

\begin{prop}[Stankova-West]  Let $c$ be a bottom cell in $Y$, and let 
$T\in S_Y(213)$ have its bottom element in cell $c$.  Then $T$
respects the $(213)$-decomposition $Y_{213}(c)=\mathcal{A}_c\otimes
\mathcal{B}_c$, and hence it decomposes as
$T=T|_{\mathcal{A}_c}\otimes T|_{\mathcal{B}_c}$. Conversely, if $T\in
S_Y$ respects this $(213)$-decomposition, and if the restrictions
$T|_{\mathcal{A}_c}$ and $T|_{\mathcal{B}_c}$ are each
$(213)$-avoiding on $\mathcal{A}_c$ and $\mathcal{B}_c$, respectively,
then $T$ is $(213)$-avoiding on all of $Y$.
\label{213-avoiding $T$ respect 213-decomposition}
\end{prop}

For some $T\in S_Y(213)$, it is possible that different bottom cells
$c_i$'s induce different (213)-decompositions
$T=T|_{\mathcal{A}_{c_i}}\otimes T|_{\mathcal{B}_{c_i}}$. However, for
only one $(213)$-decomposition $T$'s bottom element $\hbar$ is in
$\mathcal{A}_{c}$'s bottom left corner; we call this {\it the minimal
$(213)$-decomposition of $T$} since all other $(213)$-decompositions
will have $\hbar$ somewhere further to the left and hence their
components $\mathcal{A}_{c_i}$ will contain properly
$\mathcal{A}_{c}$.  For example, Figure~\ref{213-Decomposition} shows
the minimal $T_{213}(c)=(15234)\otimes (35124)$ and a non-minimal
$T_{213}(b)=(6152347)\otimes (132)$ decompositions of
$T=(8,10,6,1,5,2,3,4,7,9)$. Proposition~\ref{213-avoiding $T$ respect
213-decomposition} implies that every $T\in S_Y(213)$ respects its
minimal $(213)$-decomposition. More generally,

\begin{defn}{\rm For any transversal $T\in S_Y$, 
the (213)-decomposition of $T$ whose $\mathcal{A}_c$-component
is contained properly in the $\mathcal{A}$-components of any other
$(213)$-decomposition of $T$ is called the {\it minimal
$(213)$-decomposition of }$T$. If all $(213)$-decompositions of $T$ are
trivial, i.e. $T=T|_{\mathcal{A}_c}$, we say that $T$ is {\it
$(213)$-indecomposable.} }
\end{defn} 

When it is irrelevant which bottom cell $c$ induces some
$(213)$-decomposition of $T$, we shall drop $c$ from the notation, e.g.
$T=T|_{\mathcal{A}}\otimes (T|_{\mathcal{BL}}+T|_{\mathcal{BR}})$.

For a general transversal $T\in S_Y(\sigma)$ which (213)-decomposes as
$T= T|_{\mathcal{A}}\otimes T|_{\mathcal{B}}$, it is evidently true
that $T|_{\mathcal{A}}$ and $T|_{\mathcal{B}}$ each avoid $\sigma$
on ${\mathcal{A}}$ and ${\mathcal{B}}$, respectively. The converse
is false in general: $T|_{\mathcal{A}}\times T|_{\mathcal{B}}\in
S_{\mathcal{A}}(\sigma)\times S_{\mathcal{B}}(\sigma)$ does not
imply $T\in S_Y(\sigma)$. Yet, for special cases of $\sigma$,
the converse is true. We have seen in Proposition \ref{213-avoiding
$T$ respect 213-decomposition} that $\sigma=(213)$ is such a special
case. Another one is $\sigma=(123)$.

\begin{lem}
If $T\in S_Y$ has a $(213)$-decomposition $T=T|_{\mathcal{A}}\otimes
T|_{\mathcal{B}}$ on $Y$, and if each subtransversal avoids $(123)$:
$T|_{\mathcal{A}}\in S_{\mathcal{A}}(123)$ and
$T|_{\mathcal{B}}\in S_{\mathcal{B}}(123)$, then $T\in S_Y(123)$.
\label{123-avoiding converse}
\end{lem}
 
\noindent{\sc Proof:} Suppose that $T$ has a $(123)$-subsequence
$(\alpha\beta\gamma)$ landing inside $Y$. Because of the hypothesis on
the two subtransversals, this pattern must involve elements from both
subdiagrams $\mathcal{A}$ and $\mathcal{B}$. Since $\mathcal{A}$ is
entirely below $\mathcal{B}$, the element $\alpha$ (which plays the
role of ``1'') must come from $\mathcal{A}$. But since $\alpha$ is
also the leftmost element of the pattern, it eliminates any
participation coming from $\mathcal{BL}$. This forces the last element
$\gamma$ to come from $\mathcal{BR}$
(cf. Fig.~\ref{213-Decomposition}a.) Yet, no two elements of
$\mathcal{A}$ and $\mathcal{BR}$ can participate in any pattern in
$Y$: by construction of the $45^{\circ}$ diagonal of $\mathcal{A}$,
$\mathcal{BR}$ is entirely to the right and above $\mathcal{A}$,
forcing the rows of $\mathcal{A}$ and the columns of $\mathcal{BR}$,
to intersect outside $Y$.

We conclude that a $(123)$-pattern is impossible, so that $T\in
S_Y(123)$. \qed

\bigskip
We shall see below that the $(213)$-decompositions of any $T\in S_Y$
on $Y$ are preserved by $(213)\!\rightarrow \!(123)$- and
$(123)\!\rightarrow \!(213)$ moves on $T$, which will allow us to
prove eventually the desired inequality $|S_Y(123)|\geq |S_Y(213)|$.

\subsection{Special $(213)$-decompositions}
Fix a $(213)$-decomposition $T=T|_{\mathcal{A}_c}\otimes
T|_{\mathcal{B}_c}$ of a transversal $T\in S_Y$. If
$\mathcal{BR}_c=\emptyset$, we obtain a generalization of the
decomposability Definition \ref{decomposable} of a permutation
$\sigma\in S_n$. To keep up with the previous conventions, we say in
this case that the transversal $T$ is {\it decomposable}, and also
write $T=(T|_{\mathcal{BL}_c}\big|T|_{\mathcal{A}_c})$. The two blocks
of $T$ are arranged in a northwest/southeast fashion.

\smallskip
On the other hand, if $\mathcal{BL}_c=\emptyset$, not only the given
transversal decomposes as $T=T|_{\mathcal{A}_c}\otimes
T|_{\mathcal{BR}_c}$, but {\bf any} transversals $T^{\prime}$ of $Y$
respects this decomposition. Indeed, in this case, $c$ is the bottom
left corner cell of $Y$ and $Y$ contains a $0$-critical point $P$, and
hence any transversal $T^{\prime}\in S_Y$ has this $0$-splitting with
respect to $P$. Here the two blocks of $T^{\prime}$ are arranged in a
southwest/northeast fashion. 

\smallskip Thus, all decompositions in this paper are
$(213)$-decompositions or special cases of it.

\subsection{$\sigma\!\rightarrow\!\tau$ moves on $T\in S_Y$} 
Since each $(213)\!\rightarrow\!(123)$ move decreases the number of
inversions in $T$, any sequence of $(213)\!\rightarrow\!  (123)$ moves
eventually {\it terminates} with some $T^{\prime}\in
S_Y(213)$. Similarly, a sequence of $(123)\!\rightarrow\!(213)$ moves
terminates with some $T^{\prime\prime}\in S_Y(123)$.

\begin{conj}
Starting with a transversal $T\in S_Y$, all sequences of
$(213)\rightarrow (123)$ moves terminate in the same transversal
$T^{\prime}\in S_Y(213)$.
\end{conj}

The conjecture, if proven, would give a well-defined map
$S_Y\rightarrow S_Y(213)$, which could be restricted to a map
$\xi:S_Y(123)\rightarrow S_Y(213)$. To show then that $\xi$ is
surjective, we would start with any $T^{\prime}\in S_Y(213)$ and apply
any sequence of $(123)\rightarrow (213)$ moves on $T^{\prime}$ until
it terminates with some $T\in S_Y(123)$. Reversing the sequence of
moves would yield a sequence of $(213)\rightarrow (123)$ moves on $T$
that terminates in $T^{\prime}$. The definition of $\xi$ would then
give $\xi(T)=T^{\prime}$ and hence $\xi$ would be surjective, from
where $|S_Y(123)|\geq |S_Y(213)|$.

\smallskip
In the next Subsection \ref{definition of phi:S_Y(213) to S_Y(123)} we
proceed in a different way by defining a {\it section} of the
conjectured map $\xi$, i.e. a map $\psi:S_Y(213)\rightarrow S_Y(123)$
such that $\xi\circ\psi=\text{id}_{S_Y(213)}$. As opposed to $\xi$,
the map $\psi$ will be given by a {\it specific} sequence of
$(123)\rightarrow (213)$ moves which can be retraced back. The latter
will then readily imply $|S_Y(213)|\leq |S_Y(123)|$.

Analogously to the {\it first} and {\it second} subsequences of $T\in
S_Y$ in Definition \ref{first and second subsequences}, in working with
$(123)$-avoidance we will need the following terminology.

\begin{defn}
{\rm Let $T\in S_Y$. The {\it primary} subsequence $\dot{T}$ of $T$
consists of all elements which are {\bf not} $(12)$-dominated in
$T$. The {\it secondary} subsequence $\ddot{T}$ of $T$ consists of all
elements of $T$ which are $(12)$-dominated by something in $\dot{T}$
and by nothing in $T\backslash \dot{T}$.}
\end{defn}

In particular, if $T\in S_Y(123)$ where $Y$ is a {\it square}, then
$T$ is the disjoint union of its two decreasing subsequences $\dot{T}$
and $\ddot{T}$: $\,\,T=\dot{T}\sqcup\ddot{T}$.

\subsection{Definition of $\psi:S_Y(213)\rightarrow S_Y(123)$}
\label{definition of phi:S_Y(213) to S_Y(123)}
We define $\psi$ by {\it induction on the size of $Y$}. When $Y$ is a
single cell, $\psi$ is the identity map. Suppose we have defined
$\psi$ for all Young diagrams of size $<n$. Fix a Young diagram $Y$ of
size $n\geq 2$ and $T\in S_Y(213)$. Throughout this subsection, we
shall refer to the element in the bottom row of $Y$ as $\hbar$ and
denote by $Y\!/\!_{\hbar}$ the reduction along it. There are two cases
to discuss, depending on whether $T$ is (213)-decomposable or not.

\begin{figure}[h]
\labellist
\small\hair 2pt
\pinlabel $\scriptstyle{B}$ at -30 608
\pinlabel $\scriptstyle{C}$ at 98 608
\pinlabel $\scriptstyle{A}$ at 29 518
\pinlabel $\scriptstyle{\hbar}$ at 5 459
\pinlabel $\scriptstyle{A}$ at 194 518
\pinlabel $\scriptstyle{B+C}$ at 289 608
\pinlabel $\scriptstyle{\psi}$ at 341 567
\pinlabel $\scriptstyle{\psi(B+C)}$ at 500 608
\pinlabel $\scriptstyle{\psi(A)}$ at 419 518
\pinlabel $\scriptstyle{\tilde{B}}$ at 610 608
\pinlabel $\scriptstyle{\tilde{C}}$ at 740 608
\pinlabel $\scriptstyle{\psi(A)}$ at 688 518
\endlabellist
\begin{center}
\includegraphics[width=5.5in]{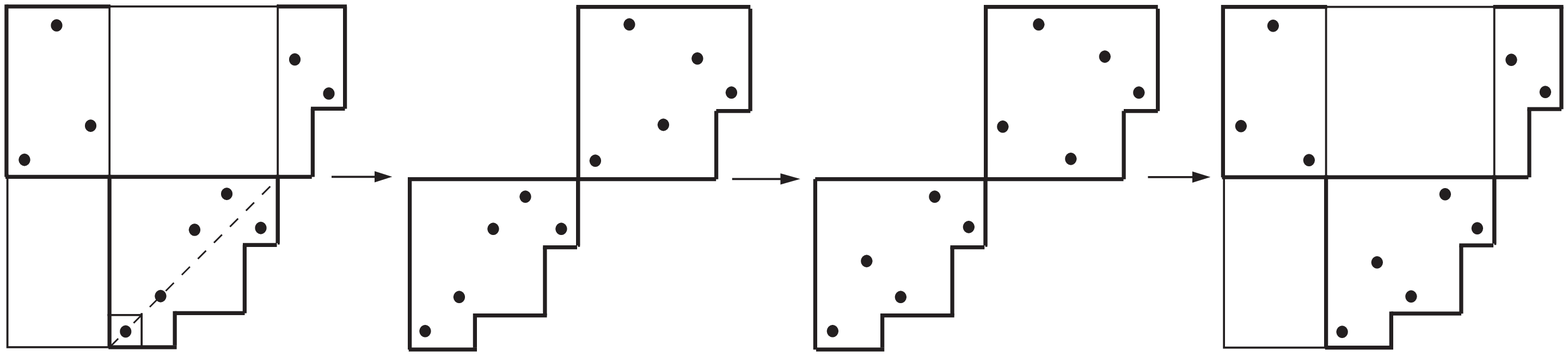}
\caption{$\psi(A\otimes (B+C))=\psi(A)\otimes \psi(B+C)=
\psi(A)\otimes (\tilde{B}+\tilde{C})$}
\label{Definition Phi Case 1}
\end{center} 
\end{figure}
\vspace*{-5mm}
\subsubsection{Case 1} $T$ has a {\it non-trivial}
$(213)$-decomposition; so consider its minimal $(213)$-decomposition
$T=T|_{\mathcal{A}_{\hbar}}\otimes(T|_{\mathcal{BL}_{\hbar}}+
T|_{\mathcal{BR}_{\hbar}})$, denoted for simplicity as $T=A\otimes
(B+C)$.  Since both $\mathcal{A}_{\hbar}$ and $\mathcal{B}_{\hbar}$
are of sizes $<n$, $\psi(A)$ and $\psi(B+C)$ are
well-defined by induction. We can further split
$\psi(B+C)=\tilde{B}+\tilde{C}$ where
$\tilde{B}=\psi(B+C)|_{\mathcal{BL}_{\hbar}}$ and
$\tilde{C}=\psi(B+C)|_{\mathcal{BR}_{\hbar}}$ occupy respectively the
same columns as the original $B$ and $C$.  We define $\psi(T)$ to be
the (213)-decomposable transversal of $Y$ given by
$\psi(T)=\psi(A)\otimes(\tilde{B}+\tilde{C})$. For example, the first
arrow in Figure~\ref{Definition Phi Case 1} signifies the
$(213)$-decomposition $Y_{213}(\hbar)$ which combines subboards
$B$ and $C$; the second arrow applies $\psi$ to each of $A$ and $B+C$,
and the third arrow splits $\psi(B+C)$ as $\tilde{B}$ and $\tilde{C}$
back into the original Young diagram $Y$.

\subsubsection{Case 2} 
\label{subsubsection Case 2}
\begin{figure}[h]
\labellist
\small\hair 2pt
\pinlabel $\scriptstyle{(Y\!/_{\!\hbar})^{\prime}}$ at 64 597
\pinlabel $\scriptstyle{(Y\!/_{\!\hbar})^{\prime\prime}}$ at 139 597
\pinlabel $\scriptstyle{D}$ at 108 563
\pinlabel $\scriptstyle{\hbar}$ at 10 459
\pinlabel $\scriptstyle{\psi|_{(Y\!/_{\!\hbar})^{\prime}}}$ at 205 572
\pinlabel $\scriptstyle{\hbar}$ at 245 459
\pinlabel $\scriptstyle{\mathcal{M}}$ at 278 530
\pinlabel $\scriptstyle{\eta}$ at 422 568
\pinlabel $\scriptstyle{D^{\prime\prime}}$ at 358 637
\pinlabel $\scriptstyle{D^{\prime}}$ at 265 637
\pinlabel $\scriptstyle{\dot{D^{\prime}}}$ at 234 670
\pinlabel $\scriptstyle{\ddot{D^{\prime}}}$ at 258 670
\pinlabel $\scriptstyle{\ddot{D^{\prime}}}$ at 470 670
\pinlabel $\scriptstyle{\hbar}$ at 533 459
\pinlabel $\scriptstyle{\psi(D)}$ at 345 562
\pinlabel $\scriptstyle{J_k}$ at 446 667
\pinlabel $\scriptstyle{D^{\prime\prime}}$ at 576 637
\endlabellist
\begin{center}
\includegraphics[width=4in]{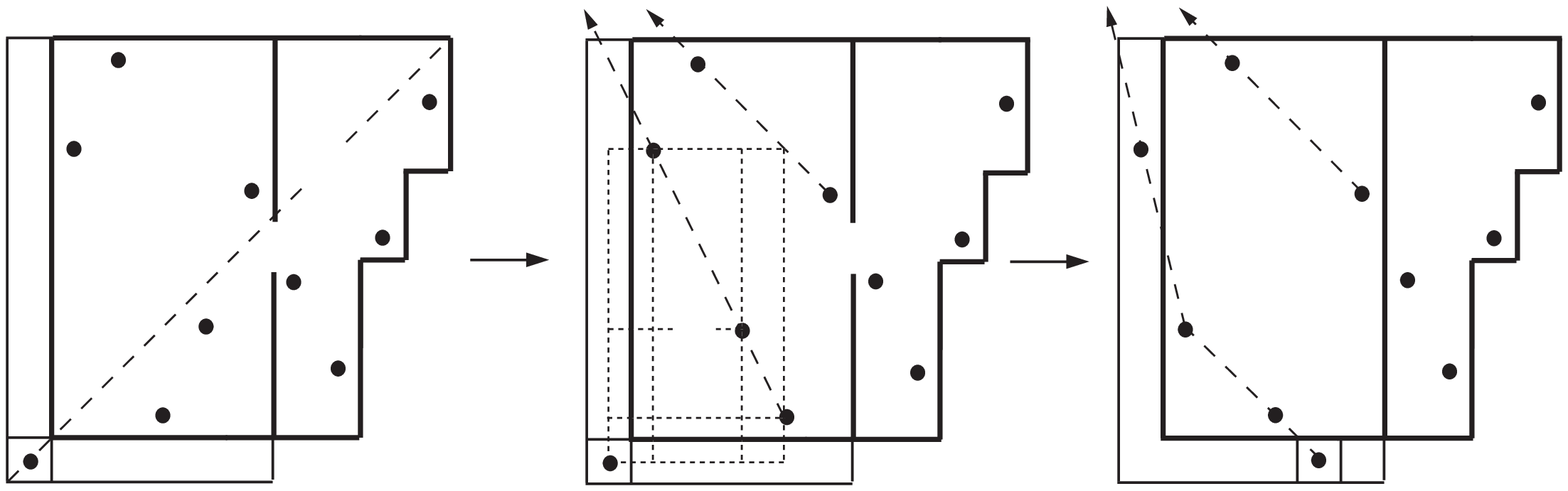}
\caption{Definition of $\psi$ in Case 2:
$T\rightarrow T_1\rightarrow \psi(T)=\eta(T_1)$}
\label{Definition Phi Case 2}
\end{center} 
\end{figure}

$T$ is (213)-indecomposable. As remarked
earlier, this can happen only if $\hbar$ is in the bottom left corner
cell of $Y$ and $Y$ has no $0$-critical points. Consider the
reduction $Y\!/_{\!\hbar}$, and let $D=T|_{Y\!/\!_{\hbar}}$ be its transversal
obtained from $T$ by removing $\hbar$. $Y\!/_{\!\hbar}$ breaks up into two
parts: the rectangle $(Y\!/_{\!\hbar})^{\prime}$ which lies over the
bottom row of $Y$, and the remaining subboard
$(Y\!/_{\!\hbar})^{\prime\prime}$ (cf. Fig.~\ref{Definition Phi Case 2}a.)
We define $\psi(T)$ in two steps: $T\rightarrow T_1\rightarrow
\psi(T)$.

Since $Y\!/\!_{\hbar}$ is of size $n-1$ and $D\in
S_{Y\!/\!_{\hbar}}(213)$, by induction $\psi(D)$ is defined as an
element of $S_{Y\!/\!_{\hbar}}(123)$. Let $T_1$ be the transversal of
$Y$ obtained from $\psi(D)$ by prepending $\hbar$ in its bottom left
corner; symbolically, $T_1=\phantom{}_{\hbar}[\psi(D)]$
(cf. Fig.~\ref{Definition Phi Case 2}b.)

\smallskip
The partition
$Y\!/\!_{\hbar}=(Y\!/\!_{\hbar})^{\prime}+(Y\!/\!_{\hbar})^{\prime\prime}$
induces a partition of the transversal
$\psi(D)=D^{\prime}+D^{\prime\prime}$. Since $D^{\prime}$ is a
$(123)$-avoiding partial transversal of the {\it rectangle}
$(Y\!/\!_{\hbar})^{\prime}$, then $D^{\prime}$ splits into its {\it
primary} $\dot{D^{\prime}}$ and {\it secondary} $\ddot{D^{\prime}}$
decreasing subsequences, as in Fig.~\ref{Definition Phi Case 2}b.  Let
$\mathcal{M}$ be the (square) submatrix of $Y$ induced by $\hbar$ and
$\ddot{D^{\prime}}$, and let
$T_1|_{\mathcal{M}}=\phantom{}_{\hbar}[\ddot{D^{\prime}}]$ be
$\mathcal{M}$'s transversal induced by $T_1$. For instance,
Figure~\ref{Definition Phi Case 2}b shows $T_1|_{\mathcal M}=(1432)$
and depicts $\mathcal{M}\cong M_4$ via dotted lines. If $k-1$ is the
length of $\ddot{D^{\prime}}$ ($k\geq 1$), then $\mathcal{M}\cong M_k$
and $T_1|_{\mathcal{M}} \approx (1,k,k-1,\ldots,3,2)$.  Let
$\eta(T_1)$ be the transversal of $Y$ obtained from $T_1$ by replacing
$T_1|_{\mathcal{M}}\mapsto J_k$ (cf. Fig.~\ref{Definition Phi Case
2}c.)

Set $\psi(T):=\eta(T_1)$ as the desired transversal in $S_Y$, and
define the map $\psi:S_Y(213)\rightarrow S_Y$ as the composition
$\psi=\eta \circ (\phantom{}_{\hbar}[\psi|_{Y/_{\hbar}}])$.

\subsection{Properties of $\psi:S_Y(213)\rightarrow S_Y(123)$}
\begin{prop}
For any $Y$ the map $\psi:S_Y(213)\rightarrow S_Y$ satisfies:

\begin{enumerate}
\item $\psi$ is a sequence of $(123)\!\rightarrow\!(213)$ moves;

\item $\psi$ maps to $S_Y(213)$;

\item $\psi$ is injective. 
\end{enumerate}
\label{properties of phi}
\end{prop}

A key idea in the proof of Proposition~\ref{properties of phi} is the
following lemma:

\begin{lem} 
Given $T\in S_Y$, all $(213)\!\rightarrow\!(123)$ and
$(123)\!\rightarrow\!(213)$ moves on $T$ respect any
$(213)$-decomposition of $T$.  Consequently, a sequence of such moves
preserves $(213)$-decomposability and 
$(213)$-indecomposability of transversals.
\label{preservation}
\end{lem}

\noindent
{\sc Proof:} Suppose $T\in S_Y$ has a $(213)$-decomposition
$T=A\otimes (B+C)$. The proof of Proposition~\ref{213-avoiding $T$
respect 213-decomposition} implies that a $(213)$-pattern inside $Y$
can involve elements only inside $A$ or only inside $B+C$.  Thus, a
$(213)\!\rightarrow\!  (123)$ move occurs entirely in $A$ or in $B+C$,
and hence it respects any $(213)$-decomposition of $T$. Using the
proof of Lemma \ref{123-avoiding converse}, the same reasoning shows
that $(123)\!\rightarrow\!(213)$ moves also respect any
$(213)$-decomposition of $T$.

It remains to prove that the two types of moves map indecomposable to
indecomposable transversals. To the contrary, if, say, a
$(123)\!\rightarrow\!(213)$ move maps an indecomposable $T$ to a
decomposable $T^{\prime}$, then the reverse
$(213)\!\rightarrow\!(123)$ move would map $T^{\prime}\rightarrow T$
and hence violate the preservation of $(213)$-decompositions proved
above. Thus, the two types of moves preserve $(213)$-indecomposability
too. \qed

\medskip
The properties of $\psi$ claimed in Proposition \ref{properties of
phi} are trivial for size $1$ Young diagrams. In the proof of
Proposition \ref{properties of phi}, we assume by induction that the
map $\psi$ satisfies the three required properties on all Young
diagrams of size smaller than $n$, and we fix a transversal $T\in
S_Y(213)$ for some $Y$ of size $n$.

\subsubsection{Proof of Proposition \ref{properties of phi}, Part
  (1):}\label{part 1} Suppose $T$ has a non-trivial
$(213)$-decomposition, so take the minimal such decomposition
$T=A\otimes (B+C)$. This is Case 1 of $\psi$'s definition $\psi$,
where $\psi$ consists of a move inside $A$ and a move inside
$B+C$. Hence $\psi$ respects this minimal $(213)$-decomposition. By
induction, $\psi(A)$ and $\psi(B+C)$ are each obtained (independently)
by $(123)\!\rightarrow\!(213)$ moves. Consequently, $\psi(T)$ is
obtained by the composition of all of these
$(123)\!\rightarrow\!(213)$ moves.
\vspace*{-5mm}
\begin{figure}[h]
\labellist
\small\hair 2pt
\pinlabel $\scriptstyle{(Y\!/_{\!\hbar})^{\prime}}$ at 41 579
\pinlabel $\scriptstyle{\mathcal{R}\,\,=}$ at -49 557
\pinlabel $\scriptstyle{\beta_1}$ at 12 603
\pinlabel $\scriptstyle{\beta_2}$ at 48 528
\pinlabel $\scriptstyle{\beta_3}$ at 66 495
\pinlabel $\scriptstyle{\alpha_1}$ at 31 638
\pinlabel $\scriptstyle{\alpha_2}$ at 71 595
\pinlabel $\scriptstyle{\dot{D^{\prime}}}$ at -4 668
\pinlabel $\scriptstyle{\ddot{D^{\prime}}}$ at -29 668
\pinlabel $\scriptstyle{\hbar}$ at -18 459
\pinlabel $\scriptstyle{\beta_1}$ at 144 613
\pinlabel $\scriptstyle{\beta_2}$ at 210 528
\pinlabel $\scriptstyle{\beta_3}$ at 228 495
\pinlabel $\scriptstyle{\alpha_1}$ at 192 638
\pinlabel $\scriptstyle{\alpha_2}$ at 233 595
\pinlabel $\scriptstyle{\dot{D^{\prime}}}$ at 158 668
\pinlabel $\scriptstyle{\hbar}$ at 164 459
\pinlabel $\scriptstyle{\beta_1}$ at 307 613
\pinlabel $\scriptstyle{\beta_2}$ at 336 528
\pinlabel $\scriptstyle{\beta_3}$ at 391 495
\pinlabel $\scriptstyle{\alpha_1}$ at 356 638
\pinlabel $\scriptstyle{\alpha_2}$ at 395 595
\pinlabel $\scriptstyle{\dot{D^{\prime}}}$ at 321 668
\pinlabel $\scriptstyle{\hbar}$ at 360 459
\pinlabel $\scriptstyle{\beta_1}$ at 470 613
\pinlabel $\scriptstyle{\beta_2}$ at 500 528
\pinlabel $\scriptstyle{\beta_3}$ at 536 495
\pinlabel $\scriptstyle{\alpha_1}$ at 519 638
\pinlabel $\scriptstyle{\alpha_2}$ at 558 595
\pinlabel $\scriptstyle{\dot{D^{\prime}}}$ at 485 668
\pinlabel $\scriptstyle{\hbar}$ at 540 459
\pinlabel $\scriptstyle{J_k}$ at 449 667
\endlabellist
\begin{center}
\includegraphics[width=4.5in]{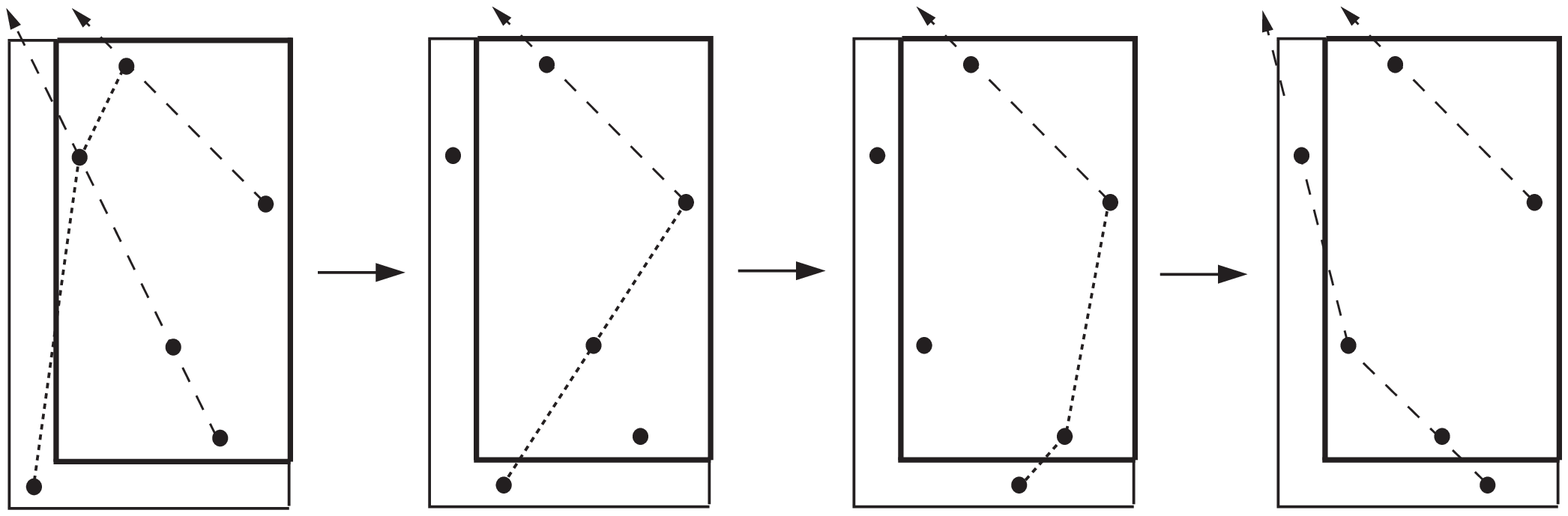}
\caption{$(\hbar,\beta_1,\beta_2,\beta_3)\!\rightarrow\! (\beta_1,\hbar,
\beta_2,\beta_3)\!\rightarrow\! (\beta_1,\beta_2,\hbar,\beta_3)
\!\rightarrow\! (\beta_1,\beta_2,\beta_3,\hbar)$ in $\mathcal{R}$}
\label{123-213 Moves}
\end{center} 
\end{figure}

Suppose now that $T$ is $(213)$-indecomposable. This is Case 2 of
$\psi$'s definition. By induction on $D$'s size, $\psi(D)$ and
therefore $T_1$ are obtained by $(123)\!\rightarrow\!(213)$ moves
inside $D$. It remains to show that $\psi(T_1|_{\mathcal M})=J_k$ can
be obtained from $T_1|_{\mathcal M}$ via such moves too. Recall that
each $\beta\in\ddot{D}^{\prime}$ is $(12)$-dominated in the rectangle
$(Y\!/\!_{\hbar})^{\prime}$ by some $\alpha\in\dot{D}^{\prime}$. Thus,
as long as $\hbar$ is before $\beta$ ($\hbar$ still in the bottom row
of $Y$), then $(\hbar \beta \alpha)$ is a $(123)$-pattern inside
$(Y\!/\!_{\hbar})^{\prime}$, and hence in $Y$. So the move
$(\hbar\beta\alpha)\!\mapsto\! (\beta\hbar\alpha)$ is a
$(123)\!\rightarrow\!(213)$ move which leaves $\hbar$ still in the
bottom row of $Y$. Let
$\ddot{D}^{\prime}=(\beta_1,\beta_2,...,\beta_{k-1})\!\!\searrow$.
Then $(\hbar,\beta_1,\beta_2,...,\beta_{k-1})\approx
(1,k,k-1,...,3,2)$. Using the above reasoning, we can switch $\hbar$
consecutively with each of the $\beta_i$'s via some
$(123)\!\rightarrow\!(213)$ move. For example, Figure~\ref{123-213
Moves} depicts this situation in the rectangular part
$\mathcal{R}=\phantom{}_{\hbar}[(Y\!/\!_{\hbar})^{\prime}]$ of $Y$
projecting onto $Y$'s bottom row, where each
$(123)\!\rightarrow\!(213)$ move is marked with a dotted line. Hence,
the sequence
\[(\hbar,\beta_1,\beta_2,...,\beta_{k-1})\rightarrow (\beta_1,\hbar,
\beta_2,...,\beta_{k-1})\rightarrow
(\beta_1,\beta_2,\hbar,...,\beta_{k-1})
\rightarrow \cdots \rightarrow
(\beta_1,\beta_2,...,\beta_{k-1},\hbar)\] is a composition of
$(123)\rightarrow (213)$ moves, and so is $\psi$. \qed

\subsubsection{Proof of Proposition \ref{properties of phi}, Part
(2):} \label{part 2} In Case 1 of $\psi$'s definition, by induction
  $\psi(A)$ and $\psi(B+C)$ are both $(123)$-avoiding. Lemma
  \ref{123-avoiding converse} implies that no new $(123)$-pattern can
  be introduced in the $(213)$-decomposition $\psi(T)=
  \psi(A)\otimes(\tilde{B}+\tilde{C})$.  We conclude that $\psi(T)\in
  S_Y(123)$.

In Case 2 of $\psi$'s definition, by induction $\psi(D)$ avoids
$(123)$. The only $(123)$-patterns in $Y$ before applying $\eta$ can
occur because of $\hbar$ being prepended to $D$'s bottom left corner,
and hence any such pattern can appear only in the rectangle
$\mathcal{R}=\phantom{}_{\hbar}[(Y\!/\!_{\hbar})^{\prime}]$
(cf. Fig.~\ref{123-213 Moves}.) A $(123)$-pattern in $\mathcal{R}$ is
of the form $(\hbar\beta\alpha)$ with $\alpha\in \dot{D}^{\prime}$,
$\beta\in \ddot{D}^{\prime}$ and $\alpha$ $(12)$-dominates $\beta$ in
$\mathcal{R}$. The map $\eta$ eliminates all these $(123)$-patterns by
shifting $\ddot{D}^{\prime}$ horizontally to the left and $\hbar$ to
the right until $\hbar$ is after all of $\ddot{D}^{\prime}$. Thus,
after $\eta$ is applied, $\hbar$ cannot participate in any more
$(123)$-patterns in $\psi(T)$.

\begin{figure}[h]
\labellist
\small\hair 2pt
\pinlabel $\scriptstyle{(Y\!/_{\!\hbar})^{\prime}}$ at 211 635
\pinlabel $\scriptstyle{(Y\!/_{\!\hbar})^{\prime\prime}}$ at 290 635
\pinlabel $\scriptstyle{\alpha}$ at 171 519
\pinlabel $\scriptstyle{\alpha_1}$ at 207 519
\pinlabel $\scriptstyle{\beta}$ at 227 597
\pinlabel $\scriptstyle{\eta}$ at 193 537
\pinlabel $\scriptstyle{\gamma}$ at 304 615
\pinlabel $\scriptstyle{{D^{\prime\prime}}}$ at 280 566
\pinlabel $\scriptstyle{\hbar}$ at 225 459
\pinlabel $\scriptstyle{(Y\!/_{\!\hbar})^{\prime}}$ at 429 635
\pinlabel $\scriptstyle{(Y\!/_{\!\hbar})^{\prime\prime}}$ at 508 635
\pinlabel $\scriptstyle{\alpha}$ at 389 519
\pinlabel $\scriptstyle{\alpha_1}$ at 425 519
\pinlabel $\scriptstyle{\beta}$ at 407 597
\pinlabel $\scriptstyle{\eta}$ at 411 537
\pinlabel $\scriptstyle{\gamma}$ at 521 615
\pinlabel $\scriptstyle{{D^{\prime\prime}}}$ at 498 566
\pinlabel $\scriptstyle{\hbar}$ at 443 459
\pinlabel $\scriptstyle{(Y\!/_{\!\hbar})^{\prime}}$ at -5 635
\pinlabel $\scriptstyle{(Y\!/_{\!\hbar})^{\prime\prime}}$ at 74 635
\pinlabel $\scriptstyle{\alpha}$ at -45 519
\pinlabel $\scriptstyle{\alpha_1}$ at -9 519
\pinlabel $\scriptstyle{\beta_1}$ at 459 566
\pinlabel $\scriptstyle{\beta}$ at 45 603
\pinlabel $\scriptstyle{\eta}$ at -23 537
\pinlabel $\scriptstyle{\gamma}$ at 88 615
\pinlabel $\scriptstyle{{D^{\prime\prime}}}$ at 64 566
\pinlabel $\scriptstyle{\hbar}$ at 9 459
\endlabellist
\begin{center}
\includegraphics[width=4.5in]{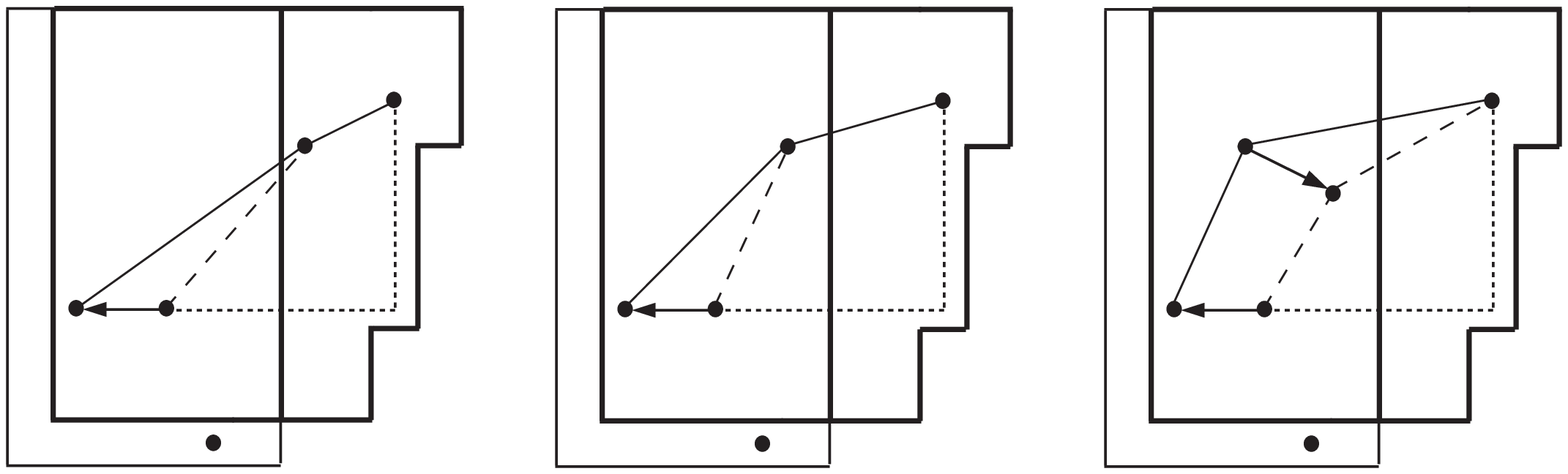}
\caption{$\psi:S_Y(213)\rightarrow\in S_Y(123)$}
\label{Proof of Phi Case 2}
\end{center} 
\end{figure}

It remains to show that $\eta$ has not created any new
$(123)$-patterns $(\alpha\beta\gamma)$ which do not involve
$\hbar$. Since $\eta$ preserves $\dot{D}^{\prime}$ and
$D^{\prime\prime}$, such $(123)$-pattern must involve an element of
$\eta(\ddot{D}^{\prime})$; in fact, $\alpha\in\eta(\ddot{D}^{\prime})$
since the elements of $\eta(\ddot{D}^{\prime})$ do not
$(12)$-dominate anything and hence they can play only the role of
``$1$'' in a $(123)$-pattern. Since
$\eta(\ddot{D}^{\prime})\!\!\searrow$ and
$\dot{D}^{\prime}\!\!\searrow$, at most one element of each can
participate in this $(123)$-pattern. Finally, at most one element of
$D^{\prime\prime}$ can participate too; indeed, suppose two elements
$\beta$ and $\gamma$ of $D^{\prime\prime}$ participate in the
$(123)$-pattern $(\alpha\beta\gamma)$ with $\alpha\in
\eta(\ddot{D}^{\prime})$ (cf. Fig.~\ref{Proof of Phi Case 2}a.) If
$\alpha_1=\eta^{-1}(\alpha)$, then $(\alpha_1\beta\gamma)$ would also
be a $(123)$-pattern in $\psi(D)$ since $\alpha_1$ is a horizontal
shift of $\alpha$ to the right but still inside
$(Y|_{\hbar})^{\prime}$ and before $\beta,\gamma\in
(Y|_{\hbar})^{\prime\prime}$. This is a contradiction with the
inductive assumption that $\psi(D)$ is $(123)$-avoiding. We conclude
that at most one element of $D^{\prime\prime}$ can participate in the
$(123)$-pattern, and therefore $(\alpha\beta\gamma)$ is formed by
$\alpha\in \eta(\ddot{D}^{\prime})$, $\beta\in \dot{D}^{\prime}$, and
$\gamma\in D^{\prime\prime}$.

If $\alpha_1=\eta^{-1}(\alpha)\in \ddot{D}^{\prime}$ is {\it before}
$\beta$, then $\alpha_1$ is $(12)$-dominated by $\beta$
(cf. Fig.~\ref{Proof of Phi Case 2}b) and so
$(\alpha_1\beta\gamma)\approx(123)$ in $\psi(D)$, a contradiction. If
$\alpha_1$ is {\it after} $\beta$, then $\alpha_1$ must be
$(12)$-dominated by some $\beta_1\in \dot{D}^{\prime}$, $\beta_1\not =
\beta$ (cf. Fig.~\ref{Proof of Phi Case 2}c.)  Since $\alpha_1$ is
after $\beta$, and $\beta_1$ is after $\alpha$, then $\beta_1$ is
after $\beta$. Because $\beta,\beta_1\in
\dot{D}^{\prime}\!\!\searrow$, it follows that
$(\beta\beta_1)\!\!\searrow$. Further, $\gamma\in D^{\prime\prime}$,
hence $\gamma$ comes after both $\beta$ and $\beta_1$. Finally,
$\beta_1<\beta<\gamma$ implies that
$(\alpha_1\beta_1\gamma)\!\!\nearrow$ is a $(123)$-subpattern of
$\psi(D)$ landing in $Y$: indeed, the intersection of $\gamma$'s
column and $\alpha_1$'s row is the same as the intersection of
$\gamma$'s column and $\alpha$'s row, and the latter is inside $Y$ by
the assumption that $(\alpha\beta\gamma)$ is a $(123)$-pattern in
$Y$. The existence of such $(\alpha_1\beta_1\gamma)$ contradicts
$(123)$-avoidance of $\psi(D)$.

Therefore, $\eta$ gets rid of all $(123)$-patterns involving $\hbar$ and
does not introduce any new $(123)$-patterns, so that $\eta(T_1)$ is
$(123)$-avoiding and $\psi(T)\in S_Y(123)$. \qed

\subsubsection{Proof of Proposition \ref{properties of phi}, Part
(3):}\label{part 3} Let $T_2\in \psi(S_Y(213))\subset S_Y(123)$ for
  some Young diagram $Y$ of size $n$, and let $T\in S_Y(213)$ be any
  preimage of $T_2$, i.e. $\psi(T)=T_2$. We will show that $T$ can be
  recovered uniquely from $T_2$.

\smallskip
{\it Case 1.} Suppose that $T_2$ is $(213)$-{\it decomposable}, and
let $T_2=\tilde{A}\otimes (\tilde{B}+\tilde{C})$ be $T_2$'s minimal
$(213)$-decomposition. By Property 1 of $\psi$,
$T\stackrel{\psi}{\mapsto}T_2$ is a sequence of
$(123)\!\rightarrow\!(213)$ moves; inverting each of these moves, we
obtain a sequence of $(213)\!\rightarrow\!(123)$ moves that takes
$T_2\mapsto T$. By Lemma~\ref{preservation}, $T$ respects the
$(213)$-decomposition of $T_2$; moreover, the induced
$(213)$-decomposition $T=A\otimes(B+C)$ is also minimal,
i.e. $\hbar\in T$ is in the bottom left corner of $A$
(cf. Fig.~\ref{Definition Phi Case 1}.) Thus,
$\psi(A\otimes(B+C))=\tilde{A}\otimes(\tilde{B}+\tilde{C})$ where
$\psi(A)=\tilde{A}$ and $\psi(B+C)=\tilde{B}+\tilde{C}$ by the
$\psi$'s definition in Case 1. By induction, $\psi$ is injective on
smaller size Young diagrams, so that $A$ and $B+C$ can be recovered
from $\tilde{A}$ and $\tilde{B}+\tilde{C}$.  Finally, since the
decomposition of $T$ is determined by $T_2$, $B$ and $C$ themselves
can be recovered uniquely from $B+C$.  We conclude that $T$ can be
recovered uniquely from $T_2$.

\smallskip
{\it Case 2.} Suppose that $T_2$ is $(213)$-indecomposable. As in Case
1, $\psi$ must have preserved this property, i.e. $T$ is also
$(213)$-indecomposable. Since $T\in S_Y(213)$, this implies that
$\hbar$ in $T$ is in the bottom left corner of $Y$. By $\psi$'s
definition in Case 2, there is an intermediate
$T_1=\phantom{}_{\hbar}[\psi(D)]$ such that $\psi(T)=\eta(T_1)=T_2$
(cf. Fig.~\ref{Definition Phi Case 2}.)  We will first show that $T_1$
is recoverable from $T_2$.

In $T_2$ we can uniquely determine $D^{\prime\prime}$ as the
subtransversal in the part of $Y$ that does not project on the bottom
row of $Y$. The remainder $\tilde{D}$ of $T_2$ projects onto the
bottom row of $Y$ and lies in the rectangle $\mathcal{R}$. Since $T_2$
avoids $(123)$, $\tilde{D}$ splits into its primary and secondary
subsequences, $\tilde{D}_1$ and $\tilde{D}_2$, respectively.  Note
that $\hbar$ in $\tilde{D}$ is $(12)$-dominated: being in the bottom
row of $\mathcal{R}$, the only way for $\hbar$ not to be
$(12)$-dominated is to be in the rightmost (bottom) cell of
$\mathcal{R}$; but then $T_2$ would be decomposable, contradicting our
assumption in this case. Thus, $\hbar\in \tilde{D}_2$.

Now consider $T_1$. By $\eta$'s definition in Case 2,
$T_1=\phantom{}_{\hbar}[\ddot{D}^{\prime}+\dot{D}^{\prime}+D^{\prime\prime}]$,
where $\eta$ fixes $\dot{D}^{\prime}$ and $D^{\prime\prime}$, slides
$\ddot{D}^{\prime}$ to the left, and slides $\hbar$ to the right until
$\hbar$ is after $\ddot{D}^{\prime}$. Sliding $\ddot{D}^{\prime}$ to
the left leaves all of its elements $(12)$-dominated by some elements
in $\dot{D}^{\prime}$, and as we argued above, it makes $\eta(\hbar)$
also $(12)$-dominated in $T_2$. In other words,
$\eta(\dot{D}^{\prime})=\dot{D}^{\prime}=\tilde{D}_1$ and
$\eta(_{\hbar}[\ddot{D}^{\prime}])=\tilde{D}_2$. Thus, to recover
$T_1$ from $T$, we keep $D^{\prime\prime}$ and $\tilde{D}_1$, and
switch horizontally the places of $\tilde{D}_2\backslash{\hbar}$ and
$\hbar$. Note that at this point $\hbar\in T_1$ must be in the bottom
left corner of $Y$ by $\psi$'s definition in Case 2.

To recover $T$ from $T_1$, note that by induction $\psi$ is injective
on $Y\!/\!_{\hbar}$, so that $\psi(D)$ in $T_1$ could have come only
from one transversal $D$; appending $\hbar$ at the bottom left corner of
$D$ gives the unique preimage $T=\phantom{}_{\hbar}[D]\in
S_Y(213)$. \qed

\subsection{Conclusions}
Subsubsections~\ref{part 1}-3 complete inductively the proof of
Proposition~\ref{properties of phi}. The latter implies that
$\psi:S_Y(213)\hookrightarrow S_Y(123)$ is injective for all 
Young diagrams $Y$. Thus, $|S_Y(213)|\leq |S_Y(123)|$ and
Theorem~\ref{SWOS3} is proven. \qed

\medskip
We leave the following questions to the reader for further study. For
which pairs of permutations $\sigma$ and $\tau$ in $S_k$ can a map
$\psi_Y:S_Y(\sigma)\rightarrow S_Y(\tau)$ be well-defined via
$\sigma\rightarrow \tau$ moves? What properties does $\psi$ possess in
such cases?

\section{Strict Inequalities $|S_Y(213)|<|S_Y(123)|$}
\label{strict |SY(213)|<|SY(123)|}

Below we refer to the notation from the definition of the map
$\psi:S_Y(213)\hookrightarrow S_Y(123)$ in Section~\ref{SY(213) <
SY(123)}; in particular, $\psi(T)=\eta(T_1)=T_2$ for any
$(213)$-indecomposable $T\in S_Y(213)$.

\begin{lem} If $Y$ has an $i$-critical point with $i\geq 2$ and no
$0$- and $1$-critical points, some $(213)$-indecomposable 
 $T_2\in S_Y(123)$ is not invertible under $\eta$ and hence under
 $\psi$.
\label{lemma 123>213}
\end{lem}

\begin{figure}[h]
\labellist
\small\hair 2pt
\pinlabel $\scriptstyle{\alpha}$ at 28 610
\pinlabel $\scriptstyle{\beta}$ at 189 625
\pinlabel $\scriptstyle{\bar{T}\,\,\text{on}\,\,\bar{Y}}$ at 123 565
\pinlabel $\scriptstyle{\hbar}$ at 45 459
\pinlabel $\scriptstyle{\alpha}$ at 261 610
\pinlabel $\scriptstyle{\beta}$ at 414 638
\pinlabel $\scriptstyle{\mathcal{R}}$ at 330 479
\pinlabel $\scriptstyle{\hbar}$ at 279 459
\endlabellist
\begin{center}
\includegraphics[width=2.5in]{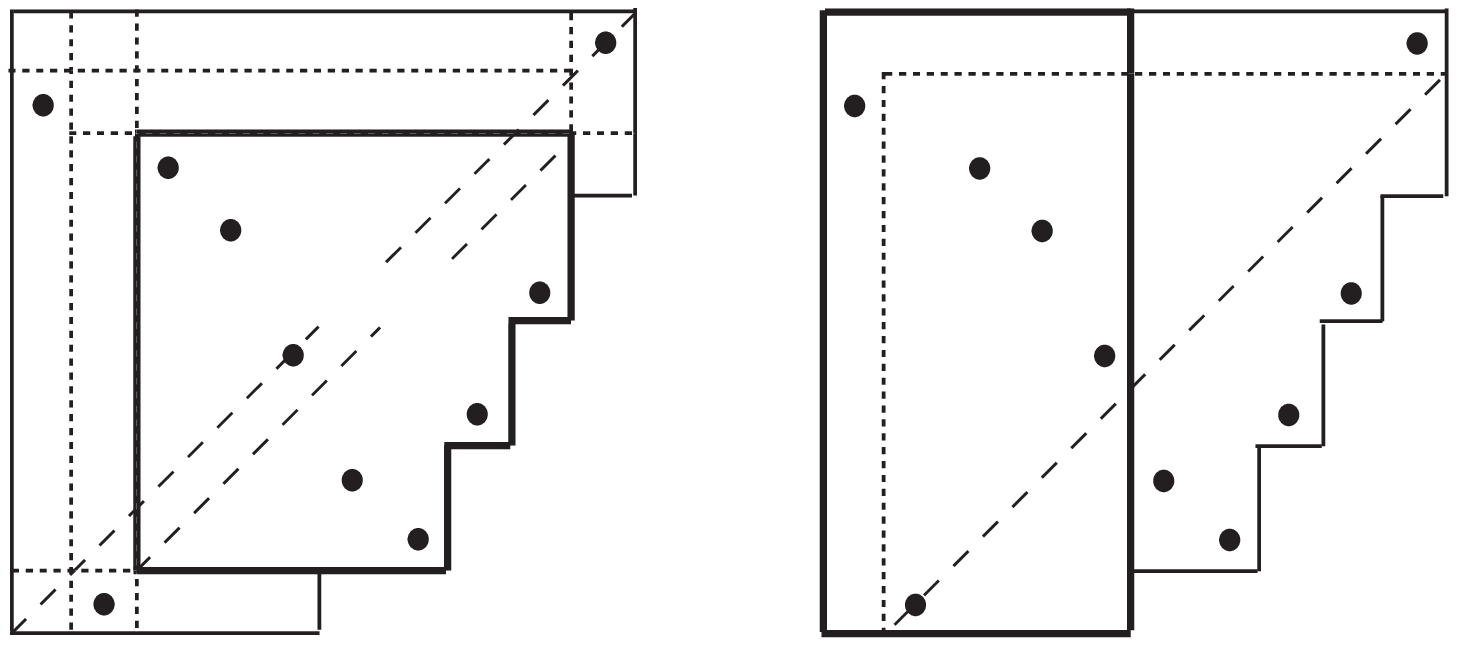}
\caption{$T_2\in S_Y(123)\backslash \psi(S_Y(213))$}
\label{Sufficient Condition 2}
\end{center} 
\end{figure}

\noindent
{\sc Proof:} Since $i\geq 2$, the size of $Y$ is $n\geq 4$. Place
$\alpha$ in position $(2,1)$, $\hbar$ in $(n,2)$ and $\beta$ in
$(1,n)$, and set $\overline{Y}=Y/\{\alpha,\hbar,\beta\}$. The
hypotheses on $Y$ imply that $\overline{Y}$ is non-empty and has no
$0$-critical points.  Since $(12)\sim_s(21)$ and since there is
obviously exactly 1 transversal of $\overline{Y}$ that avoids $(21)$
(namely, the diagonal transversal), there is also exactly 1
transversal $\overline{T}$ of $\overline{Y}$ that avoids $(12)$
(cf. Fig.~\ref{Sufficient Condition 2}a.) Thus, $T_2=\{\alpha,\hbar,
\overline{T},\beta\}$ is a transversal of $Y$. We claim that $T_2\in
S_Y(123)\backslash \psi(S_Y(213))$.

To show that $T_2$ avoids $(123)$, note that the position of $\alpha$
in $Y$ precludes it from participating in any such pattern. Moreover,
$\hbar$ and $\beta$ cannot simultaneously participate in a
$(123)$-pattern since $\beta$'s column and $\hbar$'s row do not
intersect inside the non-square $Y$.  Yet, at most 1 element from
$\overline{T}$ can participate in a $(123)$-pattern due to the
$(12)$-avoidance of $\overline{T}$. This does not leave enough
elements of $T_2$ to participate in a $(123)$-pattern in $Y$.

Next, in any $(213)$-decomposition of $T_2=A\otimes (B+C)$, $A$
contains $\hbar$ and hence the $1$-diagonal $d_1(Y)$, which starts
from $\hbar$ (cf. Fig.~\ref{Sufficient Condition 2}b.) The hypotheses
on $Y$ and the position of $\hbar$ imply that $d_1(Y)$ does not
intersect the border of $Y$ until goes through the rightmost column of
$Y$ and stops underneath $\beta$'s cell. This forces the
subtransversal $A$ to involve the second row of $Y$ and hence to
contain $\alpha$, as well as the rightmost column of $Y$ and hence to
contain $\beta$, i.e. $A=T_2$ and the $(213)$-decomposition of $T_2$ is
trivial. Therefore, $T_2$ is $(213)$-indecomposable.

From Lemma~\ref{preservation}, if a preimage $T\in S_Y(213)$ of $T_2$
existed under $\psi$, then $T$ would also be $(213)$-indecomposable
and by Case 2 of $\psi$'s definition: $T\stackrel{\psi}{\rightarrow}
T_1\stackrel{\eta}{\rightarrow}T_2$.  In particular, $T$ and $T_1$
would have $\hbar$ is their bottom left corners. But $\alpha\in T_2$
is $(12)$-dominated only by $\beta$, and $\beta$ does not project on
the bottom row of $Y$, hence $\alpha$ is {\bf not} $(12)$-dominated in
the rectangle $\mathcal{R}$, hence $\alpha\in \dot{D}^{\prime}$
(cf. Fig.~\ref{123-213 Moves}.)  Since $\eta$ fixes
$\dot{D}^{\prime}$, inverting $\eta$ would leave $\alpha$ fixed in the
first column of $T_1$. This precludes $\hbar$ from occupying the
bottom left corner in $T_1$, a contradiction.  We conclude that $T_2$
is not invertible under $\eta$ and $\psi$, and hence
$T_2\not\in\psi(S_Y(213))$. \qed

\medskip
\begin{prop} $|S_Y(213)|<|S_Y(123)|$ if and only if $Y$ has 
an $i$-critical point with $i\geq 2$.
\label{iff 213<123}
\end{prop}

\noindent
{\sc Proof:} As in the proof of Proposition~\ref{strict inequality
  312>321}, for any permutation $\sigma$ we can split $Y$ and its
  transversals with respect to any $0$- and $1$-critical points:
\[|S_Y(\sigma)|=|S_{U_1}(\sigma)|\cdot |S_{U_2}(\sigma)|\cdots 
|S_{U_k}(\sigma)|,\] where each $U_j$ is either square or contains
only $i$-critical points with $i\geq 2$. If the original $Y$ contains
only $0$- or $1$-critical points, then all $U_j$'s are square with
$|S_{U_j}(213)|=|S_{U_j}(123)|$, so that $|S_Y(213)|=|S_Y(123)|$.

If $Y$ does contain some high $i$-critical points with $i\geq 2$, in
addition to the square $U_j$'s, there will be at least one other $U_m$
with such a high critical point. Lemma~\ref{lemma 123>213}
implies strict inequalities for all non-square $U_j$'s in our
decomposition. In particular,
$|S_{U_m}(213)|<|S_{U_m}(123)|$ and therefore
$|S_Y(213)|<|S_Y(123)|$. \qed

\medskip
This completes the proof of Theorem~\ref{summary 312>321>213}. \qed

\section{Strict Wilf-ordering for $(213|\tau)$, $(123|\tau)$ and $(312|\tau)$}
\label{strict Wilf-ordering}
Subsection~\ref{strategy} gives a strategy for proving that for any
permutation $\tau$: \[|S_n(213|\tau)|\lneqq |S_n(123|\tau)|\lneqq
|S_n(312|\tau)|\,\,\text{for}\,\,n\gg 1.\] Since each Young diagram
$Y_m$ has an $(m-2)$-critical point, Theorem~\ref{summary 312>321>213}
implies $|S_{Y_m}(213)|\lneqq |S_{Y_m}(123)|\lneqq |S_{Y_m}(312)|$ for
$m\geq 5$. This fulfills the first step (SF1) of the strategy.  The
other step (SF2) is provided by the following construction.

\begin{lem}
\label{lemma saturation of Y_n}
Given a permutation $\tau\in S_k$, for any $n\geq 2k+2$ there is a
partial transversal $T_n$ of $M_n$ which saturates $Y_{n-2k}$ with
respect to $\tau$.
\end{lem}

\noindent{\sc Proof:} Take two copies $\tau_1$ and $\tau_2$ of $\tau$
and arrange them in a southwest/northeast diagonal fashion within a
square matrix $M_{2k}$ (cf. Fig.~\ref{saturation of Y_n}a.) Insert a
row and column through the middle of $M_{2k}$ so that the resulting
$\mathcal{M}\cong M_{2k+1}$ has an empty separating row and column
between $\tau_1$ and $\tau_2$. Place $\mathcal{M}$ in the bottom right
corner of $M_n$ for $n\geq 2k+2$ (cf. Fig.~\ref{saturation of Y_n}b.)
We claim that the partial transversal $T_n$ of $M_n$ produced by the
two copies $\tau_1$ and $\tau_2$ in $\mathcal{M}$ saturates $Y_{n-2k}$
with respect to $\tau$.
\begin{figure}[h]
\labellist
\small\hair 2pt
\pinlabel ${\tau_1}$ at 58 520
\pinlabel ${\tau_2}$ at 129 592
\pinlabel ${W^{\prime}}$ at 253 685
\pinlabel $\scriptstyle{w_1}$ at 279 565
\pinlabel $\scriptstyle{w_2}$ at 368 655
\pinlabel ${c}$ at 369 566
\pinlabel ${\tau_1}$ at 327 519
\pinlabel ${\tau_2}$ at 416 609
\pinlabel ${\mathcal{M}}$ at 410 531
\pinlabel ${M_n}$ at 236 500
\pinlabel ${W^{\prime}}$ at 543 684
\pinlabel $\scriptstyle{w_1}$ at 569 565
\pinlabel $\scriptstyle{w_2}$ at 659 655
\pinlabel ${W=Y_{n-2k}}$ at 767 680
\pinlabel $\scriptstyle{w_1}$ at 785 638
\pinlabel $\scriptstyle{w_2}$ at 802 655
\pinlabel ${?}$ at 587 620
\pinlabel ${?}$ at 604 620
\pinlabel ${?}$ at 587 637
\pinlabel ${?}$ at 604 637
\pinlabel ${?}$ at 315 620
\pinlabel ${?}$ at 297 620
\pinlabel ${?}$ at 315 637
\pinlabel ${?}$ at 297 637
\endlabellist
\begin{center}
\includegraphics[width=6in]{Saturation-of-Yn}
\caption{$T_n=\tau_1\oplus \tau_2$ saturates $W=Y_{n-2k}$ with respect
to $\tau$}
\label{saturation of Y_n}
\end{center} 
\end{figure}

To see this, denote by $w_1$ and $w_2$ the cells of $Y$ diagonally to
the left and above the $k\times k$ matrices of $\tau_1$ and
$\tau_2$. Then $w_1$ and $w_2$ are white cells with respect to $\tau$
and the partial transversal $T_n$. The existence of $w_1$ and $w_2$
gives $n\geq 2k+2$. If $c$ is the central cell of $\mathcal{M}$, then
the initial white subboard $W^{\prime}$ is of the union $Y_{w_1}\cup
Y_{w_2}$ of two rectangles, plus possibly some more white cells within
the rectangle $\mathcal{M}_{\bar{c}}$ (these cells are depicted by
``?''  in Fig.~\ref{saturation of Y_n}b-c). However, the reduction of
$W^{\prime}$ along $\tau_1\cup \tau_2$ deletes {\it all} cells in
$\mathcal{M}_{\bar{c}}$, and the consequent removal of the (blue)
central row and column of $\mathcal{M}$ leaves the white diagram
$W=Y_{n-2k}$ (cf. Fig.~\ref{saturation of Y_n}d).  By definition,
$T_n$ saturates $Y_{n-2k}$ with respect to $\tau$ in $M_n$. \qed

\medskip
When $n\geq 2k+5$, then the saturated $Y_{n-2k}$ satisfies
(SF1). Combining with (SF2),
\[|S_{Y_{n-2k}}(213)|\cdot |\bar{S}_{M_n\backslash
Y_{n-2k}}(\tau)|\lneqq |S_{Y_{n-2k}}(123)|\cdot
  |\bar{S}_{M_n\backslash Y_{n-2k}}(\tau)|\lneqq
  |S_{Y_{n-2k}}(312)|\cdot |\bar{S}_{M_n\backslash
    Y_{n-2k}}(\tau)|\]
\[\Rightarrow\,\,|S_n(213|\tau)|\lneqq |S_n(123|\tau)|\lneqq
|S_n(312|\tau)|\,\,\text{for}\,\,n\geq 2k+5.\]  This completes the proof of
Corollary~\ref{Wilf-ordering-corollary}. \qed
\begin{figure}[h]
\labellist
\small\hair 2pt
\pinlabel ${W^{\prime}}$ at 47 683
\pinlabel ${W^{\prime}}$ at 363 674
\pinlabel ${W=Y_4}$ at 198 665
\pinlabel ${W_1}$ at 524 672
\endlabellist
\begin{center}
\includegraphics[width=4in]{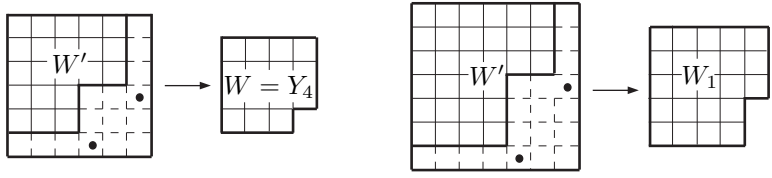}
\caption{$T^{\prime}$ saturates $W=Y_4$ with respect
to $\tau=(1)$}
\label{Example 3}
\end{center} 
\end{figure}

\begin{ex}
{\rm As an illustration of the above inequalities, let $\tau=(1)$ and
consider $(213|1)\preceq (123|1) \preceq (312|1)$. When $n=6,7$:
\begin{eqnarray}\label{S_6}
|S_6(3241)|=512<|S_6(2341)|=|S_6(4231)|=513,\\\label{S_7}
|S_7(3241)|=2740<|S_7(2341)|=2761<|S_7(4231)|=2762.
\end{eqnarray}
Let $T^{\prime}$ be a partial transversal of $M_6$ that saturates a
$W\subset M_6$ with respect to $\tau=(1)$. Then $|W|\leq 5$ with
$|W|=5$ if and only if $T^{\prime}$ consists of a single element in
the bottom right corner of $M_6$; in such a case $W=M_5$. Thus,
$|S_W(213)|=|S_W(123)|=|S_W(312)|$ for all $W\subset M_6$ except
$W=Y_4$, where the $2$-critical point of $Y_4$ implies the inequality
$|S_{Y_4}(213)|=12<|S_{Y_4}(123)|=|S_{Y_4}(312)|=13$. On the other
hand, it is easy to verify that the only $T^{\prime}$ that saturates
$Y_4$ in $M_6$ consists of two elements placed in positions $(4,6)$
and $(6,4)$ (cf. Fig.~\ref{Example 3}a-b.) Thus, the Splitting
Formulas for $S_6(3241)$, $S_6(2341)$ and $S_6(4231)$ have all but one
equal summands:
\[12\cdot 1=|S_{Y_4}(213)|\cdot |S|\lneqq |S_{Y_4}(123)|\cdot
|S|= S_{Y_4}(312)|\cdot |S|=13\cdot 1,\] where
$S=\bar{S}_{M_6\backslash Y_4}(1)$. This explains the difference of
$1$ between the quantities in (\ref{S_6}).

\smallskip
The analogous partial transversal $T^{\prime\prime}$ in $M_7$ (whose
two elements are placed in $(5,7)$ and $(7,5)$) saturates $Y_5$ with
respect to $\tau=(1)$. The $3$-critical point of $Y_5$ implies the
following inequalities, where $S=\bar{S}_{M_7\backslash Y_5}(1)$:
\[37\cdot 1=|S_{Y_5}(213)|\cdot |S|
\lneqq 41\cdot 1=|S_{Y_5}(123)|\cdot |S| \lneqq 42\cdot 1=
|S_{Y_5}(312)|\cdot |S|.\] This explains the difference of $1$ between
$|S_7(2341)|$ and $|S_7(4231)|$ in (\ref{S_7}). Further, $Y_5$,
$W_1\cong Y(5,5,5,4,4)$, its transpose $W_1^t=W_2\cong Y(5,5,5,5,3)$
and $Y_4$ are saturated in $M_7$ by correspondingly 1, 1, 1, and 9
partial transversals of $M_7$.  ($W_1$ is depicted in
Fig.~\ref{Example 3}d).  On all other induced Young subdiagrams of
$M_7$, $(213)$ and $(123)$ are equally restrictive.  Therefore, the
Splitting Formulas give the remaining difference of $21$ in
(\ref{S_7}):
\begin{eqnarray*}
|S_7(2341)|-|S_7(3241)|&=&\sum_{W\in\{Y_5,W_1,W_2,Y_4\}}
\big(|S_{W}(123)|-|S_{W}(213)|\big)\cdot
|\bar{S}_{M_7\backslash W}(1)|\\
&=& (41-37)\cdot 1+(37-33)\cdot 1+(37-33)\cdot 1+(13-12)\cdot 9=21. 
\end{eqnarray*}
}
\end{ex}

\section{Avoidance on Young Diagrams with Extreme Critical
  Indices}\label{difference}

\subsection{The sets $|S_{Y_n}(\sigma)|$ and the Catalan numbers.}
For Young diagrams $Y$ with higher $i$-critical points, it is
interesting to find out by how much $(312)$ and $(321)$ are less
restrictive than $(321)$ and $(213)$, respectively. Below we answer
this for the diagram $Y=Y_n$ with {\it highest} critical index
$i=n-2$, and leave the general question to the reader.

\begin{figure}[h]
\labellist
\small\hair 2pt
\pinlabel ${M_n}$ at 107 586
\pinlabel ${\mathcal{B}_c}$ at 93 556
\pinlabel ${a}$ at 117 511
\pinlabel ${b}$ at 135 513
\pinlabel ${c}$ at 135 530
\pinlabel ${Y_n}$ at 224 586
\pinlabel ${\mathcal{A}_c\times\mathcal{B}_c}$ at 334 586
\pinlabel ${\mathcal{B}_c}$ at 336 556
\pinlabel ${c}$ at 307 530
\pinlabel ${=}$ at 167 538
\pinlabel ${+}$ at 282 542
\endlabellist
\begin{center}
\includegraphics[width=2.5in]{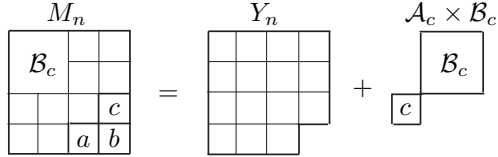}
\caption{$|S_n(213)|=|S_{Y_n}(213)|+|S_{n-2}(213)|$}
\label{SYn(213)}
\end{center} 
\end{figure}

\begin{prop} 
$|S_{Y_n}(213)|=c_n-c_{n-2}$ for $n\geq 2$.
\label{S_{Y_n}(213)}
\end{prop}
\noindent{\sc Proof:} This follows from Corollary~1 of the
 Row-Decomposition in Stankova-West \cite{Stankova-West}. Paraphrasing
 into the notation in the current paper, let $a,b,c$ be the three
 bottom right corner cells of $M_n$ as in Fig.~\ref{SYn(213)}. Then
 $M_n\backslash{\{b\}}=Y_n$. On the other hand, reducing $M_n$ along
 $a$ gives $M_n\big/_{\!\!a}=M_{n-1}$ whose right bottom cell is
 $c$. The minimal non-trivial $(213)$-decomposition of this $M_{n-1}$
 is obtained with respect to $c$: $(M_{n-1})_{213}(c)=\mathcal{A}_c
 \times \mathcal{B}_c=\{c\}\times M_{n-2}$. Thus, the {\it
 row-decomposition} formula for $S_{M_n}(213)$ in \cite{Stankova-West}
 reads: $|S_{M_n}(213)|=|S_{Y_n}(213)|+|S_{\mathcal{A}_c\times
 \mathcal{B}_c}(213)|$, from where
 $|S_{Y_n}(213)|=|S_{n}(213)|-|S_{n-2}(213)|=c_n-c_{n-2}$. \qed

\vspace*{3mm}

\smallskip
In the following, we keep the notation $b$ for the bottom right cell
of $M_n$, which is missing from $Y_n$. We shall enumerate
$S_{Y_n}(321)$ and $S_{Y_n}(312)$ by finding out how each differs as a
set from $S_n(321)$ and $S_n(312)$, respectively.

\begin{prop}
$|S_{Y_n}(321)|=c_n-1$ for $n\geq 2$.
\label{S_{Y_n}(321)}
\end{prop}
\noindent{\sc Proof:} Fix $T\in S_{Y_n}(321)$. Adding the cell $b$ to
$Y_n$ induces a transversal $T^{\prime}$ on $M_n$, which also avoids
$(321)$ on $M_n$. Indeed, if $(\alpha\beta\gamma)$ were a
$(321)$-pattern of $T^{\prime}$ in $M_n$, then $(\alpha\beta\gamma)$
lands on $\gamma$'s cell $d$. Since $d$ is dotted in $T^{\prime}$ on
$M_n$, it is also dotted in $T$ on $Y_n$, i.e.  $d\not = b$.  But then
$(\alpha\beta\gamma)$ is a $(321)$-pattern of $T$ landing on $d$ in
$Y_n$, a contradiction with $T\in S_{Y_n}(321)$.

Thus, we have a natural inclusion map
$\iota:S_{Y_n}(321)\hookrightarrow S_n(321)$. The reasoning above also
shows that the only transversals $T^{\prime}\in S_n(321)$ not hit by
$\iota$ are those with dotted $b$. However, in order to
avoid $(321)$ on $M_n$, a dot in $b$ implies that the
rest of $T^{\prime}$ is increasing, and there is only one such
transversal, namely, $T^{\prime}=(2,3,...,n,1)$. We conclude that
$S_n(321)=\iota(S_{Y_n}(321))\sqcup \{T^{\prime}\}$.
\[\Rightarrow\,\,|S_{Y_n}(321)|=|S_n(321)|-1=c_n-1\,\,
\text{for any}\,\,n\geq 2.\qed\]

\begin{prop}
$|S_{Y_n}(312)|=2c_n-3c_{n-1}$ for $n\geq 2$.
\label{S_{Y_n}(312)}
\end{prop}

\noindent
{\sc Proof:} As indicated above, we describe how $S_{Y_n}(312)$
differs as a set from $S_n(312)$.

\begin{figure}[h]
\labellist
\small\hair 2pt
\pinlabel ${b}$ at 460 513
\endlabellist
\begin{center}
\includegraphics[width=2.8in]{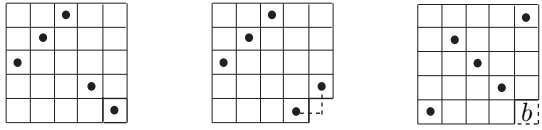}
\caption{Examples of the difference between $S_n(312)$ and $S_{Y_n}(312)$}
\label{appendix1}
\end{center} 
\end{figure}

\smallskip
On the one hand, $S_n(312)$ contains transversals of $M_n$ with a
dotted $b$. Since $b$ cannot participate in any $(312)$-pattern, we
can reduce $M_n$ along $b$ to obtain $M_{n-1}$ without any further
restrictions, and hence the number of transversals in question equals
$|S_{n-1}(312)|$. None of these transversals is in $S_{Y_n}(312)$
because $Y_n$ cannot have a dot in the missing $b$
(cf. Fig.~\ref{appendix1}a.) Thus, $|S_n(312)\backslash
S_{Y_n}(312)|=c_{n-1}$.

On the other hand, $S_{Y_n}(312)$ contains transversals of $Y_n$ for
which a $(312)$-subsequence lands outside $Y_n$ (on $b$). As we shall
see below in Lemma \ref{(312) landing on b}, the number of these
transversals is $c_n-2c_{n-1}$, and none of them is in $S_n(312)$
because of the $(312)$-pattern in $M_n$ (cf. Fig.~\ref{appendix1}b.)
Thus, $|S_{Y_n}(312)\backslash S_n(312)|=c_n-2c_{n-1}$.

All other transversals of $S_n(312)$ and $S_{Y_n}(312)$ are identical:
they don't have an element in $b$, and they don't have a
$(312)$-pattern landing on $b$ (cf. Fig.~\ref{appendix1}c.)
Summarizing,
\[|S_{Y_n}(312)|=|S_n(312)|-c_{n-1}+(c_n-2c_{n-1})=2c_n-3c_{n-1}\,\,
\text{for any}\,\,n\geq 2. \qed\]

Incidentally, we have shown the strict inequality
$|S_{Y_n}(312)|>|S_{Y_n}(321)|$ for $n\geq 5$ (proven in an indirect
way in Example~\ref{Example Y_n}). Indeed, from Propositions
\ref{S_{Y_n}(321)}-\ref{S_{Y_n}(312)}, for $n\geq 5$:
\begin{eqnarray*}
|S_{Y_n}(312)|-|S_{Y_n}(321)|&=&(2c_n-3c_{n-1})-(c_n-1) =
\frac{(n-5)n(2n-2)!}{(n+1)!(n-2)!}+1\geq 1. 
\end{eqnarray*}

\subsubsection{Claims in the Proof of Proposition~\ref{S_{Y_n}(312)}}
\label{Claims}

\begin{figure}[h]
\labellist
\small\hair 2pt
\pinlabel $\scriptstyle{A}$ at 54 547
\pinlabel $\scriptstyle{B}$ at 116 634
\pinlabel $\scriptstyle{C}$ at 208 700
\pinlabel $\scriptstyle{D_1}$ at 209 546
\pinlabel $\scriptstyle{D_2}$ at 209 633
\pinlabel $\scriptstyle{D_3}$ at 118 546
\pinlabel $\scriptstyle{D_4}$ at 118 700
\pinlabel $\scriptstyle{D_5}$ at 55 700
\pinlabel $\scriptstyle{D_6}$ at 55 633
\pinlabel $\scriptstyle{\beta}$ at 105 671
\pinlabel $\scriptstyle{\delta}$ at 81 616
\pinlabel $\scriptstyle{\alpha}$ at 171 512
\pinlabel $\scriptstyle{\gamma}$ at 244 583
\pinlabel $\scriptstyle{b}$ at 273 490
\pinlabel $\scriptstyle{T|_A}$ at 389 546
\pinlabel $\scriptstyle{T|_B}$ at 458 631
\pinlabel $\scriptstyle{T|_C}$ at 538 695
\pinlabel $\scriptstyle{\alpha}$ at 504 513
\pinlabel $\scriptstyle{\gamma}$ at 568 583
\endlabellist
\begin{center}
\includegraphics[width=3.3in]{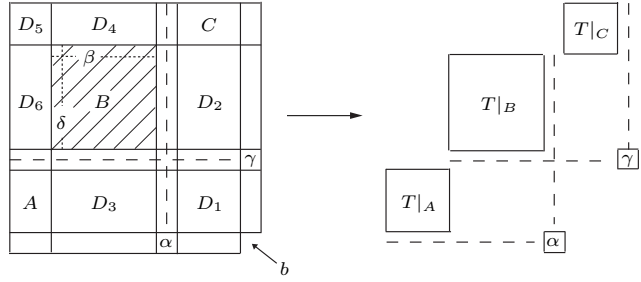}
\caption{Splitting of $T\in S_{Y_n}(312)$}
\label{9 Squares}
\end{center}
\end{figure}

\begin{lem}
The number of all transversals in $S_{Y_n}(312)$ with a
$(312)$-subsequence landing outside $Y_n$ (on $b$) is $c_n-2c_{n-1}$.
\label{(312) landing on b}
\end{lem}
\noindent{\sc Proof:} Let $T\in S_{Y_n}(312)$, and let $\alpha$ and
$\gamma$ denote the elements of $T$ in the bottom row and in the
rightmost column of $Y_n$, respectively. Because $Y_n$ misses
$b$, $\alpha\not = \gamma$ (cf. Fig.~\ref{9 Squares}a).

Suppose $T$ contains a $(312)$-subsequence which doesn't land in
$Y_n$, hence lands on $b$. Thus, for some $\beta\in T$,
$(\beta\alpha\gamma)\approx (312)$. Since $\beta$ is before $\alpha$
and above $\gamma$, without loss of generality, we can replace $\beta$
by the largest element of $T$ before $\alpha$; symbolically, $\beta:=
\max \{t\in T_{\bar{\alpha}}\}$. Symmetrically, let $\delta$ be the
leftmost element of $T$ higher than $\gamma$. It is possible that
$\delta=\beta$; if not, $(\delta\beta\alpha\gamma)\approx(3412)$ is a
subsequence of $T$ not landing in $Y_n$.
Let $B$ be the rectangle in $Y_n$ defined by $\beta$'s and $\gamma$'s
rows, and $\delta$'s and $\alpha$'s columns such that $B$ includes $\beta$
and $\delta$, but excludes $\alpha$ and $\gamma$. Let $A$ be the
rectangle below and to the left of $B$, excluding $\alpha$'s and
$\gamma$'s rows; and symmetrically, let $C$ be the rectangle to the
right and above $B$, excluding $\alpha$'s and $\gamma$'s columns.  

\begin{claim}
Except for $\alpha$ and $\gamma$, the transversal $T$ is concentrated
in rectangles $A$, $B$ and $C$.
\label{concentration of $T$}
\end{claim}
\noindent{\sc Proof:} $Y_n$ splits as a disjoint union of 9
rectangles, plus $\alpha$'s and $\gamma$'s rows and
columns. Figure~\ref{9 Squares}a depicts all these rectangles. The
definitions of $\beta$ and $\delta$ imply that rectangles $D_4$, $D_5$
and $D_6$ are empty. In order for the pair $(\beta\alpha)$ not to be
completed to a $(312)$-pattern in $Y_n$, rectangles $D_1$ and $D_2$
must also be empty. Symmetrically, in order for the pair
$(\delta\gamma)$ not to be completed to a $(312)$-pattern in $Y_n$,
rectangles $D_3$ and $D_1$ must be empty. Thus,
$T\backslash\{\alpha,\gamma\}$ is concentrated in $A$, $B$ and
$C$. \qed

\smallskip
We conclude that $T$ induces transversals on the rectangles $A$, $B$
and $C$, and since the latter are disjoint, they must be {\it
squares}. Thus, $T$ splits into an increasing sequence of 3
independent subtransversal $T|_A$, $T|_B$ and $T|_C$, with $\alpha$
inserted in the bottom row of $Y_n$ so that its column is between $B$
and $C$, and $\gamma$ is inserted in the rightmost column of $Y_n$ so
that its row is between $A$ and $B$. Finally, the assumption that $T$
contains a $(312)$-subsequence not landing inside $Y_n$ was translated
above in the existence $\beta\in B$, i.e. the square $B$ is of size at
least 1. Conversely,

\begin{claim}
If $T$ is a transversal of $Y_n$ satisfying the above description
(depicted also in Fig.~\ref{9 Squares}b), and such that the 3
subtransversals $T|_A$, $T|_B$ and $T|_C$ each avoid $(312)$ on the
respective squares $A$, $B$ and $C$, then the whole transversal $T$
avoids $(312)$ on $Y_n$, and has a $(312)$-subsequence not landing in
$Y_n$.
\label{Conversely}
\end{claim}

\noindent{\sc Proof:} Consider the reduction
$Y_n\big/_{{\!\{\alpha,\gamma\}}}=M_{n-2}$, along whose diagonal the
squares $A$, $B$ and $C$ are arranged (in increasing order). It is
evident that there can be no $(312)$-pattern in $M_{n-2}$ containing
elements from different squares. Since $T|_A$, $T|_B$ and $T|_C$ each
avoid $(312)$, any $(312)$-pattern in $T$ on $Y_n$ must contain
$\alpha$ and/or $\gamma$.  But $\alpha$ and $\gamma$ cannot
participate simultaneously in any pattern landing inside $Y_n$ because of
the missing cell $b$. Hence, only one of $\alpha$ and $\gamma$
can participate in a $(312)$-pattern in $Y_n$.

Since $\alpha$ can play only the role of ``1'', it can participate
only in a $(312)$-pattern of the form $(\xi \alpha \nu)$, where
$(\xi\nu)\!\!\searrow$, $\xi$ is before $\alpha$ and $\nu$ is after
$\alpha$. Yet, this arrangement is not possible since everything
before $\alpha$ is smaller than everything after $\alpha$: $A\oplus B
< C$, with the exception of $\gamma$, so no such pattern is possible.
``Transposing'' this argument, one concludes that $\gamma$
cannot participate in a $(312)$-pattern in $Y_n$ either.

Therefore, $T\in S_{Y_n}(312)$. Finally, since $B$ is of size at
least 1, let $\beta\in T|_B$. Then $(\beta\alpha\gamma)\approx (312)$
landing on $b$. \qed

\medskip
Claims~\ref{concentration of $T$}-\ref{Conversely} establish a 1-1
correspondence between the transversals $T\in S_{Y_n}(312)$ that do
{\bf not} induce transversals in $S_n(312)$ due to their
$(312)$-subsequence landing on $b$, and the diagrams in Figure~\ref{9
Squares}b. Therefore, each element of $S_{Y_n}(312)\backslash
S_n(312)$ is uniquely determined by the size of the squares $A$, $B$
and $C$, and the choice of (312)-avoiding transversals $T|_A, T|_B$
and $T|_C$. Below, the sum of sizes $|A|+|B|+|C|=n-2$ accounts for
$\alpha,\gamma\not \in A\cup B\cup C$.
\begin{eqnarray*}
S_{Y_n}(312)\backslash S_n(312)&\cong&
\!\!\bigsqcup_{\begin{tabular}{c}
\tiny{$|A|\!+\!|B|\!+\!|C|\!=\!n\!-\!2$}\\[-1.5mm]\tiny{$|B|\!\geq\!
1$}\end{tabular}}\!\! S_A(312)\times S_B(312)\times S_C(312)\\
\Rightarrow\,\,|S_{Y_n}(312)\backslash
S_n(312)|&=&\!\!\sum_{\begin{tabular}{c}
\tiny{$i\!+\!j\!+\!k\!=\!n\!-\!2$}\\[-1.5mm]\tiny{$j\!\geq\!
1$}\end{tabular}}\!\!c_ic_jc_k=c_n-2c_{n-1}. 
\end{eqnarray*}
The last equality was obtained using the well-known relation
${c_k=\sum_{l+m=k-1}c_lc_m}$ for the Catalan numbers. This completes
the proof of Lemma~\ref{(312) landing on b}. \qed

\subsection{The sets $|S_{St_n^3}(\tau)|$ and the Fibonacci Numbers}
In this subsection, we consider the other extreme situation of a {\it
non-decomposable} Young diagram $Y$: having {\it lowest} critical
indices $i=2$. This is $Y=St^3_n$ for $n\geq 4$, which is the {\it
smallest non-decomposable} Young diagram of size $n$. The last
description is also satisfied by the squares $M_n$ with $n\leq 3$, and
we set $St^3_n:=M_n$ for $n\leq 3$. This new notation and
Theorem~\ref{summary 312>321>213} imply that $(123)$ and $(312)$ are
equinumerant on $St^3_n$ for all $n\geq 1$, so that we can state the
following

\begin{prop}{\rm $|S_{St_n^3}(213)|=2^{n-3}(n+2)$ for $n\geq 2$ and 
\[|S_{St_n^3}(123)|=|S_{St_n^3}(312)|=f_{2n-1}=
\frac{1}{\sqrt{5}}\big(\psi^{2n-1}-\psi^{-(2n-1)}\big)\,\,\text{for}\,\,n\geq
1,\] where $f_n$ is the $n$-th Fibonacci number ($f_1=f_2=1$) and
$\psi=(1+\sqrt{5})/2$.}
\label{St^3}
\end{prop}

\begin{figure}[h]
\labellist
\small\hair 2pt
\pinlabel $\scriptstyle{a}$ at -63 583
\pinlabel $\scriptstyle{b}$ at -44 584
\pinlabel $\scriptstyle{c}$ at -27 583
\pinlabel $=$ at 65 633
\pinlabel $+$ at 227 633
\pinlabel ${2}$ at 97 633
\endlabellist
\begin{center}
\includegraphics[width=3.3in]{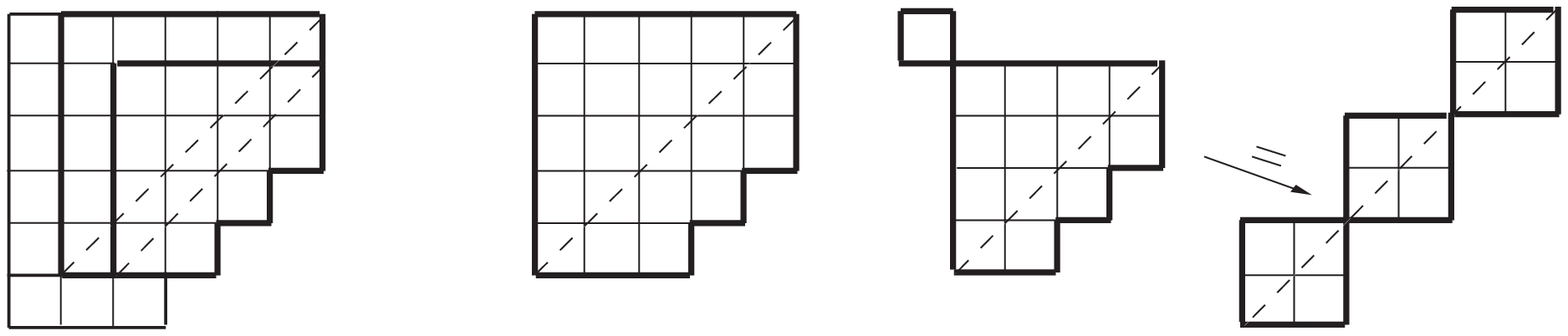}
\caption{Row-decomposition of $S_{St^3_n}(213)$}
\label{213 on St3}
\end{center}
\end{figure}

\noindent
{\sc Proof:} Let $a_n=|S_{St_n^3}(213)|$. From the row-decomposition
formula in \cite{Stankova-West} on $St^3_n$:
    \[|S_{St^3_n}(213)|=2\,|S_{St^3_{n-1}}(213)|+
    |S_{St^2_{n-2}}(213)|\,\,\text{(cf. Fig.~\ref{213 on St3})}.\]

\noindent Since $St^2_{n-2}$ $1$-decomposes as a product of $(n-3)$
squares $M_2$ (cf. Fig.~\ref{213 on St3}d), we have
$|S_{St^2_{n-2}}(213)|=2^{n-3}$. Thus, $a_n=2a_{n-1}+2^{n-3}$,
i.e. $a_n=2^{n-3}(n+2)$ for $n\geq 2$.

\begin{figure}[h]
\labellist
\small\hair 2pt
\pinlabel $\scriptstyle{a}$ at -63 582
\pinlabel $\scriptstyle{b}$ at -44 584
\pinlabel $\scriptstyle{c}$ at -27 582
\pinlabel ${2}$ at 79 637
\pinlabel $\scriptstyle{d}$ at 226 603
\pinlabel $\scriptstyle{e}$ at 261 601
\pinlabel ${=}$ at 58 637
\pinlabel ${+}$ at 200 637
\pinlabel ${+}$ at 421 575
\endlabellist
\begin{center}
\includegraphics[width=3.5in]{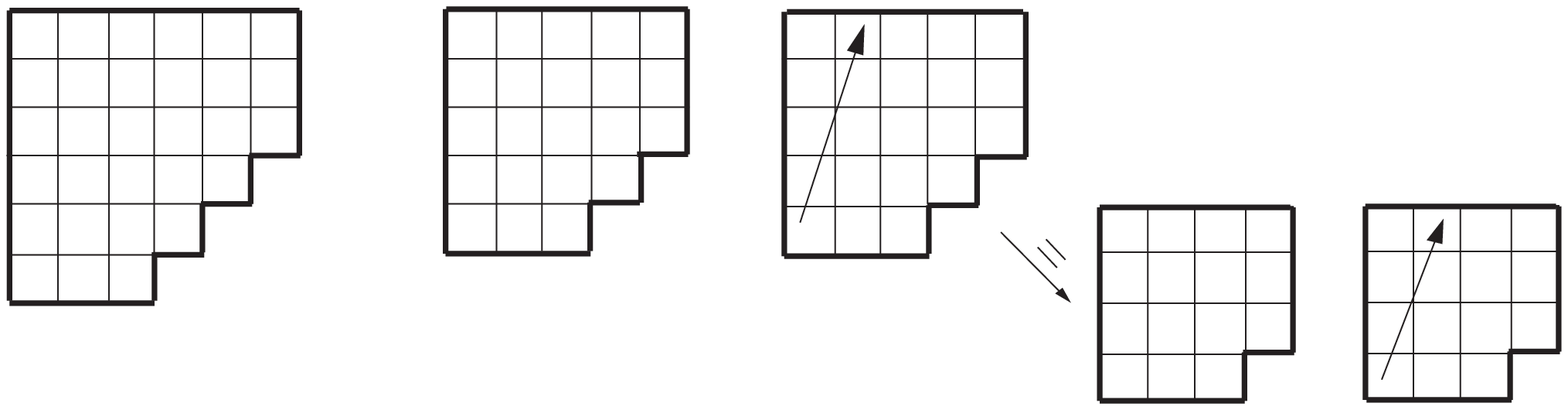}
\caption{Row-decomposition of $|S_{St^3_n}(312)|$}
\label{Golden Ratio}
\end{center}
\end{figure}

Consider now $b_n=|S_{St_n^3}(312)|$. Let $a$, $b$ and $c$ be the
bottom cells of $St^3_n$, as in Fig.~\ref{Golden Ratio}a.  Placing $1$
in $a$ or $c$ does not affect the $(312)$-avoidance in the reduction
$St_n^3/_{\!\{u\}}\cong St^3_{n-1}$ for $u=a$ or $c$
(cf. Fig.~\ref{Golden Ratio}b), and thus yields overall $2b_{n-1}$
transversals. However, placing $1$ in $b$ forces the elements in the
first two columns of the reduction $St_n^3/_{\!\{b\}}\cong St^3_{n-1}$
to form an increasing sequence (depicted by $\nearrow$ above cell $d$
in Fig.~\ref{Golden Ratio}c.)  In accordance with previous notation,
we denote the number of such $(312)$-avoiding transversals of
$St^3_{n-1}$ by $b^{\scriptscriptstyle{\nearrow}}_{n-1}$. Therefore,
$b_n=2b_{n-1}+b^{\scriptscriptstyle{\nearrow}}_{n-1}$.

To calculate $b^{\scriptscriptstyle{\nearrow}}_{n-1}$, note that $1$
can be placed only in the first cell $d$ or in the third cell $e$ of
the bottom row of $St^3_{n-1}$. The first case does not cause any
restrictions on the reduction $St_{n-1}^3/_{\!\{d\}}\cong St^3_{n-2}$
(cf. Fig.~\ref{Golden Ratio}d) and hence it produces $b_{n-2}$
transversals. Placing $1$ in $e$ reduces to
$b^{\scriptscriptstyle{\nearrow}}_{n-2}$ on
$St_{n-1}^3/_{\!\{e\}}\cong St^3_{n-2}$ (cf. Fig.~\ref{Golden
Ratio}e). Summarizing,
$b^{\scriptscriptstyle{\nearrow}}_{n-1}=b_{n-2}+
b^{\scriptscriptstyle{\nearrow}}_{n-2}$.  

Combining the two newly derived formulas, we obtain
$b_n=3b_{n-1}-b_{n-2}$, with $b_1=1$ and $b_2=2$. It is a standard
exercise to check that the odd-indexed terms in the Fibonacci sequence
satisfy the same recursive relation, and hence the desired
formula for $b_n$ involving the golden ratio $\psi$ follows. \qed

\subsection{Generalization of Stanley-Wilf limits.} 
Recall the Stanley-Wilf limits $L(\tau)=\lim_{n\rightarrow \infty}
\sqrt[n]{|S_n(\tau)|}$ for any $\tau\in S_k$. From works of Regev
\cite{Regev} and B\'{o}na \cite{Bon4}, it follows that
$L(J_k)=(k-1)^2$ and $L(213|J_k)=(k-1+\sqrt{8})^2$ for $k\geq 1$. 

From the viewpoint of the current paper,
Corollary~\ref{Wilf-ordering-corollary} has established in particular
the strict inequalities $|S_n(213|J_k)|< |S_n(123|J_k)| <
|S_n(312|J_k)|$ for any $k\geq 1$ and $n\geq 2k+5$. Hence the
Stanley-Wilf limits follow suit for $k\geq 1$:
\[L(213|J_k)=(k-1+\sqrt{8})^2 <
L(123|J_k)=L(321|J_k)=L(J_{k+3})=(k+2)^2 \leq L(312|J_k).\] It is
still an open question whether $L(312|1)<9$, and in order to complete
the above picture, it would certainly be nice to find the exact value
of $L(312|J_k)$. Conceivably, the Splitting Formula for $(312|J_k)$
from Subsection~\ref{splitting formula subsection} and other
observations in this paper might be helpful towards calculating
$L(312|J_k)$. 

With the methods so far, all known $L(\tau)$ belong to
$\mathbb{Z}[\sqrt{2}]$. However, if we generalize the definition of
Stanley-Wilf limits from the square matrices $M_n$ to using any
(proper) Young diagrams $Y$ of size $n$, we can obtain presumably a
much greater variety of limits. To this end, consider the set
$\mathcal{Y}=\cup_{n=0}^{\infty} \mathcal{Y}^n$ of all proper Young
diagrams, {\it graded} by the size $n$ of the diagrams.  Let
$\vec{Y}=\{Y^n\}$ be a sequence of (proper) Young diagrams, one per
each graded piece of $\mathcal{Y}$; we can think of $\vec{Y}$ as a
{\it path} in $\mathcal{Y}$. Define the {\it generalized Stanley-Wilf
limit} of $\tau\in S_k$ along the path $\vec{Y}$ as
\[L_{\vec{Y}}(\tau)=\lim_{n\rightarrow
  \infty}\sqrt[n]{|S_{Y^n}(\tau)|}.\] Except for the case
$\vec{Y}=\{M_n\}$ where the limits $L_{\vec{Y}}(\tau)=L(\tau)$ are
guaranteed by Stanley-Wilf Theorem, for all other paths in
$\mathcal{Y}$ the existence of $L_{\vec{Y}}(\tau)$ must be verified.

A worthwhile consequence of Proposition~\ref{St^3} is the following
\begin{cor}
For $\vec{Y}=\{St^3_n\}$,
$L_{\vec{Y}}(321)=L_{\vec{Y}}(312)=\psi^2={\frac{3+\sqrt{5}}{2}}\cdot$
\end{cor}
Two natural questions arise: for which pairs $(\vec{Y},\tau)$ do the
limits $L_{\vec{Y}}(\tau)$ exist, and what is the algebraic closure
$\overline{\mathcal{L}}$ of the set of generalized limits
$\mathcal{L}=\{L_{\vec{Y}}(\tau)\}$. As of now, we have shown that
$\overline{\mathcal{L}}\supset\mathbb{Q}(\sqrt{2},\sqrt{5})$; but are
there any other irrational or transcendental generalized
Stanley-Wilf limits $L_{\vec{Y}}(\tau)$? We leave these questions to
the reader for further study.

\section* {Acknowledgments}

The author would like to thank Miklos B\'{o}na (University of Florida)
for supplying a number of useful references and discussing his
and related works in relation to the present paper; David Moews
(Center for Communications Research, San Diego) for writing a computer
program used in this project; and Paulo de Souza (UC Berkeley) for his
help in implementing the necessary computer software.


\begin{thebibliography}{[14]}

\bibitem{Arratia} R. Arratia, On the Stanley-Wilf Conjecture for the
Number of Permutations Avoiding a Given Pattern, {\it Electronic
J. Combin.}, {\bf 6} (1999), no. 1, N1.

\bibitem{BW} E. Babson, J. West,  The permutations
$123p_4...p_t$ and $321p_4...p_t$ are Wilf-equivalent, 
Graphs Comb 16 (2000) 4, 373-380.

\bibitem{BWX} J. Backelin, J. West, G. Xin,  Wilf-equivalence
for singleton classes, {\it Proceedings of the 13th Conference on
  Formal Power Series and Algebraic Combinatorics}, Tempe, AZ, 2001.

\bibitem{Bon1} M. B\'{o}na,  Permutations avoiding certain patterns.
The case of length 4 and some generalizations, Disc. Math. 175 (1997)
55-67.

\bibitem{Bon2} M. B\'{o}na,  The Solution of a Conjecture of Wilf
and Stanley for all layered patterns, J. Combinatorial Theory,
Series A, 85 (1999) 96-104.

\bibitem{Bon3} M. B\'{o}na, Combinatorics of Permutations, Chapman \&
Hall/CRC, 2004, 135-159.

\bibitem{Bon4} M. B\'{o}na, The Limit of a Stanley-Wilf sequence is
not always rational, and layered patterns beat monotone patterns,
J. Combin. Theory Ser. A 110 (2) (2005), 223-235.

\bibitem{Kn1} D. Knuth,  {\it The Art of Computer Programming},
Vol.3, Addison-Wesley, Reading, MA, 1973.

\bibitem{Kn2} D. Knuth, Permutations, matrices, and generalized Young
tableaux, Pacific J. of Mathematics 34 (1970) 709-727.

\bibitem{Marcus} A. Marcus and J. Tardos, Excluded Permutation
Matrices and the Stanley-Wilf conjecture, J. Combin. Theory
Ser. A 107 (1) (2004), 153-160.


\bibitem{Regev} A. Regev, Asymptotic values for degrees associated
with strips of Young diagrams, {\it Advances in Mathematics}, {\bf 41}
(1981), 115-136.

\bibitem{Ri} D. Richards,  Ballot sequences and restricted
permutations, Ars Combinatoria 25 (1988) 83-86.

\bibitem{Ro} D. Rotem,  On correspondence between binary trees
and a certain type of permutation, Information Processing Letters 4
(1975), 58-61.

\bibitem{SS} R. Simion and F. Schmidt,  Restricted permutations,  
Europ. J. Combinatorics 6 (1985) 383-406.

\bibitem{St1} Z. Stankova,  Forbidden subsequences,
Disc. Math. 132 (1994) 291-316.

\bibitem{St2} Z. Stankova,  Classification of forbidden
subsequences of length 4, Europ. J. Combinatorics (1996) 17, 501-517.

\bibitem{Stankova-West} Z. Stankova and J. West, A new class of
Wilf-equivalent permutations, {\it J. Algebraic Combin.}, {\bf 15}
(2002), no. 3, 271-290.

\bibitem{We0} J. West,  Generating trees and the Catalan and
  Schr\"{o}der numbers, Disc. Math. 146 (1995) 247-262.

\bibitem{We1} J. West,  Generating trees and forbidden subsequences,
Disc. Math. 157 (1996) 363-374.
\end{thebibliography}
\end{document}